\documentclass[12pt]{article}
\usepackage{amsthm,amsmath,amssymb,amscd,verbatim}
\usepackage{graphics,pifont}
\usepackage{amssymb}
\usepackage{graphicx}
\usepackage{caption}
\date{}
\newcommand{\f}{\frac}

\begin{document}
\title{An Interesting Application of\\ the Intermediate Value Theorem: \\A Simple Proof of Sharkovsky's Theorem and \\ the Towers of Periodic Points}
\author{Bau-Sen Du\footnote {Bau-Sen Du is a retired research fellow at the Institute of Mathematics, Academia Sinica, Taiwan} \\ [.5cm] 
dubs@gate.sinica.edu.tw \\}
\maketitle
\begin{abstract}
This note is intended primarily for college calculus students right after the introduction of the Intermediate Value Theorem, to show them how the Intermediate Value Theorem is used repeatedly and straightforwardly to prove the celebrated Sharkovsky's theorem on the periods of coexistent periodic orbits of continuous maps on an interval.  Furthermore, if the maps have a periodic orbit $P$ of odd period $\ge 3$, then we find more periodic points, in the appendix, which constitute what we call {\it the towers of periodic points associated with $P$}.  These towers of periodic points have infinitely many layers.  In this note, {\large No knowledge of Dynamical Systems Theory is required}. 

In this revision, we add the result (Proposition 1 on page 9) that the periodic orbits of continuous unimodal maps on the interval $[0, 1]$ of least periods $\ge 2$ are {\it nested} in the sense that if $P$ and $Q$ are periodic orbits of a continuous unimodal map on $[0, 1]$ of least periods $\ge 2$ and if $\max P < \max Q$, then $[\min P, \max P] \subset [\min Q, \max Q]$.

\bigskip
\noindent{{\bf AMS Subject Classification}: 37E05, 37E15, 26A18}
\end{abstract}

Throughout this note, $I$ is a compact interval, and $f : I \rightarrow I$ is a continuous map.  For each integer $n \ge 1$, let $f^n$ be defined by: $f^1 = f$ and $f^n = f \circ f^{n-1}$ when $n \ge 2$.  For $x_0$ in $I$, we call the set $\big\{ x_0, f(x_0), f^2(x_0), \cdots \big\}$ the orbit of $x_0$ with respect to $f$ and call $x_0$ a periodic point of $f$ with least period $m$ or a period-$m$ point of $f$ if $f^m(x_0) = x_0$ and $f^i(x_0) \ne x_0$ when $0 < i < m$.  If $f(x_0) = x_0$, then we call $x_0$ a fixed point of $f$.  

In this note, we demonstrate how the Intermediate Value Theorem is applied repeatedly and straightforwardly to prove (1) of the following celebrated Sharkovsky's theorem on the periods of coexistent periodic orbits of continuous maps on the interval $I$:   

\noindent
{\bf Theorem (Sharkovsky{\bf{\cite{sh1}}})}  {\it The Sharkovsky's ordering of the natural numbers is as follows: $$3 \prec 5 \prec 7 \prec 9 \prec \cdots \prec 2 \cdot 3 \prec 2 \cdot 5 \prec 2 \cdot 7 \prec 2 \cdot 9 \prec \cdots \prec 2^2 \cdot 3 \prec 2^2 \cdot 5 \prec 2^2 \cdot 7 \prec 2^2 \cdot 9 \prec \cdots$$ $$\prec 2^i \cdot 3 \prec 2^i \cdot 5 \prec 2^i \cdot 7 \prec 2^i \cdot 9 \prec \cdots \prec 2^{n+1} \prec 2^n \prec \cdots \prec 2^3 \prec 2^2 \prec 2  \prec 1.$$   Then the following three statements hold:
\begin{itemize}
\item[\rm{(1)}] 
If $f$ has a period-$m$ point and if $m \prec n$, then $f$ also has a period-$n$ point.

\item[\rm{(2)}]
For each positive integer $n$ there exists a continuous map from $I$ into itself that has a period-$n$ point but has no period-$m$ point for any $m$ with $m \prec n$.

\item[\rm{(3)}]
There exists a continuous map from $I$ into itself that has a period-$(2^i)$ point for $i = 0, 1, 2, \ldots$ but has no periodic point of any other period.
\end{itemize}}

To make this note self-contained, we include the following two well-known results.

\noindent
{\bf Lemma 1.}
If $f^n(x_0) = x_0$, then the least period of $x_0$ with respect to $f$ divides $n$.

\noindent
{\it Proof.}
Let $m$ denote the least period of $x_0$ with respect to $f$ and write $n = km + r$ with $0 \le r < m$.  Then $x_0 = f^n(x_0) = f^{km+r}(x_0) = f^r\big(f^{km}(x_0)\big) = f^r(x_0)$.  Since $m$ is the smallest positive integer such that $f^m(x_0) = x_0$, we must have $r = 0$.  Therefore, $m$ divides $n$.  

\noindent
{\bf Lemma 2.}
{\it Let $i, j, k, m, n$ and $s$ be positive integers.  Then the following statements hold:
\begin{itemize}
\item[\rm{(i)}]
If $x_0$ is a periodic point of $f$ with least period $m$, then it is a periodic point of $f^n$ with least period $m/(m, n)$, where $(m, n)$ is the greatest common divisor of $m$ and $n$. 

\item[\rm{(ii)}]
If $x_0$ is a periodic point of $f^n$ with least period $k$, then it is a periodic point of $f$ with least period $kn/s$, where $s$ divides $n$ and is relatively prime to $k$.  In particular, if $f^{2^{i-1}}$ has a period-$(2 \cdot j)$ point for some $i \ge 2$ and $j \ge 1$, then $f$ has a period-$(2^i \cdot j)$ point.  
\end{itemize}}

\noindent
{\it Proof.} (i) Let $x_0$ be a period-$t$ point of $f^n$.  Then $m$ divides $nt$ since $x_0 = (f^n)^t(x_0) = f^{nt}(x_0)$.  So, $\f m{(m, n)}$ divides $\f n{(m, n)} \cdot t$.  Since $\f m{(m, n)}$ and $\f n{(m, n)}$ are coprime, $\f m{(m, n)}$ divides $t$.  Furthermore, $(f^n)^{(m/(m, n))}(x_0) = (f^m)^{(n/(m, n))}(x_0) = x_0$.  Thus, $t$ divides $\f m{(m, n)}$.  This shows that $t = \f m{(m, n)}$.

(ii) Since $x_0 = (f^n)^k(x_0) = f^{kn}(x_0)$, the least peirod of $x_0$ under $f$ is $\f {kn}s$ for some positive integer $s$.  By (i), $(\f {kn}s)/\big((\f {kn}s), n\big) = k$.  So, $\f ns = \big((\f ns)k, n\big)$\,\, (which is an integer) \,\, $= \big((\f ns)k, (\f ns)s\big) = (\f ns)(k, s)$.  This shows that $s$ divides $n$ and $(s, k) = 1$.

\vspace{.1in}

The version of the Intermediate Value Theorem we are going to use is the following:

\noindent
{\bf Intermediate Value Theorem} 
{\it Let $g : [a, b] \longrightarrow \mathcal R$ be a continuous real-valued map such that $g(a)g(b) < 0$.  Then there exists a point $a < c < b$ such that $g(c) = 0$.}

\vspace{.1in}

Following {\bf{\cite{str}}}, we first prove the following three statements: 
\begin{itemize}
\item[(a)] if $f$ has a period-$m$ point with $m \ge 2$, then $f$ has a period-2 point and a fixed point; 

\item[(b)] if $f$ has a period-$m$ point with $m \ge 3$ and odd, then $f$ has a period-$(m+2)$ point; and 

\item[(c)] if $f$ has a period-$m$ point with $m \ge 3$ and odd, then $f$ has periodic points of all even periods.  
\end{itemize}

Let $P$ be a period-$m$ orbit of $f$ with $m \ge 2$ and let $e = f^{m-1}(\min P)$.  
{\large 
$$
\text{We choose the points} \,\,\, \min P \,\,\, \text{and} \,\,\, f^{m-1}(\min P) \,\,\, \text{in the orbit} \,\,\, P \,\,\, \text{to start with}.
$$}
They seem to be the right choices because they make the subsequent proofs of (a), (b) and (c) become straightforward.  See also Appendix.

{\bf We now prove (a), (b) and (c)}

Since $\min P$ is a period-$m$ point of $f$ with $m \ge 2$, we have $f(e) = \min P < e$.  So, $f(e) \ne e$.  
$$
\text{If} \,\,\, f(x) < e = f^{m-1}(\min P) \,\,\, \text{for all} \,\,\, \min P \le x < e = f^{m-1}(\min P),
$$ 
then, since $\min P \le f^i(\min P)$ for all $i \ge 1$, we have 
$$
(\min P \le) \,\,\, f^i(\min P) < e \,\,\, \text{for all} \,\,\, i \ge 1,
$$
contradicting the fact that $f^{m-1}(\min P) = e$.  Since $f(e) = \min P$, we have $f\big([\min P, e]\big) \supset [\min P, e] \supset \{ e \}$.  Let $v$ be any point in $[\min P, e)$ such that $f(v) = e = f^{m-1}(\min P)$.  Then 
$$
f(v) - v = e - v > 0 \,\,\, \text{and} \,\,\, f(e) - e = \min P - e < 0
$$ 
which imply that $f$ has a fixed point $z$ in $(v, e)$.  

Now suppose $m \ge 3$.  Then $\min P < v$.  Since 
$$
f^2(\min P) - \min P > 0 \,\,\, \text{and} \,\,\, f^2(v) - v = \min P - v < 0,
$$ 
$$\text{the point} \,\,\, y = \max \big\{ \min P \le x \le v : f^2(x) = x \big\} \,\,\, \text{exists}.
$$
If $f(y) = y$, then since 
$$
f(y) - z = y - z < 0 \,\,\, \text{and} \,\,\, f(v) - z = f^{m-1}(\min P) - z > 0,
$$
there is a point $\tilde x$ in $(y, v)$ such that $f(\tilde x) = z$.  Consequently, since 
$$
f^2(\tilde x) - \tilde x = z - \tilde x > 0 \,\,\, \text{and} \,\,\, f^2(v) - v = \min P - v < 0,
$$ 
there is a point $\tilde y$ in $(\tilde x, v) \, (\subset (y, v))$ such that $f^2(\tilde y) = \tilde y$ which contradicts the maximality of $y$ in $[\min P, v]$.  Therefore, $f(y) \ne y$ and $y$ is a period-2 point of $f$.  This confirms (a).

\vspace{.1in}

For the proofs of (b) and (c), assume that $m \ge 3$ is odd.  Let 
$$
z_0 = \min \big\{ v \le x \le z : f^2(x) = x \big\}.
$$  
Then since $f^2(v) = \min P < v$, we have 
{\large 
$$
\qquad\quad\qquad\qquad\,\,\, f^2(x) < x < z_0 \,\,\, \text{for all} \,\,\, v \le x < z_0. \qquad\qquad\qquad\qquad\qquad (*)
$$}
If $f^2(x) < z_0$ for all $\min P \le x < v$, then $f^2(x) < z_0$ for all $\min P \le x < z_0$.  Consequently, 
$$
(\min P \le) \,\,\, f^{2i}(\min P) < z_0 \,\,\, \text{for all} \,\,\, i \ge 1,
$$
contradicting the fact that $(f^2)^{(m-1)/2}(\min P) = e > z \ge z_0$.  So, $f^2(\hat x) \ge z_0$ for some $\min P \le \hat x < v$.  Since $f^2(\min P)$ and $z_0$ are perodic points of $f$ with different periods, $f^2(\min P) \ne z_0$ and, since $$f^2(\hat x) - z_0 \ge 0 \,\,\, \text{and} \,\,\, f^2(v) - z_0 = \min P - z_0 < 0,$$ the point 
$$
d = \max \big\{ \min P \le x \le v : f^2(x) = z_0 \big\} \,\,\, (> \min P \,\,\, \text{since} \,\,\, f^2(\min P) \ne z_0)
$$ 
exists and, since $f^2(v) = \min P < z_0$, we have $f^2(x) < z_0$ for all $d < x < v$.  This, combined with the above (*), implies that 
{\large 
$$
\qquad\qquad\qquad\qquad\,\,\, f^2(x) < z_0 \,\,\, \text{for all} \,\,\, d < x < z_0. \qquad\qquad\qquad\qquad\qquad\,\,\, (**)
$$}  
If $f(\bar x) = z$ for some $d < \bar x < z_0$, then $f^2(\bar x) = z \ge z_0$ which contradicts the fact that $f^2(x) < z_0$ for all $d < x < z_0$.  Since $f(v) = e = f^{m-1}(\min P) > z$, we have 
{\large 
$$
\quad\qquad\qquad f(x) > z \,\,\, (\ge z_0 > f^2(x)) \,\,\, \text{for all} \,\,\, d < x < z_0, \qquad\qquad\qquad\quad\,\, (***)
$$}
Consequently, $f(z_0) \ge z \ge z_0$ and $f(d) \ge z \ge z_0$.  Now since 
$$
f^2(d) - d = z_0 - d > 0 \,\,\, \text{and} \,\,\, f^2(v) - d = \min P - d < 0,
$$ 
there is a point $u$ in $(d, v)$ such that $f^2(u) = d$.  Let 
{\large 
$$
u_1 = \min \big\{ d \le x \le v : f^2(x) = d \big\}.
$$}
In the following, 

{\large 
for the proof of (b), we consider the interval $[u_1, v]$ (instead of taking the period-2 point $y$ obtained in the proof of (a) and considering the interval $[y, v]$ as in {\bf\cite{du2}}) and show that, for each $n \ge 1$, the point 
$$
p_{m+2n} = \min \big\{ u_1 \le x \le v : f^{m+2n}(x) = x \big\}
$$ 
exists and is a period-$(m+2n)$ point of $f$ while 

for the proof of (c), we consider the interval $[d, u_1]$ and show that, for each $n \ge 1$, the point 
$$c_{2n} = \min \big\{ d \le x \le u_1 : f^{2n}(x) = x \big\}
$$ 
exists and is a period-$(2n)$ point of $f$.}    

\vspace{.1in}

Note that $f^3(u_1) = f\big(f^2(u_1)\big) = f(d) \ge z \ge z_0$ (by (***)) and, for each $j \ge 5$ and odd, $f^j(u_1) = f^{j-4}\big(f^4(u_1)\big) = f^{j-4}(z_0) = f(z_0) \ge z \ge z_0$.  Therefore, 
$$
f^n(u_1) \ge z \ge z_0 \,\,\, \text{for all {\it odd}} \,\,\, n \ge 3.
$$
Now since $m \ge 3$ is odd, we have 
$$
f^{m+2}(u_1) - u_1 = f(z_0) - u_1 \ge z_0 - u_1 > 0 \,\, \text{and} \,\, f^{m+2}(v) - v = f^m\big(f^2(v)\big) - v = \min P - v < 0.
$$  
Therefore, the point 
$$
p_{m+2} = \min \big\{ u_1 \le x \le v : f^{m+2}(x) = x \big\} \,\,\, \text{exists and} \,\,\, f(p_{m+2}) > z \ge z_0 > p_{m+2}.
$$  
Let $k$ denote the least period of $p_{m+2}$ with respect to $f$.  Then, $k > 1$ and, by Lemma 1, $k$ is odd.  Now we want to show that $k = m+2$.  Suppose $(3 \le) \,\,\, k < m+2$.  Then since 
$$
f^k(u_1) - v  \ge z_0 - v > 0 \,\,\, \text{and} \,\,\, f^k(p_{m+2}) - v = p_{m+2} - v < 0,
$$
there is a point $v_k$ in $(u_1, p_{m+2})$ such that $f^k(v_k) = v$.  Consequently, since 
$$
f^{k+2}(u_1) - u_1 = f(z_0) - u_1 \ge z_0 - u_1 > 0 \,\,\, \text{and} \,\,\, f^{k+2}(v_k) - v_k = f^2(v) - v_k = \min P - v_k < 0,
$$
there is a point $w_{k+2}$ in $(u_1, v_k)$ such that $f^{k+2}(w_{k+2}) = w_{k+2}$.  Similarly, since 
$$
f^{k+2}(u_1) - v  \ge z_0 - v > 0 \,\,\, \text{and} \,\,\, f^{k+2}(w_{k+2}) - v = w_{k+2} - v < 0,
$$
there is a point $v_{k+2}$ in $(u_1, w_{k+2})$ such that $f^{k+2}(v_{k+2}) = v$.  Furthermore, since 
$$
f^{k+4}(u_1) - u_1 \ge z_0 - u_1 > 0 \,\,\, \text{and} \,\,\, f^{k+4}(v_{k+2}) - v_{k+2} = f^2(v) - v_{k+2} = \min P - v_{k+2} < 0,
$$
there is a point $w_{k+4}$ in $(u_21, v_{k+2})$ such that $f^{k+4}(w_{k+4}) = w_{k+4}$.  Inductively, there exist points 
$$
u_1 < \cdots < w_{m+2} < w_m < w_{m-2} < \cdots < w_{k+4} < w_{k+2} < p_{m+2} < v
$$ 
such that $f^{k+2i}(w_{k+2i}) = w_{k+2i}$ for all $i \ge 1$.  In particular, $f^{m+2}(w_{m+2}) = w_{m+2}$ and $u_1 < w_{m+2} < p_{m+2}$, contradicting the fact that $p_{m+2}$ is the {\it smallest} point in $(u_1, v)$ which satisfies $f^{m+2}(x) = x$.  Therefore, $k = m+2$.  This establishes (b).  

Note that similar arguments as above show that the point 
$$
p_{m+4} = \min \big\{ u_1 \le x \le p_{m+2} : f^{m+4}(x) = x \big\}
$$
exists and is a period-$(m+4)$ point of $f$.  Inductively, we obtain a sequence of points 
$$
u_1 < \cdots < p_{m+2i} < \cdots < p_{m+4} < p_{m+2} < v < z_0
$$ 
such that each $p_{m+2i} = \min \{ d \le x \le v : f^{m+2i}(x) = x \}, i \ge 1$, is a period-$(m+2i)$ point of $f$.      

\vspace{.1in}

We now prove (c).  Note that $f^{2i}(d) = z_0$ for all $i \ge 1$ and recall that we have the property 
{\large 
$$
\qquad\qquad\qquad\qquad\,\,\, f^2(x) < z_0 \,\,\, \text{for all} \,\,\, d < x < z_0. \qquad\qquad\qquad\qquad\qquad\,\,\, (**)
$$}
Since 
$$
f^2(d) = z_0 \,\,\, \text{and} \,\,\, u_1 = \min \big\{ d \le x \le v : f^2(x) = d \big\} \,\, (\text{and so}, \, f^2(u_1) = d),
$$ 
$$
\text{we have} \,\,\, d < f^2(x) < z_0 \,\,\, \text{on} \,\,\, (d, u_1) \,\,\, \text{and so, by (**)}, \,\, f^4(x) = (f^2)^2(x) < z_0 \,\,\, \text{on} \,\,\, (d, u_1).
$$
Let $c_2 = \min \big\{ d \le x \le u_1 : f^2(x) = x \big\} = \min \big\{ d \le x \le v : f^2(x) = x \big\}$.  Then $f^4\big([d, c_2]\big) = f^2\big(f^2([d, c_2])\big) \supset f^2\big([c_2, z_0]\big) \supset f^2\big([u_1, z_0]\big) \supset \{ d \}$ (since $f^2(u_1) = d)$.  Since 
$$
f^4(d) = z_0,\,\, f^4(x) < z_0 \,\,\ \text{on} \,\,\, (d, u_1) \,\,\, \text{and} \,\,\, f^4\big([d, u_1]\big) \supset f^4\big([d, c_2]\big) \supset \{ d \},
$$ 
$$
\text{the point} \,\,\, u_2 = \min \big\{ d \le x \le c_2 : (f^2)^2(x) = d \big\} \,\, (\text{and so}, \, f^4(u_2) = d) \,\,\, \text{exists and}
$$ 
$$\text{we have} \,\,\, d < (f^2)^2(x) < z_0 \,\,\, \text{on} \,\,\, (d, u_2) \,\,\, \text{and so, by (**)}, \,\, f^6(x) = (f^2)^3(x) < z_0 \,\,\, \text{on} \,\,\, (d, u_2).
$$
Let $c_4 = \min \big\{ d \le x \le u_2 : f^4(x) = x \big\} = \min \big\{ d \le x \le v : f^4(x) = x \big\}$.  Then $f^6\big([d, c_4]\big) = f^2\big(f^4([d, c_4])\big) \supset f^2\big([c_4, z_0]\big) \supset f^2\big([u_1, z_0]\big) \supset \{ d \}$.  Since 
$$
f^6(d) = z_0, \,\, f^6(x) < z_0 \,\,\ \text{on} \,\,\, (d, u_2) \,\,\, \text{and} \,\,\, f^6\big([d, u_2]\big) \supset f^4\big([d, c_4]\big) \supset \{ d \}
$$ 
$$\text{the point} \,\,\, u_3 = \min \big\{ d \le x \le c_4 : (f^2)^3(x) = d \big\} \,\,\, \text{exists and}
$$ 
$$\text{we have} \,\,\, d < (f^2)^3(x) < z_0  \,\,\, \text{on} \,\,\, (d, u_3) \,\,\, \text{and so, by (**)}, \,\, f^8(x) = (f^2)^4(x) < z_0 \,\,\, \text{on} \,\,\, (d, u_3).
$$
Let $c_6 = \min \big\{ d \le x \le u_3 : f^6(x) = x \big\} = \min \big\{ d \le x \le v : f^6(x) = x \big\}$.  Proceeding in this manner indefinitely, we obtain points 
$$
d < \cdots < c_{2n} < u_{n} < \cdots < c_4 < u_2 < c_2 < u_1 < v < z_0
$$ 
such that $d < (f^2)^n(x) < z_0$ on $(d, u_{n})$ and $(f^2)^n(c_{2n}) = c_{2n}$ for all $n \ge 1$.  Since $f(x) > z \ge z_0$ on $(d, z_0)$, we have 
{\large 
$$
d < f^i(c_{2n}) < z_0 \le z < f^j(c_{2n}) \,\,\, \text{for all even} \,\,\, i \,\,\, \text{and all odd} \,\,\, j \,\,\, \text{in} \,\,\, [0, 2n].
$$}
Consequently, each $c_{2n}$ has even period $\le 2n$ with respect to $f$.  Since, for each $1 \le k < n$, $c_{2k}$ is the smallest point in $[d, u_k]$ $(\supset [d, u_n])$ which satisfies $f^{2k}(x) = x$ and $c_{2n} < c_{2k}$, we obtain that $f^{2k}(c_{2n}) \ne c_{2n}$ for all $1 \le k < n$.  So, $c_{2n}$ is a period-$(2n)$ point of $f$.  This proves (c).

Furthermore, since we have shown that, for each $1 \le i \le n$, $d < (f^2)^i(x) < z_0$ for all $x$ in $(d, u_{i})$, we obtain that 
{\large 
$$
\qquad\qquad\,\,\,\, d < f^{2i}(x) < z_0 \,\,\, \text{for all} \,\,\, x \,\,\, \text{in} \,\,\, (d, u_{n}) \,\,\, \text{and all} \,\,\, 1 \le i \le n.\qquad\qquad \hspace{.3mm}(\dagger)
$$}
This fact will be used in the Appendix.

In the following, we show that, for each $n \ge 1$, the interval $[d, u_{n}]$ contains no periodic points of $f$ of even periods $< 2n$. 
\begin{multline*}
$$
\text{The arguments we used here will also be used repeatedly in the Appendix to show} \\ 
\text{the non-existence of periodic points of} \,\,\, f \,\,\, \text{of certain even periods in some intervals}.
$$
\end{multline*}
\indent Suppose, for some $n \ge 2$, the interval $[d, u_{n}]$ contained a periodic point $p$ of $f$ whose least period $2k_n$ is even and $< 2n$.  So, $2(k_n+1) \le 2n$ and $u_{n} \le u_{k_n+1}$.  Then since $f^{2i}(d) = z_0$ for all $i \ge 1$, we have $f^{2k_n}\big([d, p]\big) \supset [p, z_0] \supset \{ v \}$ and so, there is a point $\nu$ in $[d, p]$ such that $f^{2k_n}(\nu) = v$.  Therefore, $f^{2k_n+2}\big([d, \nu]\big) \supset [\min P, z_0] \supset \{ d \}$.  This implies the existence of a point $u_{k_n+1}^*$ in $[d, \nu]$ such that $f^{2k_n+2}(u_{k_n+1}^*) = d$.  Since $d < u_{k_n+1}^* < \nu < p < u_{n} \le u_{k_n+1} < v$, this contradicts the fact that $u_{k_n+1}$ is the smallest point in $[d, v]$ that satisfies $f^{2k_n+2}(x) = d$.  Thus, we have shown that, for each $n \ge 1$, the interval $[d, u_n]$ contains no periodic point of $f$ whose least period is {\it even} and $< 2n$.  This fact will be used below in the Appendix.  In particular, for each $n \ge 1$, the point $c_{2n}$ is a period-$(2n)$ of $f$.    

Note that, by combining the results in the above proofs of (b) and (c), we obtain that, if $f$ has a period-$m$ orbit $P$ with $m \ge 3$ and odd, then there exist points $v, z, z_0, d$ and, for each $n \ge 1$, a point $u_n$, a period-$(m+2n)$ point $p_{m+2n}$ and a period-$(2n)$ point $c_{2n}$ such that 
$$
v = \min \big\{\min P \le x \le f^{m-1}(\min P) : f(v) = f^{m-1}(\min P) \big\},
$$ 
$$
z \,\,\, \text{is a fixed point of} \,\,\, f \,\,\, \text{in} \,\,\, \big[v, f^{m-1}(\min P)\big],
$$ 
$$
z_0 = \min \big\{ v \le x \le z : f^2(x) = x \big\}, \,\, d = \max \big\{ \min P \le x \le v : f^2(x) = z_0 \big\} \, (> \min P),
$$ 
$$
f^2(x) < z_0 \le z < f(x) \,\,\, \text{for all} \,\,\ d < x < z_0,
$$ 
$$
\text{for each} \,\,\, n \ge 1, \,\, u_n = \min \big\{ d \le x \le v : f^{2n}(x) = d \big\} \,\,\, \text{and for all} \,\,\, 0 \le i \le n-1,
$$ 
$$
d < f^{2i}(c_{2n}) < z_0 \le z < f^{2i+1}(c_{2n}),\,\, \text{and}
$$ 
$$
\min P < d < \cdots < c_4 < u_2 < c_2 < u_1 < \cdots < p_{m+4} < p_{m+2} < v < z_0 < f^{m-1}(\min P).
$$

\vspace{.1in}
{\bf We now prove (1), (2) and (3) of Sharkovsky's theorem.}

If $f$ has a period-$m$ point with $m \ge 3$ and odd, then it follows from (b) that $f$ has a period-$(m+2)$ point and, from (c) that $f$ has periodic points of all even periods.   

If $f$ has a period-$(2 \cdot m)$ point with $m \ge 3$ and odd, then, by Lemma 2(i), $f^2$ has a period-$m$ point.  It follows from the above \big(or by (b) and (c)\big) that $f^2$ has a period-$(m+2)$ point and a period-$(2 \cdot 3)$ point. If $f^2$ has a period-$(m+2)$ point, then, by Lemma 2(ii), 
$$
f \,\,\, \text{has either a period-}(m+2) \,\,\, \text{point or a period-}\big(2 \cdot (m+2)\big) \,\,\, \text{point}.
$$  
If $f$ has a period-$(m+2)$ point, then it follows from (c) that $f$ has a period-$\big(2 \cdot (m+2)\big)$ point.  In either case, $f$ has a period-$\big(2 \cdot (m+2)\big)$ point.  On the other hand, if $f^2$ has a period-$(2 \cdot 3)$ point, then, by Lemma 2(ii), $f$ has a period-$(2^2 \cdot 3)$ point.  This shows that if $f$ has a period-$(2 \cdot m)$ point with $m \ge 3$ and odd, then $f$ has a period-$\big(2 \cdot (m+2)\big)$ point and a period-$(2^2 \cdot 3)$ point.

Now if $f$ has a period-$(2^k \cdot m)$ point with $m \ge 3$ and odd and if $k \ge 2$, then, by Lemma 2(i), $f^{2^{k-1}}$ has a period-$(2 \cdot m)$ point.  It follows from the previous paragraph that $f^{2^{k-1}}$ has a period-$\big(2 \cdot (m+2)\big)$ point and a period-$(2^2 \cdot 3)$ point.  So, by Lemma 2(ii), $f$ has a period-$\big(2^k \cdot (m+2)\big)$ point and a period-$(2^{k+1} \cdot 3)$ point.  

Furthermore, if $f$ has a period-$(2^i \cdot m)$ point with $m \ge 3$ and odd and if $i \ge 0$, then, by Lemma 2(i), $f^{2^i}$ has a period-$m$ point.  For each $\ell \ge i$, by Lemma 2(i), $f^{2^\ell} = (f^{2^i})^{2^{\ell-i}}$ has a period-$m$ point and so, by (a), $f^{2^\ell}$ has a period-2 point.  This implies, by Lemma 2(ii), that $f$ has a period-$(2^{\ell+1})$ point for each $\ell \ge i$.  

Finally, if $f$ has a period-$(2^k)$ point for some $k \ge 2$, then, by Lemma 2(i), $f^{2^{k-2}}$ has a period-4 point.  By (a), $f^{2^{k-2}}$ has a period-2 point.  By Lemma 2(ii), $f$ has a period-$(2^{k-1})$ point and hence, by induction, $f$ has a period-$(2^j)$ point for each $j = 1, 2, \cdots, k-2$.  Furthermore, it follows from (a) that $f$ has a fixed point.  This completes the proof of (1).

As for the proofs of (2) and (3), there are very elegant examples in {\bf\cite{a2, bh}}.  Here we present different examples (see {\bf\cite{du3}} for more examples).  We consider the tent map $T(x) = 1 - |2x - 1|$ and the doubly truncated tent family $\widehat T_{a,b}(x)$, where $0 < a < b < 1$, defined on $[0, 1]$ by 
$$
\widehat T_{a,b}(x) = \begin{cases}
               b, & \text{if} \,\,\, T(x) > b; \cr
               T(x), &  \text{if} \,\,\, a \le T(x) \le b; \cr
               a, & \text{if} \,\,\, T(x) < a.\cr
               \end{cases}
$$
Note that the relationship between $T(x)$ and $\widehat T_{a,b}(x)$ is that the periodic orbits of $\widehat T_{a,b}(x)$ are also periodic orbits of $T(x)$ with the same periods and, conversely, the periodic orbits of $T(x)$ which lie entirely in the interval $[a, b]$ are also periodic orbits of $\widehat T_{a,b}(x)$ with the same periods.  Consequently, if $Q_k$ is a period-$k$ orbit of $T(x)$, then it is also a period-$k$ orbit of $\widehat T_{\min Q_k, \max Q_k}(x)$.  By (1), $\widehat T_{\min Q_k, \max Q_k}(x)$ has a period-$\ell$ orbit for each $\ell$ with $k \prec \ell$.  In other words, the interval $[\min Q_k, \max Q_k]$ contains a period-$\ell$ orbit of $T(x)$ for each $\ell$ with $k \prec \ell$.  Since, for each integer $k \ge 1$, the equation $T^k(x) = x$ has exactly $2^k$ distinct solutions in $[0, 1]$, $T(x)$ has {\it finitely many} period-$k$ orbits.  Among these {\it finitely many} period-$k$ orbits, let 
$$
P_k \,\,\, \text{be one with the {\it smallest} diameter} \,\,\, \max P_k - \min P_k.
$$  
For each $x$ in $I$, let $\widehat T_k(x) = \widehat T_{a_k,b_k}(x)$, where $a_k = \min P_k$ and $b_k = \max P_k$.  Then it is easy to see that, for each $k \ge 1$, $\widehat T_k(x)$ has exactly one period-$k$ orbit (i.e., $P_k$) but has no period-$j$ orbit for any $j$ with $j \prec k$ in the Sharkovsky ordering.  This establishes (2).  

Clearly, $T(x)$ has a unique period-2 orbit, i.e., $\{ \frac 25, \frac 45 \}$.  For every periodic orbit $P$ of $T(x)$ with least period $\ge 3$, it follows from (a) that $\widehat T_{\min P, \max P}(x)$ has a period-2 orbit.  So, $\min P \le \frac 25 < \frac 45 \le \max P$.  Now let $Q_3$ be any period-3 orbit of $T(x)$ of smallest diameter.  Then $[\min Q_3, \max Q_3]$ contains finitely many period-6 orbits of $T(x)$.  If $Q_6$ is one of smallest diameter, then $[\min Q_6, \max Q_6]$ contains finitely many period-12 orbits of $T(x)$.  We choose one, say $Q_{12}$, of smallest diameter and continue the process inductively.  Let 
$$
q_0 = \sup \big\{\min Q_{2^n \cdot 3} : n \ge 0 \big\} \,\,\, \text{and} \,\,\, q_1 = \inf \big\{\max Q_{2^n \cdot 3} : n \ge 0 \big\}
$$
Then $q_0 \le \frac 25 < \frac 45 \le q_1$.  Let $\widehat T_\infty(x) = \widehat T_{q_0,q_1}(x)$ for all $0 \le x \le 1$.  If $\widehat T_\infty(x)$ had a period-$(2^n \cdot m)$ orbit for some $n \ge 0$ and some odd $m \ge 3$, then, by (1), $\widehat T_\infty(x)$ has a period-$(2^{n+1} \cdot 3)$ orbit, say $\widehat Q_{2^{n+1} \cdot 3}$.  Since $\widehat Q_{2^{n+1} \cdot 3} \subset [q_0, q_1] \subsetneq [\min Q_{2^{n+1} \cdot 3}, \max Q_{2^{n+1} \cdot 3}]$, $\widehat Q_{2^{n+1} \cdot 3}$ is also a period-$(2^{n+1} \cdot 3)$ orbit of $T(x)$ with diameter {\it smaller} than that of $Q_{2^{n+1} \cdot 3}$.  This is a contradiction.  So, $\widehat T_\infty(x)$ has no periodic orbit of period not a power of 2.  On the other hand, for each $k \ge 0$, the map $T(x)$ has finitely many period-$(2^k)$ orbits.  If each such orbit had an {\it exceptional} point which is not in the interval $[q_0, q_1]$, then it is clear that we can find an $n \ge 1$ such that the interval $[\min Q_{2^n \cdot 3}, \max Q_{2^n \cdot 3}]$ contains none of these {\it exceptional} points which implies that $[\min Q_{2^n \cdot 3}, \max Q_{2^n \cdot 3}]$ contains no period-$(2^k)$ {\it orbits} of $T(x)$.  Consequently, the map $\widehat T_{s_n,t_n}(x)$, where $s_n = \min Q_{2^n \cdot 3}, t_n = \max Q_{2^n \cdot 3}$, has no period-$(2^k)$ {\it orbits} and yet it has a period-$(2^n \cdot 3)$ orbit, i.e., $Q_{2^n \cdot 3}$. This contradicts (1).  Therefore,the map $\widehat T_\infty(x)$ is an example for (3).

Note that it is an easy consequence of the following result that, for each positive integer $m$, the map $\widehat T_{\min P_m, \max P_m}(x)$ defined above and the map 
$T_{h(m)}(x) = \min\big\{h(m), 1-|2x-1|\big\}$, where $h(m) = \min \big\{\max Q: \, Q$ is a period-$m$ orbit of $f \big\}$, defined in {\bf\cite{bh}}, although look different, have exactly one and the same period-$m$ orbit, i.e., $P_m$, and have, for each $m \prec n \ge 2$, the same collection of period-$n$ orbits.  

\noindent
{\bf Proposition 1.}
{\it Let $f : [0, 1] \rightarrow [0, 1]$ be a continuous map such that $f$ is increasing on $[0, 1/2]$ and decreasing on $[1/2, 1]$.  Suppose $f$ has a periodic orbit of least period $\ge 2$.  Then the periodic orbits of $f$ of least periods $\ge 2$ are nested in the sense that if $P$ and $Q$ are periodic orbits of $f$ of least periods $\ge 2$ with $\max P < \max Q$ then $[\min P, \max P] \subset [\min Q, \max Q]$.}

\noindent
{\it Proof.}
Let $P$ be a period-$n$ orbit of $f$ with $n \ge 2$.  Then, since $f$ is monotonic on $[0, 1/2]$ and on $[1/2, 1]$, we have 
$$
\min P \le 1/2 \le \max P.
$$
Let $p$ be a point in $P$ such that $f(p) = \min P$.  Then $1/2 \le p \le \max P$.  Since $f$ is decreasing on $[1/2, 1]$, we obtain that 
$\min P \le f(\max P) \le f(p) = \min P$.  This forces $f(\max P) = \min P$.  

Now let $P$ and $Q$ be distinct periodic orbits of $f$ of least periods $\ge 2$ with $\max P < \max Q$.  Then $1/2 \le \max P < \max Q$.  Since $f$ is decreasing on $[1/2, 1]$, we have 
$$
(\min Q =) \, f(\max Q) < f(\max P) \,\, (= \min P < \max P < \max Q).$$
This establishes that $[\min P, \max P] \subset [\min Q, \max Q]$.
$\hfill\square$

Note that Proposition 1 applies to the logistic family $f_\alpha(x) = \alpha x(1-x)$, the tent family $T_\beta(x) = \beta\big(1 - |2x-1|\big)$, the truncated tent family $T_h(x) = \min \big\{ h, 1 - |2x-1| \big\}$ and the above doubly truncated tent family $\widehat T_{a,b}(x)$. 

\noindent
{\bf Acknowledgement}\\
This work was partially supported by the Ministry of Science and Technology of Taiwan.

\pagebreak
\noindent
{\bf Appendix.}
\begin{center}
{\bf The Towers of Periodic Points of $f$ Associated with the Periodic Orbit $P$}
\end{center}

Suppose $f$ has a period-$m$ orbit $P$ with $m \ge 3$ and odd.  In this appendix, we are going to build a tower (among other towers) by using periodic points of $f$ which lie in the interval $[\min P, \max P]$ and call it the {\it basic} tower of periodic points of $f$ associated with $P$.  This {\it basic} tower is divided into 5 sections which are contained respectively in the intervals: 
$$
[\min P, d], \,\,\, [d, u_1], \,\,\, [u_1, v], \,\,\, [v, \bar u_1'] \,\,\, \text{and} \,\,\, [\bar u_1', z_0],
$$
where $u_1 = \min \big\{ d \le x \le v: f^2(x) = d \big\}$ and $\bar u_1' = \max \big\{ v \le x \le z_0: f^2(x) = d \big\}$ (the $'$ corresponds to taking the maximum of a set). It is built one layer on top the other in a {\it systematic and recursive} way and consists of countably infinitely many layers.  For each $n \ge 2$, the $n^{th}$ layer in each section consists of countably infinitely many {\it compartments}.  Each compartment consists of 3 monotonic sequences of periodic points of $f$ and the convex hulls of these monotonic sequences are pairwise disjoint ({\it the convex hull of a set is the smallest compact interval that contains the set}).  Each sequence of periodic points of $f$ in higher layers may contain finitely many (but not all) periodic points in the lower layers.  The foundation of this {\it basic} tower which is the first layer also serves to provide a simple proof of Part (1) of the celebrated Sharkovsky's cycle coexistence theorem.  
\[
\begin{array}{l}
\textrm{The construction of these various towers (including the {\it basic} tower) is relatively easy.}\\
\textrm{The hard part is to determine the least periods of these periodic points on each layer.}
\end{array}  
\]
See Lemmas 5 $\&$ 6.
$$\aleph \qquad\qquad\qquad \aleph \qquad\qquad\qquad \aleph \qquad\qquad\qquad \aleph \qquad\qquad\qquad \aleph$$  
\indent We shall follow the notations used in the main text.  Furthermore, for any distinct points $a$ and $b$ in $I$, 
$$
\text{we denote} \,\,\, [a : b] \,\,\, \text{as the compact interval with} \,\,\, a \,\,\, \text{and} \,\,\, b \,\,\, \text{as the endpoints}.
$$
In the following, our discussions are roughly divided into two parts depending on the relative locations of points with respect to the point $v$: 

Suppose $a$ and $w$ are two distinct points in $[\min P, v]$ or in $[v, z_0]$ such that $f^{n}(a) \in \big\{ z_0, f(z_0), f(d) \big\}$ and $f^n(w) = w$.  

If both $a$ and $w$ lie in the interval $[\min P, v]$, then we have $f^n\big([a : w]\big) \supset [w, z_0] \supset \{ v \}$.  So, $f^n(\nu) = v$ for some $\nu$ in $[a : w]$.  Consequently, 
$$
f^{n+2}(a) - a \ge z_0 - a > 0 \,\,\, \text{and} \,\,\, f^{n+2}(\nu) - \nu = f^2(v) - \nu = \min P - \nu < 0.
$$
Therefore, there is a point $w^*$ between $a$ and $w$ such that $f^{n+2}(w^*) = w^*$.  

If both $a$ and $w$ lie in the interval $[v, z_0]$, then we {\it may not} have $f^n\big([a : w]\big) \supset \{ v \}$.  Fortunately, we have the fact that 
{\large 
$$
\qquad\quad\qquad\qquad\,\,\, f^2(x) < x < z_0 \,\,\, \text{for all} \,\,\, v \le x < z_0. \qquad\qquad\qquad\qquad\qquad (*)
$$}
So, 
$$
f^{n+2}(a) - a \ge z_0 - a > 0 \,\,\, \text{and} \,\,\, f^{n+2}(w) - w = f^2(w) - w < 0.
$$
Therefore, there is a point $w^*$ between $a$ and $w$ such that $f^{n+2}(w^*) = w^*$. 

By repeating this process indefinitely, we obtain the following result:

\noindent
{\bf Lemma 3.}
{\it Assume that both $n$ and $k$ are positive integers such that $n-k$ is even and $\ge 0$.  Let $a$ and $w$ be distinct points in one of the intervals: $[\min P, v]$ and $[v, z_0]$.  
{\large 
$$
\text{Suppose} \,\,\, f^{k}(a) \in \big\{ z_0, f(z_0), f(d) \big\} \,\,\, \text{and} \,\,\, f^n(w) = w.
$$}
Then the following hold:
\begin{itemize}
\item[{\rm (1)}]
If $a < b$, then, for each $i \ge 0$, the point $p_{n+2i} = \min \big\{ a \le x \le b: f^{n+2i}(x) = x \big\}$ exists and \, $a < \cdots < p_{n+6} < p_{n+4} < p_{n+2} < p_n \le w < b$.  Furthermore, if all periodic points of $f$ in $[a : b]$ of odd periods have least periods $\ge k$, then, for each $i \ge 0$, the point $p_{n+2i}$ is a period-$(n+2i)$ point of $f$;

\item[{\rm (2)}]
If $b < a$, then, for each $i \ge 0$, the point $q_{n+2i} = \max \big\{ a \le x \le b: f^{n+2i}(x) = x \big\}$ exists and \, $b < w \le q_n < q_{n+2} < q_{n+4} < q_{n+6} < \cdots < a$.  Furthermore, if all periodic points of $f$ in $[a : b]$ of odd periods have least periods $\ge k$, then, for each $i \ge 0$, the point $q_{n+2i}$ is a period-$(n+2i)$ point of $f$. \\  
$$\aleph \qquad\qquad\qquad \aleph \qquad\qquad\qquad \aleph \qquad\qquad\qquad \aleph \qquad\qquad\qquad \aleph$$
\end{itemize}}

Throughout this paper, we shall call 
$$
\text{a point} \,\,\, x_0 \,\,\, \text{in} \,\,\, [\min P, z_0] \,\,\,  d\text{-point if its orbit} \,\,\, O_f(x_0) \,\,\, \text{contains the point} \,\,\, d.
$$
The $d$-points play fundamental role in the construction of various towers of periodic points of $f$ associated with $P$.  

\noindent
{\bf Lemma 4.}
{\it For each $n \ge 1$, let 
$$
u_n = \min \big\{ d \le x \le v : f^{2n}(x) = d \big\} \,\,\, \text{and}\quad
$$ 
$$
\bar u_n' = \max \big\{ v \le x \le z_0 : f^{2n}(x) = d \big\}.\qquad\,
$$
Note the relationship between the subscript $n$ of $u_n$ ($\bar u_n'$ respectively) and the supserscript $2n$ of $f^{2n}(x)$ in the definition of $u_n$ ($\bar u_n'$ respectively).  It denotes the $n^{th}$ point of the sequence $<u_n>$ ($<\bar u_n'>$ respectively) which are needed for the tower-building process to continue from the first layer to the second.  

\noindent
Then $d < \cdots < u_3 < u_2 < u_1 < v < \bar u_1' < \bar u_2' < \bar u_3' < \cdots < z_0$ and the following hold:
\begin{itemize}
\item[(1)]
For any $x_0$ in $(d, u_n)$, we have 
$$
d < f^{2i}(x_0) < z_0 \,\,\, \text{for all} \,\,\, 1 \le i \le n.
$$
Furthermore, in $[d, u_n]$, $f$ has no periodic points of odd periods $\le 2n+1$, nor has periodic points of even periods $\le 2n$ except period-$(2n)$ points (we use this kind of phrasing for a reason which will become apparant in Lemmas 8 $\&$ 12 below);

\item[(2)]
For any $x_0$ in $[\bar u_n', z_0)$, we have 
$$
v < f^{2n-2}(x_0) < \cdots < f^4(x_0) < f^2(x_0) < x_0 < z_0 \,\,\, \text{and}
$$ 
$$
d \le f^{2n}(x_0) < f^{2n-2}(x_0) < \cdots < f^2(x_0) < x_0 < z_0 < f^j(x) \,\,\, \text{for all odd} \,\,\, 1 \le j \le 2n+1;
$$

\item[(3)]
In $[\bar u_n', z_0)$, $f$ has no periodic points of even periods $\le 2n$, nor has periodic points of odd periods $\le 2n+1$.
\end{itemize}}

\noindent
{\it Proof.}
It is easy to see that $d < \cdots < u_3 < u_2 < u_1 < v < \bar u_1' < \bar u_2' < \bar u_3' < \cdots < z_0$.  By arguing as that $(\dagger)$ in the proof of (c) in the main text, we easily obtain that 
{\large 
$$
\qquad\,\,\, d < f^{2i}(x) < z_0 \,\,\, \text{for all} \,\,\, x \,\,\, \text{in} \,\,\, (d, u_n) \cup (\bar u_n', z_0) \,\,\, \text{and all} \,\,\, 1 \le i \le n. \quad\,\,\,\, (\dagger)
$$}
This, combined with the fact 
{\large 
$$
\quad\qquad\qquad f(x) > z \,\,\, (\ge z_0 > f^2(x)) \,\,\, \text{for all} \,\,\, d < x < z_0, \qquad\qquad\qquad\quad\,\, (***)
$$}
implies that $f$ has {\it no} periodic points of odd periods $\le 2n+1$ in $(d, u_n) \cup (\bar u_n', z_0)$.

On the other hand, let $x_0$ be any point in $[\bar u_n', z_0)$.  Suppose $f^{2k}(x_0) \le v$ for some $1 \le k < n$.  Then, since $f^2(z_0) = z_0$, we obtain that $f^{2k}\big([x_0, z_0]\big) \supset [v, z_0]$ and so, $f^{2k+2}\big([x_0, z_0]\big) \supset [\min P, z_0]$.  Let $\bar \nu_{k+1}$ be a point in $[x_0, z_0]$ such that $f^{2k+2}(\bar \nu_{k+1}) = \min P$ and let $\bar u_{k+1}^*$ be a point in $[\bar \nu_{k+1}, z_0]$ such that $f^{2k+2}(\bar u_{k+1}^*) = d$. Then $\bar \nu_{k+1} < \bar u_{k+1}^* < z_0$.  Since $k+1 \le n$ and $2(k+1) \le 2n$, we have $\bar u_n' \le x_0 \le \bar \nu_{k+1} < \bar u_{k+1}^* \le \bar u_{k+1}' \le \bar u_n'$.  This is a contradiction.  So, $f^{2i}(x) > v$ for all $\bar u_n' \le x < z_0$ and all $1 \le i \le n-1$.  It follows from the above $(\dagger)$ and the following (*) 
{\large 
$$
\qquad\quad\qquad\qquad\,\,\,\, f^2(x) < x < z_0 \,\,\, \text{for all} \,\,\, v \le x < z_0 \qquad\qquad\qquad\qquad\qquad (*)
$$}
that, for any $x_0$ in $[\bar u_n', z_0)$, we have 
$$
v < f^{2n-2}(x_0) < \cdots < f^4(x_0) < f^2(x_0) < x_0 < z_0 \,\,\, \text{and},
$$
by (*), we have $f^{2n}(x_0) < f^{2n-2}(x_0)$.  This, combined with $(\dagger)$ and $(***)$ above, implies that 
$$
d \le f^{2n}(x_0) < f^{2n-2}(x_0) < \cdots < f^2(x_0) < x_0 < z_0 < f^j(x) \,\,\, \text{for all odd} \,\,\, 1 \le j \le 2n+1.
$$
Consequently, we have shown (2) $\&$ (3) and part of (1).

As for the rest of (1).  It was already shown in the proof of (c) in the main text.  We include a proof here for completeness.  Suppose $f$ had a periodic point $p_{2k}^*$ of even period $2k$ with $2 \le 2k < 2n$ in $[d, u_n]$.  Then $f^{2k}\big([d, p_{2k}^*]\big) \supset [p_{2k}^*, z_0]  \supset [v, z_0]$.  So, $f^{2k+2}\big([d, p_{2k}^*]\big) \supset f^2\big([v, z_0]\big) \supset [\min P, z_0] \supset \{ d \}$.  Thus, there is a point $(d <)$ $u_{k+1}^* \, (< p_{2k}^* < u_n)$ such that $f^{2k+2}(u_{k+1}^*) = d$.  Since $k+1 \le n$, we have $u_n \le u_{k+1} \le u_{k+1}^* < u_n$.  This is a contradiction.  Therefore, $f$ has no periodic points of {\it even} periods $< 2n$ in $[d, u_n]$.

Finally, since 
$$
f^{2n}(d) - d = z_0 -d > 0 \,\,\, \text{and} \,\,\, f^{2n}(u_n) - u_n = d - u_n < 0,
$$
the point 
$$
c_{2n} = \min \big\{ d \le x \le u_n: f^{2n}(x) = x \big\}
$$
exists.  Since we have just shown that $f$ has no periodic points of even periods $< 2n$ in $[d, u_n]$, The point $c_{2n}$ must be a period-$(2n)$ point of $f$.  Therefore, $f$ has no periodic points of {\it even} periods $\le 2n$ except period-$(2n)$ points in $[d, u_n]$.  This proves (1) and hence Lemma 4.
$\hfill\square$

\noindent
{\bf Lemma 5.} 
{\it Assume that both $a$ and $b$ are distinct points in one of the intervals: $[\min P, d]$, $[d, v]$ and $[v, z_0]$.  Let $n \ge 0$ and $k \ge 0$ be integers such that $m+2n \ge 2k+1$.  
{\large 
$$
\text{Suppose} \,\,\,\,\, f^{2k+1}(a) \in \big\{ z_0, f(z_0), f(d) \big\} \,\,\, \text{and} \,\,\, f^{m+2n}(b) = \min P.
$$}
Then the following hold: 
\begin{itemize}
\item[{\rm (1)}]
If $a < b$, then, for each $i \ge 0$, the points 
$$
p_{m+2n+2i} = \min \big\{ a \le x \le b: f^{m+2n+2i}(x) = x \big\},
$$ 
$$
\quad \mu_{m,n+i} = \min \big\{ a \le x \le b: f^{m+2n+2i}(x) = d \big\}\,\,\, \text{and} \,\,
$$ 
$$\, c_{2m+2n+2i}^{\,(*)} = \min \big\{ a \le x \le b: f^{2m+2n+2i}(x) = x \big\}
$$
exist (note the relationship between the subscript of $\mu_{m,n+i}$ and the superscript of $f^{m+2n+2i}$) and
$$
a < \cdots < c_{2m+2n+6}^{\,(*)} < c_{2m+2n+4}^{\,(*)} < c_{2m+2n+2}^{\,(*)} < c_{2m+2n}^{\,(*)} < b,
$$ 
$$
a < \cdots < p_{m+2n+4} < \mu_{m,n+1} < p_{m+2n+2} < \mu_{m,n} < p_{m+2n} < b \,\,\, \text{if} \,\,\, \min P < a < b < d,
$$ 
$$
a < \cdots < \mu_{m,n+2} < p_{m+2n+2} < \mu_{m,n+1} < p_{m+2n} < \mu_{m,n} < b\,\,\, \text{if} \,\,\, d \le a < b \le v
$$
and, if $v < a < b < z_0$, then 
\begin{multline*}
$$
a < \cdots < \mu_{m,n+2} < \mu_{m,n+1} < \mu_{m,n} < b, \\ 
a < \cdots < p_{m+2n+4} < p_{m+2n+2} < p_{m+2n} \,\,\, \text{and} \\ 
p_{m+2n+2i} < \mu_{m,n+i} \,\,\, \text{for all} \,\,\, i \ge 0.
$$
\end{multline*}
Furthermore, $$\text{if} \,\,\, i \ge 0 \,\,\, \text{and} \,\,\, m+2n+2i \ge 3(2k+1), \, \text{then} \,\,\, p_{m+2n+2i} \,\,\, \text{is a period-}(m+2n+2i) \,\,\, \text{point of} \,\,\, f,$$and, if all periodic points of $f$ in $[a, b]$ of odd periods have least periods $\ge 2k+1$, then, 
\begin{itemize}
\item[{\rm (i)}]
for each $i \ge 0$, $p_{m+2n+2i}$ is a period-$(m+2n+2i)$ point of $f$; 

\item[{\rm (ii)}]
if $\min P < a < b < d$ or $d \le a < b \le v$, then, for each $i \ge n+1$, $c_{2m+2n+2i}^{\,(*)}$ is a period-$(2m+2n+2i)$ point of $f$ and

\item[{\rm (iii)}]
if $v < a < b < z_0$, then, for each $i \ge n+ \max \big\{ 1, m+1-2n \big\}$, $c_{2m+2n+2i}^{\,(*)}$ is a period-$(2m+2n+2i)$ point of $f$.
\end{itemize}

\item[{\rm (2)}]
If $b < a$, then, for each $i \ge 0$, the points 
$$
q_{m+2n+2i} = \max \big\{ b \le x \le a: f^{m+2n+2i}(x) = x \big\},\,
$$ 
$$
\,\,\,\, \mu_{m,n+i}' = \max \big\{ b \le x \le a: f^{m+2n+2i}(x) = d \big\} \,\,\, \text{and}
$$ 
$$
\, c_{2m+2n+2i}'^{\,(*)} = \max \big\{ b \le x \le a: f^{2m+2n+2i}(x) = x \big\}
$$
exist (note the relationship between the subscript of $\mu_{m,n+i}'$ and the superscript of $f^{m+2n+2i}$) and 
$$
b < c_{m+2n}'^{\,(*)} < c_{m+2n+2}'^{\,(*)} < c_{m+2n+4}'^{\,(*)} < c_{m+2n+6}'^{\,(*)} < \cdots < a,
$$ 
$$
b < q_{m+2n} < \mu_{m,n}' < q_{m+2n+2} < \mu_{m,n+1}' < q_{m+2n+4} < \cdots < a \,\,\, \text{if} \,\,\, \min P \le b < a \le d,
$$ 
$$
b < \mu_{m,n}' < q_{m+2n} < \mu_{m,n+1}' < q_{m+2n+2} < \mu_{m,n+2}' < \cdots < a \,\,\, \text{if} \,\,\, d < b < a < v
$$
and, if $v \le b < a \le z_0$, then 
\begin{multline*}
$$
b < \mu_{m,n}' < \mu_{m,n+1}' < \mu_{m,n+2}' < \cdots < a, \\ 
b < q_{m+2n} < q_{m+2n+2} < q_{m+2m+4} < \cdots < a \,\,\, \text{and} \\ 
\mu_{m,n+i}' < q_{m+2n+2i} \,\,\, \text{for all} \,\,\, i \ge 0.
$$
\end{multline*}
Furthermore, $$\text{if} \,\,\, i \ge 0 \,\,\, \text{and} \,\,\, m+2n+2i \ge 3(2k+1), \, \text{then} \,\,\, q_{m+2n+2i} \,\,\, \text{is a period-}(m+2n+2i) \,\,\, \text{point of} \,\,\, f,$$and, if all periodic points of $f$ in $[b, a]$ of odd periods have least periods $\ge 2k+1$, then, 
\begin{itemize}
\item[{\rm (i)}]
for each $i \ge 0$, $q_{m+2n+2i}$ is a period-$(m+2n+2i)$ point of $f$; 

\item[{\rm (ii)}]
if $\min P \le b < a \le d$ or $d < b < a < v$, then, for each $i \ge n+1$, $c_{2m+2n+2i}'^{\,(*)}$ is a period-$(2m+2n+2i)$ point of $f$ and

\item[{\rm (iii)}]
if $v \le b < a \le z_0$, then, for each $i \ge n+ \max \big\{ 1, m+1-2n \big\}$, $c_{2m+2n+2i}'^{\,(*)}$ is a period-$(2m+2n+2i)$ point of $f$.
\end{itemize}
\end{itemize}}

\noindent
{\bf Remark 1.}
In the above result, the periodic points, say, $c_{2m+2n+2i}^{\,(*)}$'s of even periods are {\it interspersed} with the periodic points $p_{m+2n+2i}$'s of odd periods.  To make things simple, we shall use these periodic points $p_{m+2n+2i}$'s of odd periods (because we apply odd iterates of $f$ to the endpoints $a$ and $b$) and ignore those periodic points $c_{2m+2n+2i}^{\,(*)}$'s of even periods later on to build the {\it basic} tower of periodic points of $f$ associated with $P$.

\noindent
{\it Proof.}
Assume that $\min P < a < b < d$ or $d \le a < b \le v$.  Since $f^{2k+1}(a) \in \big\{ z_0, f(z_0), f(d) \big\}$ and, $m+2n \ge 2k+1$, we have $f^{m+2n}(a) \in \big\{ z_0, f(z_0), f(d) \big\}$.  
{\large 
$$
\text{Suppose} \,\,\, \min P < a < b < d.
$$}
Then $\min P < x < d$ for all $x$ in $[a, b]$.  Since 
$$
f^{m+2n}(a) - a \ge z_0 - a > 0 \,\,\, \text{and} \,\,\, f^{m+2n}(b) - b = \min P - b < 0,
$$
the point $p_{m+2n} = \min \big\{ a \le x \le b: f^{m+2n}(x) = x \big\}$ exists.  Since $f^{m+2n}(a) \in \big\{ z_0, f(z_0), f(d) \big\}$ and $f^{m+2n}(p_{m+2n}) = p_{m+2n}$, it follows from Lemma 3 that, for each $i \ge 0$, the point $p_{m+2n+2i} = \min \big\{ a \le x \le b: f^{m+2n+2i}(x) = x \big\}$ exists.

Furthermore, for each $i \ge 0$, since (note that $p_{m+2n+2i} < b \le d$) 
$$
f^{m+2n+2i}(a) - d \ge z_0 - d > 0 \,\,\, \text{and} \,\,\, f^{m+2n+2i}(p_{m+2n+2i}) - d = p_{m+2n+2i} - d < 0,
$$
the point $\mu_{m,n+i} = \min \big\{ a \le x \le p_{m+2n+2i}: f^{m+2n+2i}(x) = d \big\}$ exists and is $< p_{m+2n+2i}$.

On the other hand, since $f^{m+2n+2i}\big([a, \mu_{m,n+i}]\big) \supset [d, z_0] \supset \{ v \}$, there is a point $\nu_{m,n+i}$ in $[a, \mu_{m,n+i}]$ such that $f^{m+2n+2i}(\nu_{m,n+i}) = v$.  Since $f^{m+2n+2i+2}(a) - a \ge z_0 - a > 0$ and $f^{m+2n+2i+2}(\nu_{m,n+i}) - \nu_{m,n+i} = f^2(v) - \nu_{m,n+i} = \min P - \nu_{m,n+i} < 0$, the point 
$$
p_{m+2n+2i+2} = \min \big\{ a \le x \le \nu_{m,n+i}: f^{m+2n+2i+2}(x) = x \big\} \,\,\, \text{exists and is} \,\, < \nu_{m,n+i} < \mu_{m,n+i}.
$$
This, combined with the above, implies that 
$$
min P < a < \cdots < \mu_{m,n+2} < p_{m+2n+4} < \mu_{m,n+1} < p_{m+2n+2} < \mu_{m,n} < p_{m+2n} < b < d.
$$  
\vspace{.05in}
\[
\left( \begin{array}{l}
\textrm{\qquad\qquad\qquad Suppose $d \le a < b \le v$.  Then $d \le x \le v$ for all $x$ in $[a, b]$.}\\
\textrm{}\\
\textrm{Since $f^{m+2n}\big([a, b]\big) \supset [\min P, z_0] \supset \{ d \}$, the point}\\
\textrm{\qquad\qquad\qquad $\mu_{m,n} = \min \big\{ a \le x \le b: f^{m+2n}(x) = d \big\}$ exists.}\\
\textrm{Since $f^{m+2n}(a) - a \ge z_0 - a > 0$ and $f^{m+2n}(\mu_{m,n}) - \mu_{m,n} = d - \mu_{m,n} < 0$, the point}\\
\textrm{\qquad\qquad\qquad $p_{m+2n} = \min \big\{ a \le x \le \mu_{m,n}: f^{m+2n}(x) = x \big\}$ exists and is $< \mu_{m,n}$.}\\
\textrm{Since $f^{m+2n}\big([a, p_{m+2n}]\big) \supset [p_{m+2n}, z_0] \supset \{ v \}$, there is a point}\\
\textrm{$\nu_{m,n}$ in $[a, p_{m+2n}]$ such that $f^{m+2n}(\nu_{m,n}) = v$.  Since $f^{m+2n+2}(a) - d \ge z_0 - d > 0$}\\
\textrm{and $f^{m+2n+2}(\nu_{m,n}) - d = f^2(v) - d = \min P - d < 0$, the point}\\
\textrm{\qquad $\mu_{m,n+1} = \min \big\{ a \le x \le \nu_{m,n}: f^{m+2n+2}(x) = d \big\}$ exists and is $< \nu_{m,n} < p_{m+2n}$.}\\
\textrm{Since $f^{m+2n+2}(a) - a \ge z_0 - a > 0$ and $f^{m+2n+2}(\mu_{m,n+1}) - \mu_{m,n+1} = d - \mu_{m,n+1} < 0$,}\\
\textrm{the point $p_{m+2n+2} = \min \big\{ a \le x \le \mu_{m,n+1}: f^{m+2n+2}(x) = x \big\}$ exists and is $< \mu_{m,n+1}$.}\\
\textrm{By induction, we have}\\
\textrm{$d \le a < \cdots < p_{m+2n+4} < \mu_{m,n+2} < p_{m+2n+2} < \mu_{m,n+1} < p_{m+2n} < \mu_{m,n} < b \le v$.}
\end{array} \right)  
\]

\indent Note that we have assumed that $\min P < a < b < d$ or $d \le a < b \le v$.  Since $f^{2m+2n}(a) - a \ge z_0 - a > 0$ and $f^{2m+2n}(b) - b = f^m\big(f^{m+2n}(b)\big) - b = f^m(\min P) - b = \min P - b < 0$, the point $c_{2m+2n}^{\,(*)} = \min \big\{ a \le x \le b: f^{2m+2n}(x) = x \big\}$ exists.  Since $f^{2m+2n}(a) \in \big\{ z_0, f(z_0) \big\}$ and $f^{2m+2n}(c_{2m+2n}^{\,(*)}) = c_{2m+2n}^{\,(*)}$, it follows from Lemma 3 that, for each $i \ge 0$, the point $c_{2m+2n+2i}^{\,(*)} = \min \big\{ a \le x \le b: f^{2m+2n+2i}(x) = x \big\}$ exists and $a < \cdots < c_{2m+2n+6}^{\,(*)} < c_{2m+2n+4}^{\,(*)} < c_{2m+2n+2}^{\,(*)} < c_{2m+2n}^{\,(*)} < b$.  (Note that these periodic points $c_{2m+2n+2i}^{\,(*)}$'s are {\it interspersed} with the periodic points $p_{m+2n+2i}$'s of odd periods.  To make things simple, we shall ignore these periodic points $c_{2m+2n+2i}^{\,(*)}$'s and use only the periodic points $p_{m+2n+2i}$'s of odd periods later on to build the {\it basic} tower of periodic points of $f$ associated with $P$. We include them here for the interested readers).

Now, we want to find the least periods of $p_{m+2n+2i}$'s with respect to $f$.  For any periodic point $p$ of $f$, let $\ell(p)$ denote the least period of $p$ with respect to $f$.  Suppose, for some $j \ge 0$, $\ell(p_{m+2n+2j}) < m+2n+2j$.  By Lemma 1, $\ell(p_{m+2n+2j})$ divides $m+2n+2j$ and so, is odd.
\begin{itemize}
\item[{\rm (i)}]
If $\ell(p_{m+2n+2j}) \ge 2k+1$, then, since $\ell(p_{m+2n+2j})$ is odd and $f^{2k+1}(a) \in \big\{ z_0, f(z_0), f(d) \big\}$, we have $$
f^{\ell(p_{m+2n+2j})}(a) \in \big\{ z_0, f(z_0), f(d) \big\} \,\, \text{and} \,\, f^{\ell(p_{m+2n+2j})}(p_{m+2n+2j}) = p_{m+2n+2j}.
$$
Since $m+2n+2j > \ell(p_{m+2n+2j})$, it follows from Lemma 3 that there is a periodic point $w_{m+2n+2j}$ of $f$ in $[a, p_{m+2n+2j})$ such that $f^{m+2n+2j}(w_{m+2n+2j}) = w_{m+2n+2j}$.  Since $a < w_{m+2n+2j} < p_{m+2n+2j} < b$, this contradicts the minimality of $p_{m+2n+2j}$ in $[a, b]$.  This shows that if $j \ge 0$ is such that $\ell(p_{m+2n+2j}) < m+2n+2j$, then $\ell(p_{m+2n+2j}) < 2k+1$.

\item[{\rm (ii)}]
Suppose $\ell(p_{m+2n+2j}) < 2k+1$ and let $r = (m+2n+2j)/\ell(p_{m+2n+2j})$.  By Lemma 1, $r \ge 1$ is odd and, since $m+2n+2j > \ell(p_{m+2n+2j})$, $r \ge 3$.
$$
\text{Suppose} \,\,\, (m+2n+2j)/(2k+1) \ge 3.
$$
Then $m+2n+2j - 2\ell(p_{m+2n+2j}) = (r-2)\ell(p_{m+2n+2j})$ is odd and $\ge 3(2k+1) - 2\ell(p_{m+2n+2j}) = (2k+1) + 2[2k+1 - \ell(p_{m+2n+2j})] > 2k+1$.  Since 
$$
f^{(m+2n+2j) - 2\ell(p_{m+2n+2j})}(a) \in \big\{ z_0, f(z_0) \big\} \, \text{and} \, f^{(m+2n+2j) - 2\ell(p_{m+2n+2j})}(p_{m+2n+2j}) = p_{m+2n+2j},
$$
and since $m+2n+2j > m+2n+2j - 2\ell(p_{m+2n+2j})$, it follows from Lemma 3 that there is a periodic point $w_{m+2n+2j}$ of $f$ in $[a, p_{m+2n+2j})$ such that $f^{m+2n+2j}(w_{m+2n+2j}) = w_{m+2n+2j}$.  Since $a < w_{m+2n+2j} < p_{m+2n+2j} < b$, this contradicts the minimality of $p_{m+2n+2j}$ in $[a, b]$.  
\end{itemize}
By combining the above (i) and (ii), we obtain that if $i \ge 0$ is an integer such that $m+2n+2i \ge 3(2k+1)$, then $p_{m+2n+2i}$ is a period-$(m+2n+2i)$ point of $f$.  
$$.............. \qquad .............. \qquad .............. \qquad ..............$$
\indent Assume that all periodic points of $f$ in $[a,b ]$ of odd periods have least periods $\ge 2k+1$.  Suppose $\ell(p_{m+2n+2j}) < m+2n+2j$ for some $j \ge 0$.  Then, by hypothesis, $\ell({p_{m+2n+2j}}) \ge 2k+1$.  

\pagebreak \noindent By arguing as those in (i) above, we obtain that
\begin{multline*}
$$
\text{if all periodic points of} \,\,\, f \,\,\, \text{in} \,\,\, [a, b] \,\,\, \text{of odd periods have least periods} \,\, \ge 2k+1, \, \text{then} \\ 
\text{for each} \,\,\, i \ge 0, \,\,\, \text{the point} \,\,\, p_{m+2n+2i} \,\,\, \text{is a period-}(m+2n+2i) \,\,\, \text{point of} \,\,\, f.
$$
\end{multline*}
$$.............. \qquad .............. \qquad .............. \qquad ..............$$
\indent We now determine the least periods of $c_{2m+2n+2i}^{\,*}$'s under the assumption that all periodic points of $f$ in $[a, b]$ of odd periods have least periods $\ge 2k+1$.  

Recall that, we have assumed that $\min P < a < b < d$ or $d \le a < b \le v$ and, \\ for each $i \ge 0$, $c_{2m+2n+2i}^{\,*} = \min \big\{ a \le x \le b: f^{2m+2n+2i}(x) = x \big\}$.
$$
\text{Suppose} \,\,\, i \ge 0 \,\,\, \text{and} \,\,\, m+n+i \ge 2k+2 \,\,\, \text{and let} \,\,\, r = (2m+2n+2i)/\ell(c_{2m+2n+2i}^{\,*}).
$$
If $r \ge 3$ is odd, then $r-1 \ge 2$ is even.  Since 
$$
(r-1)\ell(c_{2m+2n+2i}^{\,*}) = (r-1)[(2m+2n+2i)/r] \ge 2[(r-1)/r] \cdot (2k+2) \ge 2k+2,
$$
we can apply Lemma 3 with
$$
f^{2k+2}(a) \in \big\{ z_0, f(z_0) \big\} \,\,\, \text{and} \,\,\, f^{(r-1)\ell(c_{2m+2n+2i}^{\,*})}(c_{2m+2n+2i}^{\,*}) = c_{2m+2n+2i}^{\,*}.
$$
If $r \ge 4$ is even, then $r-2 \ge 2$ is even.  Since 
$$
(r-2)\ell(c_{2m+2n+2i}^{\,*}) = (r-2)[(2m+2n+2i)/r] \ge 2[(r-2)/r] \cdot (2k+2) \ge 2k+2,
$$
we can apply Lemma 3 with
$$
f^{2k+2}(a) \in \big\{ z_0, f(z_0) \big\} \,\,\, \text{and} \,\,\, f^{(r-2)\ell(c_{2m+2n+2i}^{\,*})}(c_{2m+2n+2i}^{\,*}) = c_{2m+2n+2i}^{\,*}.
$$
\indent In either case, we obtain a periodic point $w_{2m+2n+2i}$ of $f$ in $[a, c_{2m+2n+2i}^{\,*})$ such that $f^{2m+2n+2i}(w_{2m+2n+2i}) = w_{2m+2n+2i}$.  This contradicts the minimality of $c_{2m+2n+2i}^{\,*}$ in $[a, b]$.  So, if $r > 2$, then $\ell(c_{2m+2n+2i}^{\,*}) = 2m+2n+2i$.

If $r = 2$ and $\ell(c_{2m+2n+2i}^{\,*}) \, (= m+n+i \ge 2k+2)$ is even, then by applying Lemma 3 with 
$$
f^{2k+2}(a) \in \big\{ z_0, f(z_0) \big\} \,\,\, \text{and} \,\,\, f^{\ell(c_{2m+2n+2i}^{\,*})}(c_{2m+2n+2i}^{\,*}) = c_{2m+2n+2i}^{\,*},
$$
we obtain the same contradiction.  So, in this case, $\ell(c_{2m+2n+2i}^{\,*}) = 2m+2n+2i$.

\noindent
Consequently, if $i \ge n+1$ and $m+n+i$ is even, then $m+n+i \ge (m+2n)+1 \ge (2k+1)+1 \ge 2k+2$.  In this case, since we have either $r > 2$ or $r = 2$, it follows from the above that $c_{2m+2n+2i}^{\,*}$ is a period-$(2m+2n+2i)$ point of $f$.  

On the other hand, suppose $i \ge n+1$ and $m+n+i$ is odd, then $m+n+i \ge m+2n+1 \ge 2k+2$.  If $r > 2$, then it follows form the above that $c_{2m+2n+2i}^{\,*}$ is a period-$(2m+2n+2i)$ point of $f$.  If $r = 2$, then $\ell(c_{2m+2n+2i}^{\,*}) = m+n+i$.  So, in this case, $c_{2m+2n+2i}^{\,*}$ is either a period-$(2m+2n+2i)$ point of $f$ or a period-$(m+n+i)$ point of $f$.

Suppose $m+n+i$ is odd and $\ge 2k+2$ and the point $c_{2m+2n+2i}^{\,*}$ is a period-$(m+n+i)$ point of $f$.  Since $m+n+i \ge 2k+2$ and $f^{2k+1}(a) \in \big\{ z_0, f(z_0), f(d) \big\}$, we have $f^{m+n+i}(a) \in \big\{ z_0, f(z_0) \big\}$ and, since we have assumed that $\min P < a < b < d$ or $d \le a < b \le v$, we obtain that
$$
f^{m+n+i}\big([a, c_{2m+2n+2i}^{\,*}]\big) \supset \big[f^{m+n+i}(c_{2m+2n+2i}^{\,*}), f^{m+n+i}(a)\big] \supset [c_{2m+2n+2i}^{\,*}, z_0] \supset [v, z_0].
$$
Consequently, we have 
$$
f^{m+n+i+2}\big([a, c_{2m+2n+2i}^{\,*}]) \supset f^2([v, z_0]\big) \supset [\min P, z_0].
$$
Since $f^2\big([\min P, z_0]\big) \supset f^2\big([v, z_0]\big) \supset [\min P, z_0]$, we obtain that, 
$$
\text{for each} \,\,\, i \ge n+1 \,\,\, \text{such that} \,\, m+n+i \,\,\, \text{is odd},
$$ 
$n+i$ is even and $\ge n+(n+1) \ge 1$.  So, $n+i \ge 2$ and $f^{2m+2n+2i}\big([a, c_{2m+2n+2i}^{\,*}]\big) \supset f^{m+n+i-2}\big([\min P, z_0]\big) \supset f^m\big(f^{n+i-2}([\min P, z_0])\big) \supset f^m\big([\min P, z_0]\big)$ $\supset [\min P, z_0]$ $\supset \big\{ \min P \big\}$.  

Let $\delta_{2m+2n+2i}$ be a point in the interval $[a, c_{2m+2n+2i}^{\,*}]$ such that $f^{2m+2n+2i}(\delta_{2m+2n+2i}) = \min P$.  It is clear that $a < \delta_{2m+2n+2i} < c_{2m+2n+2i}^{\,*} < b$.  Since 
$$
f^{2m+2n+2i}(a) - a \ge z_0 - a > 0 \,\,\, \text{and} \,\,\, f^{2m+2n+2i}(\delta_{2m+2n+2i}) - \delta_{2m+2n+2i} = \min P - \delta_{2m+2n+2i} < 0,
$$
there exists a point $w_{2m+2n+2i}$ in $[a, \delta_{2m+2n+2i}]$ such that $f^{2m+2n+2i}(w_{2m+2n+2i}) = w_{2m+2n+2i}$.  Since $a < w_{2m+2n+2i} < \delta_{2m+2n+2i} < c_{2m+2n+2i}^{\,*} < b$ and $f^{2m+2n+2i}(w_{2m+2n+2i}) = w_{2m+2n+2i}$, this contradicts the minimality of $c_{2m+2n+2i}^{\,*}$ in $[a, b]$.  So, the point $c_{2m+2n+2i}^{\,*}$ can not be a period-$(m+n+i)$ point of $f$ and must be a period-$(2m+2n+2i)$ point of $f$.  Therefore, we have shown that 
\begin{multline*}
$$
\text{if all periodic points of} \,\,\, f \,\,\, \text{in} \,\,\, [a, b] \,\,\, \text{of odd periods have least periods} \,\, \ge 2k+1, \, \text{then}, \\ 
\text{for each} \,\,\, i \ge n+1, \, c_{2m+2n+2i}^{\,*} \,\,\, \text{is a period-}(2m+2n+2i) \,\,\, \text{point of} \,\,\, f.
$$
\end{multline*}
$$........................................................................$$
{\large 
$$
\text{Suppose} \,\,\, v < a < b < z_0.
$$}
In this case, we only need to show that, under the assumption that all periodic points of $f$ in $[a, b]$ of odd periods have least periods $\ge 2k+1$, for each $i \ge n + \max \big\{ 1, m+1 - 2n \big\}$, the point $c_{2m+2n+2i}^{\,(*)} = \min \big\{ a \le x \le b: f^{2m+2n+2i}(x) = x \big\}$ is a period-$(2m+2n+2i)$ point of $f$.  The rest statements regarding this case can be easily proved and are omitted.

By arguments similar to the case as $\min P < a < b < d$, we obtan that

\noindent
\,\,\, if $m+n+i \, (\ge 2k+2)$ is even, then $c_{2m+2n+2i}^{\,*}$ is a period-$(2m+2n+2i)$ point of $f$ and
\begin{multline*}
$$
\text{if} \,\,\, m+n+i \, (> 2k+2) \,\,\, \text{is odd, then} \,\,\, c_{2m+2n+2i}^{\,*} \,\,\, \text{is either a period-}(2m+2n+2i) \,\,\, \text{point or} \\ 
\text{a period-}(m+n+i) \,\,\, \text{point of} \,\,\, f.
$$
\end{multline*}
Now we show that, 
{\large 
$$
\text{for each} \,\, i \ge 1 \,\, \text{such that} \,\,\, 2(n+i) \ge m+1,
$$}
the point $c_{2m+4n+4i}^{\,*}$ is a period-$(2m+4n+4i)$ point of $f$ by arguing as follows (note that $m+2(n+i)=m+2n+2i$ \, is odd and $= (m+2n)+2i \ge (2k+1)+2 > 2k+2$):

Recall that, for each $i \ge 0$, $p_{m+2n+2i} = \min \big\{ a \le x \le b : f^{m+2n+2i}(x) = x \big\}$ is the {\it smallest} point in $[a, b]$ that satisfies the equation $f^{m+2n+2i}(x) = x$.  
$$
\text{Suppose, for some} \,\,\, j \ge 1 \,\,\, \text{such that} \,\,\, 2(n+j) \ge m+1, \, \ell(c_{2m+4n+4j}^{\,*}) = m+2n+2j \,\,\, \text{(is odd)}.
$$
The following arguments are different from those used in the proof of the case $\min P < a < b < d$ above and are based on the following two facts 
{\large 
$$
\qquad\quad\qquad\qquad\,\,\, f^2(x) < x < z_0 \,\,\, \text{for all} \,\,\, v \le x < z_0 \qquad\qquad\qquad\qquad\qquad (*)
$$}
and 
{\large 
$$
\quad\qquad\qquad f(x) > z \,\,\, (\ge z_0 > f^2(x)) \,\,\, \text{for all} \,\,\, d < x < z_0. \qquad\qquad\qquad\quad\,\, (***)
$$}
Since $v < a < p_{m+2n+2j} < b$, if the points $f^2(p_{m+2n+2j})$, $f^4(p_{m+2n+2j})$, $\cdots, f^{m-1}(p_{m+2n+2j})$ are all $> v$, then by (*) above, we have 
$$
v < f^{(m-1)}(p_{m+2n+2j}) < f^{(m-3)}(p_{m+2n+2j}) < \cdots < f^2(p_{m+2n+2j}) < p_{m+2n+2j} < b < z_0.
$$
But then, by (***) above, $(z_0 > b >) \, p_{m+2n+2j} = f^m(p_{m+2n+2j}) = f\big(f^{(m-1)}(p_{m+2n+2j})\big) > z \ge z_0$.  This is a contradiction.  Therefore, $f^{2s}(p_{m+2n+2j}) \le v$ for some $1 \le s \le (m-1)/2$.  Since $f^{m+2n+2j}\big([a, p_{m+2n+2j}]\big) \supset [p_{m+2n+2j}, \, z_0]$, we have 
$$
f^{m+2n+2j+2s+2}\big([a, p_{m+2n+2j}]\big) \supset f^{2s+2}\big([p_{m+2n+2j}, \, z_0]\big) \supset f^2\big([v, z_0]\big) \supset [\min P, z_0].
$$
Since $f^2\big([\min P, z_0]\big) \supset f^2\big([v, z_0]\big) \supset [\min P, z_0]$ and $m-1 \, (\ge 2s)$ is even, we obtain that, 
$$
f^{m+2n+2j+(m+1)}\big([a, p_{m+2n+2j}]\big) = f^{(m+2n+2j+2s+2)+(m-1-2s)}\big([a, p_{m+2n+2j}]\big) \supset [\min P, z_0].
$$
This implies that, since $j \ge 1$ and $2(n+j) \ge  m+1$ (so, $m+2n+2j \ge (m+2n)+2j > 2k+2$),
\begin{multline*}
$$
f^{2m+4n+4i}\big([a, p_{m+2n+2i}]\big) = f^{m+2n+2i-m-1}\big(f^{m+2n+2i+(m+1)}([a, p_{m+2n+2i}])\big) \\  
\supset f^{m+2n+2i-(m+1)}\big([\min P, z_0]\big) \supset f^m\big([\min P, z_0]\big) \supset [\min P, z_0].
$$
\end{multline*} 

Let $\delta_{2m+4n+4j}$ be a point in $[a, p_{m+2n+2j}] \, (\subset [a, b])$ such that $f^{2m+4n+4j}(\delta_{2m+4n+4j}) = \min P$.  It is clear that $a < \delta_{2m+4n+4j} < p_{m+2n+2j} < b$.  Since 
$$
f^{2m+4n+4j}(a) - a \ge z_0 - a > 0 \,\,\, \text{and} \,\,\, f^{2m+4n+4j}(\delta_{2m+4n+4j}) - \delta_{2m+4n+4j} = \min P - \delta_{2m+4n+4j} < 0,
$$
there is a point $w_{2m+4n+4j}$ in $[a, \delta_{2m+4n+4j}] \, (\subset [a, b])$ so that $f^{2m+4n+4j}(w_{2m+4n+4j}) = w_{2m+4n+4j}$.  It follows from the definition of $c_{2m+4n+4j}^{\,*}$ that $a < c_{2m+4n+4j}^{\,*} \le w_{2m+4n+4j} \, (\le \delta_{2m+4n+4j}] < p_{m+2n+2j}) < b$.  Since $p_{m+2n+2j}$ is the {\it smallest} point in $[a, b]$ that satisfies the equation $f^{m+2n+2j}(x) = x$, we obtain that $c_{2m+4n+4j}^{\,*}$ can not be a period-$(m+2n+2j)$ point of $f$.  This is a contradiction.  Therefore, we have shown that, for each $i \ge 1$ such that $2(n+i) = 2n+2i \ge m+1$, $c_{2m+4n+4i}^{\,*}$ is a period-$(2m+4n+4i)$ point of $f$.  In other words, we have shown that
$$
\text{for each {\it even}} \,\,\, i \ge \max \big\{ 1, m+1-2n \big\}, \, c_{2m+4n+2i}^{\,*}\,\,\, \text{is a period-}(2m+4n+2i) \,\,\, \text{point of} \,\,\, f.
$$
On the other hand, for each {\it odd} \, $i \ge 1$, $m+2n+i$ $(\ge 2k+2)$ is even.  It follows from the above that 
$$
\text{for each {\it odd}} \,\,\, i \ge 1, \, c_{2m+4n+2i}^{\,*}\,\,\, \text{is a period-}(2m+4n+2i) \,\,\, \text{point of} \,\,\, f.
$$
By combining the above two results, we obtain that, for each $i \ge \max \big\{ 1, m+1-2n \big\}$, $i$ even or odd, $c_{2m+4n+2i}^{\,*}$ is a period-$(2m+4n+2i)$ point of $f$ or, equivalently, 
$$
\text{for each} \,\,\, i \ge n+ \max \{ 1, m+1-2n \},
$$
the point $c_{2m+2n+2i}^{\,*}$ is a period-$(2m+2n+2i)$ point of $f$.  

This completes the proof of (1).  (2) can be proved similarly.
$\hfill\square$

\vspace{.1in}

\noindent
{\bf Lemma 6.} 
{\it Assume that both $a$ and $b$ are distinct points in one of the intervals: $[\min P, d]$, $[d, v]$ and $[v, z_0]$.  Let $n \ge k \ge 1$ be integers.  
{\large 
$$
\text{Suppose} \,\,\,\,\, f^{2k}(a) \in \big\{ z_0, f(z_0), f(d) \big\} \,\,\, \text{and} \,\,\, f^{2n}(b) = \min P.
$$}
Then the following hold: 
\begin{itemize}
\item[{\rm (1)}]
If $a < b$, then, for each $i \ge 0$, the points 
$$
\, c_{2n+2i} = \min \big\{ a \le x \le b: f^{2n+2i}(x) = x \big\}\,\,\, \text{and}
$$ 
$$
\,\,\, \mu_{n+i} = \min \big\{ a \le x \le b: f^{2n+2i}(x) = d \big\} \qquad\quad
$$
exist (note the relationship between the subscript of $\mu_{n+i}$ and the superscript of $f^{2n+2i}$) and
$$
a < \cdots < \mu_{n+2} < c_{2n+4} < \mu_{n+1} < c_{2n+2} < \mu_{n} < c_{2n} < b \,\,\, \text{if} \,\,\, \min P < a < b < d,
$$ 
$$
a < \cdots < c_{2n+4} < \mu_{n+2} < c_{2n+2} < \mu_{n+1} < c_{2n} < \mu_{n} < b \,\,\, \text{if} \,\,\, d \le a < b \le v
$$
and, if $v < a < b < z_0$, then 
\begin{multline*}
$$
a < \cdots < \mu_{n+3} < \mu_{n+2} < \mu_{n+1} < \mu_{n} < b, \\ 
a < \cdots < c_{2n+6} < c_{2n+4} < c_{2n+2} < c_{2n} \,\,\, \text{and} \\ c_{2n+2i} < \mu_{n+i} \,\,\, \text{for all} \,\,\, i \ge 0.
$$
\end{multline*}

\item[{\rm (2)}]
If $b < a$, then, for each $i \ge 0$, the points 
$$
\, c_{2n+2i}' = \max \big\{ b \le x \le a: f^{2n+2i}(x) = x \big\}\,\,\, \text{and}
$$ 
$$\,\,\, \mu_{n+i}' = \max \big\{ b \le x \le a: f^{2n+2i}(x) = d \big\} \qquad\quad
$$ 
exist (note the relationship between the subscript of $\mu_{n+i}'$ and the superscript of $f^{2n+2i}$) and
$$
b < c_{2n}' < \mu_{n}' < c_{2n+2}' < \mu_{n+1}' < c_{2n+4}' < \cdots < a \,\,\, \text{if} \,\,\, \min P \le b < a \le d,
$$ 
$$
b < \mu_{n}' < c_{2n}' < \mu_{n+1}' < c_{2n+2}' < \mu_{n+2}' < \cdots < a \,\,\, \text{if} \,\,\, d < b < a < v
$$
and, if $v \le b < a \le z_0$, then 
\begin{multline*}
$$
b < \mu_{n}' < \mu_{n+1}' < \mu_{n+2}' < \mu_{n+3}' < \cdots < a, \\ 
b < c_{2n}' < c_{2n+2}' < c_{2n+4}' < c_{2n+6}' < \cdots < a \,\,\, \text{and} \\ 
\mu_{n+i}' > c_{2n+2i}' \,\,\, \text{for all} \,\,\, i \ge 0.
$$
\end{multline*}

\item[{\rm (3)}]
for each $i \ge 0$ such that $n+i$ is even and $n+i \ge 2k$, each of the points $c_{2n+2i}$ and $c_{2n+2i}'$ is period-$(2n+2i)$ point of $f$.  

\item[{\rm (4)}]
for each $i \ge 0$ such that $n+i$ is odd and $n+i \ge 2k$, each of the points $c_{2n+2i}$ and $c_{2n+2i}'$ is either a period-$(2n+2i)$ point of $f$ or a period-$(n+i)$ point of $f$.  

\item[{\rm (5)}]
if both $a$ and $b$ are points in one of the intervals $[\min P, d]$ and $[d, v]$, then, for each $i \ge 0$ such that $n+i$ is odd and $n+i \ge \max \big\{ 2k, m+2 \big\}$, each of the points $c_{2n+2i}$ and $c_{2n+2i}'$ is a period-$(2n+2i)$ points of $f$ while, if both $a$ and $b$ are in the interval $[v, z_0]$, then, for each $i \ge 0$ such that $n+i$ is odd and $n+i \ge m+2k$, each of the points $c_{2n+2i}$ and $c_{2n+2i}'$ is a period-$(2n+2i)$ point of $f$.
\end{itemize}}

\noindent
{\bf Remark 2.}
In the above lemma, suppose, say, $a < b$.  Since $f^{2k}(a) \in \big\{ z_0, f(z_0), f(d) \big\}$ and $f^{2n}(b) = \min P$, we have $f^{2k+1}(a) \in \big\{ z_0, f(z_0) \big\}$ and $f^{m+2n}(b) = f^m(\min P) = \min P$.  By Lemma 5, we obtain that, for each $i \ge 0$ such that $m+2n+2i \ge 3(2k+1)$, the point 
$$
p_{m+2n+2i}^* = \min \big\{ a \le x \le b: f^{m+2n+2i}(x) = x \big\}
$$
exists and is a period-$(m+2n+2i)$ point of $f$.  These periodic points $p_{m+2n+2i}^*$'s of odd periods are {\it interspersed} with the periodic points of even periods obtained in the lemma.  To make things simple, we choose only those periodic points of $f$ of even periods obtained in the lemma (because we apply even iterates of $f$ to the points $a$ and $b$) to build the {\it basic} tower of periodic points of $f$ associated with $P$ and ignore all these periodic points $p_{m+2n+2i}^*$'s of odd periods.  

\noindent
{\it Proof.}
The following proof is similar to that of Lemma 5.  Assume that $\min P < a < b < d$ or $d \le a < b \le v$.  Since $f^{2k}(a) \in \big\{ z_0, f(z_0), f(d) \big\}$ and $2n \ge 2k$, we have $f^{2n}(a) \in \big\{ z_0, f(z_0), f(d) \big\}$.
{\large 
$$
\text{Suppose} \,\,\, \min P < a < b < d.
$$}Then $\min P < x < d$ for all $x$ in $[a, b]$.  Since $f^{2n}(a) - a \ge z_0 - a > 0$ and $f^{2n}(b) - b = \min P - b < 0$, the point $c_{2n} = \min \big\{ a \le x \le b: f^{2n}(x) = x \big\}$ exists.  Since $f^{2k}(a) \in \big\{ z_0, f(z_0), f(d) \big\}$, $2n \ge 2k$ and $f^{2n}(c_{2n}) = c_{2n}$, it follows from Lemma 3 that, for each $i \ge 1$, the point $c_{2n+2i} = \min \big\{ a \le x \le b: f^{2n+2i}(x) = x \big\}$ exists.

Furthermore, for each $i \ge 0$, since $f^{2n+2i}(a) - d \ge z_0 - d > 0$ and $f^{2n+2i}(c_{2n+2i}) - d = c_{2n+2i} - d < 0$, the point $\mu_{n+i} = \min \big\{ a \le x \le c_{2n+2i}: f^{2n+2i}(x) = d \big\}$ exists and is $< c_{2n+2i}$.

On the other hand, since $f^{2n+2i}\big([a, \mu_{n+i}]\big) \supset [d, z_0] \supset \{ v \}$, there is a point $\nu_{n+i}$ in $[a, \mu_{n+i}]$ such that $f^{2n+2i}(\nu_{n+i}) = v$.  Since $f^{2n+2i+2}(a) - a \ge z_0 - a > 0$ and $f^{2n+2i+2}(\nu_{n+i}) - \nu_{n+i} = f^2(v) - \nu_{n+i} = \min P - \nu_{n+i} < 0$, the point 
$$
c_{2n+2i+2} = \min \big\{ a \le x \le \nu_{n+i}: f^{2n+2i+2}(x) = x \big\} \,\,\, \text{exists and is} \,\, < \nu_{n+i} < \mu_{n+i}.
$$
This, combined with the above, implies that 
$$
min P < a < \cdots < \mu_{n+2} < c_{2n+4} < \mu_{n+1} < c_{2n+2} < \mu_{n} < c_{2n} < b < d.
$$  
\vspace{.05in}
\[
\left( \begin{array}{l}
\textrm{\qquad\qquad Suppose $d \le a < b \le v$.  Then $d \le x \le v$ for all $x$ in $[a, b]$.}\\
\textrm{} \\
\textrm{Since $f^{2n}\big([a, b]\big) \supset \{ d \}$, the point $\mu_{n} = \min \big\{ a \le x \le b: f^{2n}(x) = d \big\}$ exists.}\\
\textrm{Since $f^{2n}(a) - a \ge z_0 - a > 0$ and $f^{2n}(\mu_{n}) - \mu_{n} = d - \mu_{n} < 0$, the point}\\
\textrm{\qquad\qquad $c_{2n} = \min \big\{ a \le x \le \mu_{n}: f^{2n}(x) = x \big\}$ exists and is $< \mu_{n}$.}\\
\textrm{Since $f^{2n}\big([a, c_{2n}]\big) \supset [c_{2n}, z_0] \supset \{ v \}$, there is a point $\nu_{n}$ in $[a, c_{2n}]$ such that}\\
\textrm{\qquad\qquad $f^{2n}(\nu_{n}) = v$.  Since $f^{2n+2}(a) - d \ge z_0 - d > 0$ and}\\
\textrm{\qquad $f^{2n+2}(\nu_{n}) - d = f^2(\nu_{n}) - d = \min P - d < 0$, the point}\\
\textrm{\qquad $\mu_{n+1} = \min \big\{ a \le x \le \nu_{n}: f^{2n+2}(x) = x \big\}$ exists and is $< \nu_{n} < c_{2n}$.}\\
\textrm{Since $f^{2n+2}(a) - a \ge z_0 - a > 0$ and $f^{2n+2}(\mu_{n+1}) - \mu_{n+1} = d - \mu_{n+1} < 0$,}\\
\textrm{the point $c_{2n+2} = \min \big\{ a \le x \le \mu_{n+1}: f^{2n+2}(x) = x \big\}$ exists and is $< \mu_{n+1}$.}\\
\textrm{By induction, we have}\\
\textrm{$d \le a < \cdots < c_{2n+4} < \mu_{n+2} < c_{2n+2} < \mu_{n+1} < c_{2n} < \mu_{n} < b \le v$.}
\end{array} \right)  
\]

\indent Note that we have assumed that $\min P < a < b < d$ or $d \le a < b \le v$.  Now, we want to find the least periods of $c_{2n+2i}$'s with respect to $f$.  For each $i \ge 0$, let $\ell(c_{2n+2i})$ denote the least period of $c_{2n+2i}$ with respect to $f$.  
$$
\text{Suppose} \,\,\, i \ge 0 \,\,\, \text{and} \,\,\, n+i \ge 2k \,\,\, \text{and let} \,\,\, r = (2n+2i)/\ell(c_{2n+2i}).
$$
If $r \ge 3$ is odd, then $r-1 \ge 2$ is even.  Since $(r-1)\ell(c_{2n+2i}) = (r-1)[(2n+2i)/r] \ge 2[(r-1)/r] \cdot 2k \ge 2k$, we can apply Lemma 3 with 
$$
f^{2k}(a) \in \big\{ z_0, f(z_0), f(d) \big\} \,\,\, \text{and} \,\,\, f^{(r-1)\ell(c_{2n+2i})}(c_{2n+2i}) = c_{2n+2i}.
$$
If $r \ge 4$ is even, then $r-2 \ge 2$ is even.  Since $(r-2)\ell(c_{2n+2i}) = (r-2)[(2n+2i)/r] \ge 2[(r-2)/r] \cdot 2k \ge 2k$, we can apply Lemma 3 with 
$$
f^{2k}(a) \in \big\{ z_0, f(z_0), f(d) \big\} \,\,\, \text{and} \,\,\, f^{(r-2)\ell(c_{2n+2i})}(c_{2n+2i}) = c_{2n+2i}.
$$
In either case, we obtain a periodic point $w_{2n+2i}$ of $f$ in $[a, c_{2n+2i})$ such that $f^{2n+2i}(w_{2n+2i}) = w_{2n+2i}$.  This contradicts the minimality of $c_{2n+2i}$ in $[a, b]$.  So, if $r > 2$, then $\ell(c_{2n+2i}) = 2n+2i$.

If $r = 2$ and $\ell(c_{2n+2i}) \, (= n+i \ge 2k)$ is even, then by applying Lemma 3 with 
$$
f^{2k}(a) \in \big\{ z_0, f(z_0), f(d) \big\} \,\,\, \text{and} \,\,\, f^{\ell(c_{2n+2i})}(c_{2n+2i}) = c_{2n+2i}.,
$$
we obtain the same contradiction.  So, in this case, $\ell(c_{2n+2i}) = 2n+2i$.

\noindent
Consequently, if $i \ge 0$ and $n+i$ is an even integer such that $n+i \ge 2k$, then since we have either $r > 2$ or $r = 2$, it follows from the above that $c_{2n+2i}$ is a period-$(2n+2i)$ point of $f$.  

On the other hand, suppose $i \ge 0$ and $n+i$ is an odd integer such that $n+i \ge 2k$.  If $r > 2$, then it follows form the above that $c_{2n+2i}$ is a period-$(2n+2i)$ point of $f$.  If $r = 2$, then $\ell(x_{2n+2i}) = n+i$.  So, if $i \ge 0$ and $n+i$ is an odd integer such that $n+i \ge 2k$, then $c_{2n+2i}$ is either a period-$(2n+2i)$ point of $f$ or a period-$(n+i)$ point of $f$.
$$.............. \qquad .............. \qquad .............. \qquad ..............$$
\indent Suppose $i \ge 0$ and, $n+i$ is odd and $\ge 2k$ and $c_{2n+2k}$ is a period-$(n+i)$ point of $f$.  Since $n+i \ge 2k$ and $f^{2k}(a) \in \big\{ z_0, f(z_0), f(d) \big\}$, we have $f^{n+i}(a) \in \big\{ z_0, f(z_0), f(d) \big\}$.  Since we have assumed that $\min P < a < b < v$ or $d \le a < b \le v$, we have 
$$
f^{n+i}\big([a, c_{2n+2i}]\big) \supset \big[f^{n+i}(c_{2n+2i}), f^{n+i}(a)\big] \supset [c_{2n+2i}, z_0] \supset [v, z_0].
$$
Consequently, we have \,\,\, $f^{n+i+2}\big([a, c_{2n+2i}]\big) \supset f^2\big([v, z_0]\big) \supset [\min P, z_0]$.  Since $f^2\big([\min P, z_0]\big) \supset f^2\big([v, z_0]\big) \supset [\min P, z_0]$, we obtain that, 
$$
\text{for each} \,\,\, i \ge 1 \,\,\, \text{such that} \,\,\, n+i \,\,\, \text{is odd and} \,\, \ge \max \big\{ m+2, 2k \big\},
$$ 
$n+i-m-2$ is even and $\ge 0$ and $f^{2n+2i}\big([a, c_{2n+2i}]\big) \supset f^{n+i-2}\big([\min P, z_0]\big) \supset \\ f^m\big(f^{n+i-m-2}([\min P, z_0])\big) \supset f^m\big([\min P, z_0]\big) \supset [\min P, z_0]$.  

Let $\delta_{2n+2i}$ be a point in $[a, c_{2n+2i}]$ such that $f^{2n+2i}(\delta_{2n+2i}) = \min P$.  It is clear that $a < \delta_{2n+2i} < c_{2n+2i} < b$.  Since $f^{2n+2i}(a) - a \ge z_0 -a > 0$ and $f^{2n+2i}(\delta_{2n+2i}) - \delta_{2n+2i} = \min P - \delta_{2n+2i} < 0$, we see that there exists a point $w_{2n+2i}$ in $[a, \delta_{2n+2i}]$ such that $f^{2n+2i}(w_{2n+2i}) = w_{2n+2i}$.  Since $a < w_{2n+2i} < \delta_{2n+2i} < c_{2n+2i}$ and $f^{2n+2i}(w_{2n+2i}) = w_{2n+2i}$, this contradicts the minimality of $c_{2n+2i}$ in $[a, b]$.  Therefore, the point $c_{2n+2i}$ can not be a period-$(n+i)$ point of $f$ and must be a period-$(2n+2i)$ point of $f$.

In summary, we have shown that, for each $i \ge 1$,  
\begin{multline*}
$$
\text{if} \,\,\, n+i \,\,\, \text{is even and} \, \ge 2k, \, \text{then the point} \,\,\, c_{2n+2i} \,\,\, \text{is a period-}(2n+2i) \,\,\, \text{point of} \,\,\, f; \vspace{.25in} \\ 
\text{if} \,\,\, n+i \,\,\, \text{is odd and} \, \ge \max \big\{2k, m+2 \big\}, \, \text{then the point} \,\,\, c_{2n+2i} \,\,\, \text{is a period-}(2n+2i) \,\,\, \text{point of} \,\,\,f.
$$
\end{multline*}
$$........................................................................$$
{\large 
$$
\text{Now suppose} \,\,\, v < a < b < z_0.
$$}In this case, let $\mu_{n} = \min \big\{ a \le x \le b : f^{2n}(x) = d \big\}$.  Then $f^{2n}\big([a, \mu_{n}]\big) \supset [d, z_0] \supset \{ v \}$.  Let $\nu_{n}$ be a point in $[a, \mu_{n}]$ such that $f^{2n}(\nu_{n}) = v$.  Since $f^{2n+2}\big([a, \nu_{n}]\big) \supset [\min P, z_0] \supset \{ d \}$, the point 
$$
\mu_{n+1} = \min \big\{ a \le x \le \nu_{n} : f^{2n+2}(x) = d \big\} = \min \big\{ a \le x \le b : f^{2n+2}(x) = d \big\}
$$
exists and is $< \nu_n < \mu_n$.  Inductively, for each $i \ge 2$, the point 
$$
\mu_{n+i} = \min \big\{ a \le x \le \mu_{n+i-1} : f^{2n+2i}(x) = d \big\} = \min \big\{ a \le x < b : f^{2n+2i}(x) = d \big\}
$$
exists and \qquad $v < a < \cdots < \mu_{n+3} < \mu_{n+2} < \mu_{n+1} < \mu_{n} < b < z_0$.  

On the other hand, since $f^{2n}(a) - a \ge z_0 - a > 0$ and $f^{2n}(b) - b = \min P - b < 0$, the point 
$$
c_{2n} = \min \big\{ a \le x \le b : f^{2n}(x) = x \big\}
$$
exists.  Therefore, we can apply Lemma 3 with 
$$
f^{2k}(a) \in \big\{ z_0, f(z_0), f(d) \big\} \,\,\, \text{and} \,\,\, f^{2n}(c_{2n}) = c_{2n}
$$
to obtain that, for each $i \ge 0$, the point 
$$
c_{2n+2i} = \min \big\{ a \le x \le b : f^{2n+2i}(x) = x \big\}
$$
exists and \qquad $v < a < \cdots < c_{2n+6} < c_{2n+4} < c_{2n+2} < c_{2n} < b < z_0$.  

Furthermore, for each $i \ge 0$, since 
$$
f^{2n+2i}(a) - a \ge z_0 - a > 0 \,\,\, \text{and} \,\,\, f^{2n+2i}(\mu_{n+i}) - \mu_{n+i} = d - \mu_{n+i} < 0,
$$
we have $c_{2n+2i} = \min \big\{ a \le x \le b: f^{2n+2i}(x) = x \big\} =$ {\small $\min \big\{ a \le x \le \mu_{n+i}: f^{2n+2i}(x) = x \big\}$} $< \mu_{n+i}$.  
$$.............. \qquad .............. \qquad .............. \qquad ..............$$
\indent Now we find the least periods of $c_{2n+2i}$'s with respect to $f$.  By following the arguments as those on the interval $[\min P, d]$ above, we can obtain that, for each $i \ge 0$,

\noindent
if $n+i$ is even and $\ge 2k$, then $c_{2n+2i}$ is a period-$(2n+2i)$ point of $f$,

\noindent
if $n+i$ is odd and $\ge 2k$, then $c_{2n+2i}$ is either a period-$(2n+2i)$ or a period-$(n+i)$ point of $f$.
$$
\text{Suppose} \,\,\, i \ge 0, \, \text{and} \,\,\, n+i \,\,\, \text{is odd and} \, \ge m+2k \,\,\, \text{and} \,\,\, c_{2n+2i} \,\,\, \text{is a period-}(n+i) \,\,\, \text{point of} \,\,\, f.
$$
If the points $f^2(c_{2n+2i})$, $f^4(c_{2n+2i})$, $\cdots$, $f^{n+i-1}(c_{2n+2i})$ are all $> v$, then, by the fact that $f^2(x) < x$ for all $v \le x < z_0$, we have 
$$
d < v < f^{n+i-1}(c_{2n+2i}) < \cdots < f^4(c_{2n+2i}) < f^2(c_{2n+2i}) < c_{2n+2i} < z_0.
$$
By the fact that $f(x) > z \ge z_0$ for all $d < x < z_0$, we obtain that $c_{2n+2i} = f\big(f^{n+i-1}(c_{2n+2i})\big) > z_0 > c_{2n+2i}$ which is a contradiction.  Therefore, there is an integer $s$ such that $1 \le s \le (n+i-1)/2$ and $f^{2s}(c_{2n+2i}) \le v$.  Note that $f^2\big([\min P, z_0]\big) \supset f^2\big([v, z_0]\big) \supset [\min P, z_0]$.  We have two cases to consider:

\noindent
(i) If $2s < 2k$, and $n+i$ is odd and $\ge m+2k$ $(\ge m+2s+2)$, then \vspace{.03in} \\ \indent $n+i-m$ $(\ge 2k \ge 2s+2)$ is even and $$
f^{n+i+m+2s+2}\big([a, c_{2n+2i}]\big) \supset f^{m+2s+2}\big([c_{2n+2i}, z_0]\big) \supset f^{m+2}\big([v, z_0]\big) \supset [\min P, z_0]. \,\,\, \text{So},
$$
\indent $f^{2n+2i}\big([a, c_{2n+2i}]\big) \supset f^{n+i-m-2s-2}\big([\min P, z_0]\big) \supset [\min P, z_0]$.  

\noindent
(ii) If $2k \le 2s \,(\le n+i-1)$, then \vspace{.03in} \\ \indent $2n+2i - 2s -2 =(n+i+1)+[(n+i-1)-2s]-2 \ge n+i-1 \ge 2s \ge 2k > 0$ and 
$$
f^{2s+2}\big([a, c_{2n+2i}]\big) \supset f^{2}\big([f^{2s}(c_{2n+2i}), z_0]\big) \supset f^{2}\big([v, z_0]\big) \supset [\min P, z_0]. \,\,\, \text{So,}
$$ 
\indent $f^{2n+2i}\big([a, c_{2n+2i}]\big) \supset f^{2n+2i-2s-2}\big(f^{2s+2}([a, c_{2n+2i}])\big) \supset f^{2n+2i-2s-2}\big([\min P, z_0]\big) \supset [\min P, z_0]$.

By combining the above two results, we obtain that, for each $i \ge 0$ such that $n+i$ is odd and $\ge m+2k$, there exists a point $\delta_{2n+2i}$ in $[a, c_{2n+2i}]$ such that $f^{2n+2i}(\delta_{2n+2i}) = \min P$.  Since 
$$
f^{2n+2i}(a) - a \ge z_0 - a > 0 \,\,\, \text{and} \,\,\, f^{2n+2i}(\delta_{2n+2i}) - \delta_{2n+2i} = \min P - \delta_{2n+2i} < 0,
$$
there is a point $w_{2n+2i}$ in $[a, \delta_{2n+2i}]$ such that $f^{2n+2i}(w_{2n+2i}) = w_{2n+2i}$.  Since $a < w_{2n+2i} < \delta_{2n+2i} < c_{2n+2i} < b$, this contradicts the minimality of $c_{2n+2i}$ in $[a, b]$.  Therefore, for each $i \ge 0$ such that $n+i$ is odd and $\ge m+2k$, the point $c_{2n+2i}$ is a period-$(2n+2i)$ point of $f$.  Since we have known that, for each $i \ge 0$ such that $n+i$ is even and $\ge 2k$, the point $c_{2n+2i}$ is a period-$(2n+2i)$ point of $f$, this proves (1).  (2) can be proved similarly.  
$\hfill\square$
$$\aleph \qquad\qquad\qquad \aleph \qquad\qquad\qquad \aleph \qquad\qquad\qquad \aleph \qquad\qquad\qquad \aleph$$
\indent We now divide the interval $[\min P, z_0]$ into the 5 subintervals: $[\min P, d]$, $[d, u_1]$, $[u_1, v]$, $[v, \bar u_1']$ and $[\bar u_1', z_0]$.  On each such subinterval, we can apply Lemmas 5 $\&$ 6 in various combinations to build various towers of periodic points of $f$ associated with $P$.  As an example, we shall apply Lemmas 5 $\&$ 6 in a recursive way to build the so-called {\it basic} tower of periodic points of $f$ associated with $P$.  

In the sequel, for the sake of clarity, we shall use the following notations: 
\begin{multline*}
$$
\, x' \,\,\, \text{to denote points obtained by taking the maximum of a set};\, \\ \tilde x, x, \breve x, \hat x \,\,\, \text{and} \,\,\, \bar x \,\,\, \text{respectively to denote points (except endpoints) in the intervals}\quad\qquad\,\,\, \\ 
\, [\min P, d], \, [d, u_1], \, [u_1, v], \, [v, \hat u_0] \, (\subset [v, \bar u_1']), \, [\bar u_1', z_0] \,\,\, \text{respectively}; \\ 
\,\,\,\,\, \mu_{m,\cdots} \,\,\, \text{and} \,\,\, u_{\cdots} \,\,\, \text{respectively to denote} \,\,\, d\text{-points in}\qquad\qquad\qquad\qquad\qquad\qquad\qquad\qquad\qquad\quad\quad\quad\,\,\, \\ 
[\min P, d) \cup (u_1, v] \cup [v, \bar u_1') \,\,\, \text{and} \,\,\, (d, u_1] \cup [\bar u_1', z_0] \,\,\, \text{respectively}; \\ 
p \,\,\, \text{and} \,\,\, q \,\,\, \text{respectively to denote periodic points of {\it odd} periods of} \,\, f \qquad\qquad\qquad\qquad\quad\,\,\,\,\, \\ 
\text{obtained by taking the minimum and maximum of a set respectively}; \\ 
c \,\,\, \text{to denote a periodic point of {\it even} period of} \,\,\, f \, \qquad\qquad\qquad\qquad\qquad\qquad\qquad\qquad\,\,\,\,\,\,\,\,\,\, \\ 
\text{obtained by taking the minimum of a set}.\qquad\qquad\qquad\qquad\qquad\qquad\qquad\qquad\,\,
$$
\end{multline*}

\noindent
{\bf $\mathsection$1. The first layer of the {\it basic} tower of periodic points of $f$ associated with $P$.}

To start with, we shall find a monotonic sequence of periodic points of $f$ (and a monotonic sequence of $d$-points, a $d$-point is a point whose orbit under $f$ contains the point $d$) on each of the 5 intervals $[\min P, d]$, $[d, u_1]$, $[\breve u_0', v]$ $(\subset [u_1, v])$, $[v, \hat u_0]$ $(\subset [v, \bar u_1'])$ and $[\bar u_1', z_0]$ which comprise the first layer of the {\it basic} tower of periodic points of $f$ associated with $P$.  Then on each of the 3 intervals: $[\min P, d], [\breve u_0', v]$ $(\subset [u_1, v])$ and $[v, \hat u_0]$ $(\subset [v, \bar u_1'])$, where we apply odd iterates of $f$ to the endpoints, we shall apply Lemma 5 once and Lemma 6 twice successively, and on each of the 2 intervals: $[d, u_1]$ and $[\bar u_1', z_0]$, where we apply even iterates of $f$ to the endpoints, we shall apply Lemma 6 once and Lemma 5 twice successively, to continue the {\it basic} tower-building process to the second and higher layers.  

\vspace{.1in}

\noindent
{\bf 1.1 On the existence of periodic points of $f$ of all odd periods $\ge m$ in $[\min P, d]$.}

In this case, we clearly have 
{\large 
$$
f(d) \in \big\{ z_0, f(z_0), f(d) \big\} \,\,\,  \text{and} \,\,\, f^m(\min P) = \min P
$$}
and trivially 
$$
\text{all periodic points of} \,\,\, f \,\,\, \text{in} \,\,\, [\min P, d] \,\,\, \text{with {\it odd} periods have least periods} \,\, \ge 1.
$$
Therefore, we can apply Lemma 5(2) to obtain that, for each $i \ge 0$, the points
$$
\, \tilde q_{m+2i} = \max \big\{ \min P \le x \le d : f^{m+2i}(x) = x \big\}
$$
and
$$
\tilde \mu_{m,i}' = \max \big\{ \min P \le x \le d : f^{m+2i}(x) = d \big\} \,\,\,
$$
exist and $\tilde q_{m+2i}$ is a period-$(m+2i)$ point of $f$.  Furthermore, we have 
$$
\min P \, \le \, \tilde q_m \, < \, \tilde \mu_{m,0}' \, < \, \tilde q_{m+2} \, < \, \tilde \mu_{m,1}' \, < \, \tilde q_{m+4} \, < \, \tilde \mu_{m,2}' \, < \, \tilde q_{m+6} \, < \, \tilde \mu_{m,3}' \, < \, \cdots \, < \, d.
$$
Note the relationship between the subscript of $\tilde \mu_{m,i}'$ and the superscript of $f^{m+2i}$ in the definition of $\tilde \mu_{m,i}'$.  We need these points $\tilde \mu_{m,i}'$'s later to continue the {\it basic} tower-building process to the next (second) layer. 

\noindent
{\bf Remark 3.}
Note that we can apply Lemma 5(2) with 
{\large 
$$
f(d) \in \big\{ z_0, f(z_0), f(d) \big\} \,\,\,  \text{and} \,\,\, f^m(\min P) = \min P
$$}
to obtain that, for each $i \ge 1$, the point $\tilde c_{2m+2i}'^{\,*} = \max \big\{ \min P \le x \le d : f^{2m+2i}(x) = x \big\}$ exists and is a period-$(2m+2i)$ point of $f$ in $[\min P, d]$ and 
$$
\min P \, < \, \tilde c_{2m+2}'^{\,*} < \, \tilde c_{2m+4}'^{\,*} < \, \tilde c_{2m+6}'^{\,*} < \, \tilde c_{2m+8}'^{\,*} < \, \cdots < \, d.
$$
As for the point $\tilde c_{2m}'^{\,*} = \max \big\{ \min P \le x \le d : f^{2m}(x) = x \big\}$, it may be a period-$m$ point of $f$.  
\vspace{.05in}
\[
\left( \begin{array}{l}
\textrm{For example, let $g : [0, 1] \longrightarrow [0, 1]$ be the continuous map defined by putting }\\
\textrm{(i) $g(x) = x + 1/2$ and (ii) $g(x) = 2 - 2x$ and let $P = \{ 0, 1/2, 1 \}$ be the unique}\\
\textrm{period-3 orbit of $g$.  Then $\min P = 0, \, d = 1/6, \, v = 1/2, \, z = 2/3$ and $g(v) = 1$.}\\
\textrm{$g$ has exactly two period-6 orbits and they all lie in the interval $[1/6, 1]$.}\\ 
\textrm{In this case, $m = 3$ and $\tilde c_6'^{\,*}$ is a period-3, but not a period-6 point of $g$.}  
\end{array} \right)  
\]
Note that, for each $i \ge 1$, the orbits $O_f(\tilde c_{2m+2i}'^{\,*})$'s of the periodic points $\tilde c_{2m+2i}'^{\,*}$'s are disjoint from the orbits $O_f(c_{2i})$'s of $c_{2i}$'s (defined in the main text) and the orbits $O_f(\bar c_{2i})$'s of $\bar c_{2i}$'s (defined below) because all $O_f(c_{2i})$'s and $O_f(\bar c_{2i})$'s are contained in the interval $(d, z_0)$ which is disjoint from the interval $[\min P, d]$.  However, these periodic points $\tilde c_{2m+2i}'^{\,*}$'s of $f$ of even periods are {\it interspersed} with those periodic points $\tilde q_{m+2i}$'s of $f$ of odd periods and, to make things simple, are not counted in the first layer of the {\it basic} tower of periodic points of $f$ associated with $P$.    

\noindent
{\bf Remark 4.}
On the other hand, since $f^{2m}(\min P) = \min P$ and $f^{2m}(\tilde \mu_{m,0}') = f^m(d) = f(z_0) \ge z_0$, the point $\tilde \xi_{2m}' = \max \big\{ \min P \le x \le \tilde \mu_{m,0}': f^{2m}(x) = d \big\}$ exists and $\min P < \tilde \xi_{2m}' < \tilde \mu_{m,0}' < d$.  By applying Lemma 6(2) with
$$
f^{m+1}(\tilde \mu_{m,0}') = f(d) \in \big\{ z_0, f(z_0), f(d) \big\} \,\,\, \text{and} \,\,\, f^{2m}(\min P) = \min P
$$
and arguments similar to Remark 2, we obtain infinitely many more periodic points of $f$ with even periods $> 2m$ and with odd periods $> 3m$ in the interval $[\tilde \xi_{2m}', \tilde \mu_{m,0}']$ $(\subset [\min P, \tilde \mu_{m,0}'])$.  

Since $f^{3m}(\min P) = \min P$ and $f^{2m+1}(\tilde \xi_{2m}') = f\big(f^{2m}(\tilde \xi_{2m}')\big) = f(d)$, the point $\tilde \xi_{3m}' = \max \big\{ \min P \le x \le \tilde \xi_{2m}': f^{3m}(x) = d \big\}$ exists and $\min P < \tilde \xi_{3m}' < \tilde \xi_{2m}' < \tilde \mu_{m,0}' < d$.  By applying Lemma 5(2) with
$$
f^{2m+1}(\tilde \xi_{2m}') = f(d) \,\,\, \text{and} \,\,\, f^{3m}(\min P) = \min P,
$$
we obtain infinitely many more periodic points of $f$ with odd periods $> 3m$ and with even periods $> 4m$ in the interval $[\tilde \xi_{3m}', \tilde \xi_{2m}']$ $(\subset [\min P, \tilde \xi_{2m}'])$.

Inductively, for each $k \ge 2$, the point $\tilde \xi_{km}' = \max \big\{ \min P \le x \le \tilde \xi_{(k-1)m}': f^{km}(x) = d \big\}$ exists and $$
\min P < \cdots < \tilde \xi_{km}' < \cdots < \tilde \xi_{3m}' < \tilde \xi_{2m}' < \tilde \mu_{m,0}' < d.
$$
By applying Lemma 5(2) or, Lemma 6(2) and Remark 2 appropriately, we obtain that, for each $k \ge 2$, there are infinitely many more periodic points of $f$ with even periods and with odd periods in the interval $[\tilde \xi_{km}', \tilde \xi_{(k-1)m}']$ $(\subset [\min P, \tilde \xi_{(k-1)m}'])$.

{\large 
However, to make things simple, we choose to {\it ignore} these periodic points of $f$ obtained in the interval $[\min P, \tilde \mu_{m,0}']$ as described in Remarks 3 $\&$ 4 above and keep only the periodic points $\tilde q_{m+2i}$'s of $f$ with odd periods in the interval $[\tilde \mu_{m,0}', d]$ in the first layer and use the points $\tilde \mu_{m,i}'$'s to continue the process into the next layer of the {\it basic} tower of periodic points of $f$ associated with $P$.} 

\vspace{.1in}

\noindent
{\bf 1.2 On the existence of periodic points of $f$ of all even periods $\ge 2$ in $[d, u_1]$ $(\subset [d, v])$.}

This is done in Part (c) in the main text.  We can also argue as follows: Since we have two points $d < v$ such that $f^2(d) = z_0$ and $f^2(v) = \min P$ so that Lemma 6(1) can be applied to obtain that, for each $n \ge 1$, the points 
$$
u_n = \min \{ d \le x \le v: f^{2n}(x) = d \} \,\,\, \text{and} \,\,\, c_{2n} = \min \big\{ d \le x \le v: f^{2n}(x) = x \big\}
$$
exist and $d < \cdots < c_6 < u_3 < c_4 < u_2 < c_2 < u_1 < v$.  Note the relationship between the subscript of $u_n$ and the superscript of $f^{2n}$ in the definition of $u_n$.  We need these points $u_n$'s later to continue the {\it basic} tower-building process to the next (second) layer.  

By Lemma 4(1), we obtain that, for each $n \ge 1$, the point 
$$
c_{2n} = \min \big\{ d \le x \le v: f^{2n}(x) = x \big\}
$$
is a period-$(2n)$ point of $f$ in $[d, u_1]$ $(\subset [d, v])$.

\noindent
{\bf Remark 5.}
On the interval $[d, u_1]$, let $\nu$ be a point such that $f^2(\nu) = v$.  Then we have 
$$
f^1(d) = f(d) \,\,\, \text{and} \,\,\, f^{m+2}(\nu) = f^m(v) = \min P.
$$
On the other hand, trivially, all periodic points of $f$ in $[d, z_0]$ of odd periods have least periods $\ge 1$.  It follows from Lemma 5(1) that, for each $i \ge 1$, the point 
$$
p_{m+2i} = \min \big\{ d \le x \le v: f^{m+2i}(x) = x \big\}
$$
exists and is a period-$(m+2i)$ point of $f$ and $d < \cdots < p_{m+8} < p_{m+6} < p_{m+4} < p_{m+2} < \nu < u_1$.  These periodic points $p_{m+2i}$'s of odd periods are {\it interspersed} with the periodic points $c_{2i}$'s of $f$ of even periods and, to make things simple, are not counted in the first layer of the {\it basic} tower of periodic points of $f$ associated with $P$

\vspace{.1in}

\noindent
{\bf 1.3 On the existence of periodic points of $f$ of all odd periods $\ge m+2$ in $[\breve u_0', v]$ $([u_1, v]$ $\subset [d, v])$.}

We have shown this result in the main text.  Here we apply Lemma 5(1) to obtain the same result.  On the interval $[u_1, v]$, we consider the point $\breve u_0' = \max \big\{ u_1 \le x \le v: f^2(x) = d \big\}$ instead of the point $u_1$ (for a reason which will become apparant below).  Then, we have
$$
f^3(\breve u_0') = f(d) \in \big\{ z_0, f(z_0), f(d) \big\} \,\,\, \text{and} \,\,\, f^{m+2}(v) = f^m(f^2(v)) = f^m(\min P) = \min P.
$$
On the other hand, it follows from the fact 
$$
\qquad\quad\qquad\qquad f(x) > z \,\,\, \big(\ge z_0 > f^2(x)\big) \,\,\, \text{for all} \,\,\, d < x < z_0, \qquad\qquad\qquad\qquad\qquad\quad\,\, (***)
$$
that $f$ has no fixed points in $(\breve u_0', v)$. So, 
$$
\text{all periodic points of} \,\,\, f \,\,\, \text{in} \,\,\, [\breve u_0', v] \,\,\, \text{with {\it odd} periods have least periods} \,\, \ge 3.
$$
By applying Lemma 5(1) with
{\large 
$$
f^3(\breve u_0') = f(d) \in \big\{ z_0, f(z_0), f(d) \big\} \,\,\, \text{and} \,\,\, f^{m+2}(v) = \min P,
$$}
we obtain that, for each $i \ge 0$, the points
$$
\qquad\quad \breve p_{m+2+2i} = \min \big\{ \breve u_0' \le x \le v : f^{m+2+2i}(x) = x \big\} \,\,\, \text{and}
$$ 
$$\, \breve \mu_{m,1+i} = \min \big\{ \breve u_0' \le x \le v : f^{m+2+2i}(x) = d \big\}
$$
exist and $\breve p_{m+2+2i}$ is a period-$(m+2+2i)$ point of $f$.  Furthermore, it follows from the choice of $\breve u_0$ that 
$$
f^2(x) < d \,\,\, \text{for all} \,\,\, \breve u_0 < x < v.
$$
Consequently, we have 
$$
u_1 \, \le \, \breve u_0' < \, \cdots \, < \, \breve p_{m+6} \, < \, \breve \mu_{m,3} \, < \, \breve p_{m+4} \, < \, \breve \mu_{m,2} \, < \, \breve p_{m+2} \, < \, \breve \mu_{m,1} \, < \, v.
$$
Note the relationship between the subscript of $\breve \mu_{m,1+i}$ and the superscript of $f^{m+2+2i}$ in the definition of $\breve \mu_{m,1+i}$.  We need these points $\breve \mu_{m,1+i}$'s later to continue the {\it basic} tower-building process to the second layer.

\noindent
{\bf Remark 6.}  
Note that, in the interval $[\breve u_0, v]$ $(\subset [u_1, v])$, we have $f^3(\breve u_0') = f(d) \ge z_0$ and $f^{m+2}(v) = \min P$ \big(compare with $f^m(\min P)$ and $f(d)$ on $[\min P, d]$\big).  As discussed in Remarks 3 $\&$ 4 above, we can apply Lemma 5(1) or, Lemma 6(1) and arguments similar to Remark 2 appropriately on the interval $[\breve u_0, v]$ to obtain infinitely many more periodic points of $f$ of even periods and of odd periods in $[\breve u_0', v]$.  As before, to make things simple, we choose to ignore all of these periodic points and keep only those periodic points $\breve p_{m+2+2i}$'s of $f$ with odd periods (in the first layer of the {\it basic} tower) and the points $\breve \mu_{m,1+i}$'s to continue the process of building the {\it basic} tower of periodic points of $f$ associated with $P$ into the second layer.

\vspace{.1in}

\noindent
{\bf 1.4 On the existence of periodic points of $f$ of all odd periods $\ge m+2$ in $[v, \hat u_0]$ $(\subset [v, \bar u_1']$ $\subset [v, z_0])$.}

On the interval $[v, \bar u_1']$, we consider the point $\hat u_0 = \min \big\{ v \le x \le z_0 : f^2(x) = d \big\}$ $(\le \bar u_1')$ instead of the point $\bar u_1'$ (for a reason which will become apparant below).  Then, we have
$$
f^{3}(\hat u_0) = f(d) \,\,\, \text{and} \,\,\, f^{m+2}(v) = f^m\big(f^2(v)\big) = f^m(\min P) = \min P.
$$
On the other hand, it follows from the fact 
$$
\qquad\quad\qquad\qquad f(x) > z \,\,\, \big(\ge z_0 > f^2(x)\big) \,\,\, \text{for all} \,\,\, d < x < z_0, \qquad\qquad\qquad\qquad\qquad\quad\,\, (***)
$$
that $f$ has no fixed points in $(v, \hat u_0)$.  So, 
$$
\text{all periodic points of} \,\,\, f \,\,\, \text{in} \,\,\, [v, \hat u_0] \,\,\, \text{with {\it odd} periods have least periods} \,\, \ge 3.
$$
By applying Lemma 5(2) with 
{\large 
$$
f^3(\hat u_0) = f(d) \in \big\{ z_0, f(z_0), f(d) \big\}  \,\,\, \text{and} \,\,\, f^{m+2}(v) = \min P,
$$}
we obtain that, for each $i \ge 0$, the points
$$
\quad \hat q_{m+2+2i} = \max \big\{ v \le x \le \hat u_0 : f^{m+2+2i}(x) = x \big\}
$$
and
$$
\, \, \hat \mu_{m,1+i}' = \max \big\{ v \le x \le \hat u_0 : f^{m+2+2i}(x) = d \big\}
$$
exist and $\hat q_{m+2+2i}$ is a period-$(m+2+2i)$ point of $f$.  Furthermore, it follows from the choice of $\hat u_0$ that 
$$
f^2(x) < d \,\,\, \text{for all} \,\,\, v < x < \hat u_0.
$$
Consequently, 
we have 
$$
v \, < \, \hat \mu_{m,1}' \, < \, \hat q_{m+2} \, < \, \hat \mu_{m,2}' \, < \, \hat q_{m+4} \, < \, \hat \mu_{m,3}' \, < \hat q_{m+6} \, < \, \cdots \, < \, \hat u_0.
$$
Note the relationship between the subscript of $\hat \mu_{m,1+i}'$ and the superscript of $f^{m+2+2i}$ in the definition of $\hat \mu_{m,1+i}'$.  We need these points $\hat \mu_{m,1+i}'$'s later to continue the {\it basic} tower-building process to the second layer.

\noindent
{\bf Remark 7.}  Note that, on the interval $[v, \hat u_0]$, we have $f^{m+2}(v) = \min P$ and $f^3(\hat u_0) = f(d) \ge z \ge z_0$ \big(compare with $f^m(\min P)$ and $f(d)$ on $[\min P, d]$ and also with $f^3(\breve u_0')$ and $f^{m+2}(v)$ on $[\breve u_0', v]$ $(\subset [u_1, v])$\big).  As discussed in Remarks 3 $\&$ 4 above, by applying Lemmas 5(2) or, 6(2) and arguments similar to Remark 2 appropriately on the interval $[v, \hat u_0]$, we can obtain infinitely many more periodic points of $f$ of even periods and of odd periods.  However, to make things simple, we choose to ignore all of these periodic points and keep only those periodic points $\hat q_{m+2+2i}$'s of $f$ with odd periods (in the first layer of the {\it basic} tower) and the points $\hat \mu_{m,1+i}'$'s to continue the process of building the {\it basic} tower of periodic points of $f$ associated with $P$ into the next layer.

\vspace{.1in}

\noindent
{\bf 1.5 On the existence of periodic points of $f$ of all even periods $\ge 4$ in $[\bar u_1', z_0]$ $(\subset [v, z_0])$.}

In the main text, we have found points $(\min P <) \,\,\, d < \cdots < u_2 < u_1 < v$ such that, 
$$
\text{for each} \,\,\, n \ge 1, \,\, u_n = \min \big\{ d \le x \le v : f^{2n}(x) = d \big\}
$$
and have shown the existence of periodic points of $f$ of all {\it even} periods in the interval $[d, u_1]$ by showing that, for each $n \ge 1$, the point 
$$
c_{2n} = \min \big\{ d \le x \le v : f^{2n}(x) = x \big\}
$$
is a period-$(2n)$ point of $f$.  'Symmetrically', we define, as in Lemma 4, points $(v <) \,\,\, \bar u_1' < \bar u_2' < \cdots \,\,\, (< z_0)$ by putting, for each $n \ge 1$, 
$$
\bar u_n' = \max \big\{ v \le x \le z_0 : f^{2n}(x) = d \big\}.
$$
{\large 
However, to find new periodic points of $f$ of even periods in $[v, z_0]$, we cannot '{\it symmetrically}' define $\bar c_{2n}'$ as $$\bar c_{2n}' = \max \{ v \le x \le z_0 : f^{2n}(x) = x \}$$ because in this way we only get $z_0$ and no new periodic points.  So, we use a different strategy:}

For each $n \ge 1$, since $f^{2n+2}\big([\bar u_n', \bar u_{n+1}']\big) \supset [d, z_0] \supset \{ d \}$, the point $\bar \omega_{n+1}^* = \min \big\{ \bar u_n' \le x \le \bar u_{n+1}' : f^{2n+2}(x) = d \big\} \, (\ge \bar u_n')$ exists.  So, $d < f^{2n+2}(x) < z_0$ for all $\bar u_n' \le x \le \bar \omega_{n+1}^*$.  This, combined with the following fact
{\large 
$$
\qquad\,\,\, d < f^{2k}(x) < z_0 \,\,\, \text{for all} \,\,\, x \,\,\, \text{in} \,\, (d, u_n) \cup (\bar u_n', z_0) \,\, \text{and all} \,\,\, 1 \le k \le n, \quad\,\,\,\, (\dagger)
$$}
implies that $d < f^{2k}(x) < z_0$ for all $\bar u_n' \le x \le \bar \omega_{n+1}^*$ and all $1 \le k \le n+1$.  Since $f^{2n+2}(\bar u_n') - \bar u_n' = z_0 - \bar u_n' > 0$ and $f^{2n+2}(\bar \omega_{n+1}^*) - \bar \omega_{n+1}^* = d - \bar \omega_{n+1}^* < 0$, the point $$\bar c_{2n+2} = \min \big\{ \bar u_n' \le x \le \bar \omega_{n+1}^* : f^{2n+2}(x) = x \big\} = \min \big\{ \bar u_n' \le x \le \bar u_{n+1}' : f^{2n+2}(x) = x \big\}$$exists.  By Lemma 4(2), $\bar c_{2n+2}$ is a period-$(2n+2)$ point of $f$ such that, for all {\it odd} \, $1 \le j \le 2n+1$, 
$$
d < f^{2n}(\bar c_{2n+2}) < f^{2n-2}(\bar c_{2n+2}) < \cdots < f^2(\bar c_{2n+2}) < \bar c_{2n+2} < z_0 \le z < f^j(\bar c_{2n+2}).
$$
This establishes the existence of periodic points $\bar c_4, \bar c_6, \bar c_8, \cdots$ of $f$ of all even periods $\ge 4$ \big(best possible because $f$ may have exactly 2 period-2 points which form a period-2 orbit, one in $[d, v]$ and the other is $> z$\big) in $[\bar u_1', z_0]$ such that 
$$
v \,\, < \,\, \bar u_1' \,\, < \,\, \bar c_4 \,\, < \,\, \bar u_2' \,\, < \,\, \bar c_6 \,\, < \,\, \bar u_3' \,\, < \,\, \bar c_8 \,\, < \,\, \bar u_4' \,\, < \,\, \cdots \,\, < \,\, z_0.
$$
{\bf Remark 8.}  On the interval $[\bar u_1', z_0]$, let $\bar \nu$ to be a point such that $f^2(\bar \nu) = v$.  Then we have 
$$
f^{1}(z_0) = f(z_0) \in \big\{ z_0, f(z_0), f(d) \big\} \,\,\, \text{and} \,\,\, f^{m+2}(\bar \nu) = \min P.
$$
Since, trivially, all periodic points of $f$ in $[\bar \nu, z_0]$ of odd periods have least periods $\ge 1$, it follows from Lemma 5(2) that, for each $i \ge 0$, the point 
$$
\bar q_{m+2+2i} = \max \big\{ \bar u_1' \le x \le z_0: f^{m+2+2i}(x) = x \big\}
$$
exists and is a period-$(m+2+2i)$ point of $f$ and $\bar u_1' < \bar \nu < \bar q_{m+2} < \bar q_{m+4} < \bar q_{m+6} < \bar q_{m+8} < \cdots < z_0$.  These periodic points $\bar q_{m+2+2i}$'s of $f$ of odd periods are clearly {\it interspersed} with the periodic points $\bar c_{2i+2}$' of $f$ of even periods and, to make things simple, are not counted in the first layer of the {\it basic} tower of periodic points of $f$ associated with $P$.  

We can combined the results on the intervals $[d, u_1], [u_1, v]$ and $[v, \hat u_0]$ $(\subset [v, \bar u_1'])$ into one on the interval $[d, \hat u_0]$ $(\subset [d, \bar u_1'])$ as follows:

On the interval $[d, \hat u_0]$, we have $f^2(d) = z_0$, $f^2(v) = \min P$, $f^2(\hat u_0) = d$ and 
$$
\text{the graph of} \,\,\, y = f^2(x) \,\,\,  \text{looks 'roughly' like a {\it skew} 'V' shape at points} \,\,\, d, v, \hat u_0
$$
\big(meaning that the graph of $y = f^2(x)$ passes through the 3 points $(d, z_0), (v, \min P), (\hat u_0, d)$ on the $x$-$y$ plane with $y$-coordinates $z_0, \min P, d$ respectively\big).  {\it Note that} we choose the 3 points $d, v, \hat u_0$ instead of the 3 points $d, v, z_0$ for a reason which will become apparent later on.  

By combining the results on the intervals $[d, u_1]$, $[\breve u_0', v]$ $(\subset [u_1, v])$ and $[v, \hat u_0]$ $(\subset [v, \bar u_1'])$, we obtain that, with $x$-coordinates moving from point $d$ via point $v$ to point $\hat u_0$, for each $i \ge 1$,
\begin{itemize}
\item[{\rm (a)}]
(even periods) $c_{2i} = \min \big\{ d \le x \le v: f^{2i}(x) = x \big\}$ is a period-{\small $(2i)$} point of $f$;

\item[{\rm (b)}]
(odd periods) $\breve p_{m+2i} = \min \big\{ \breve u_0' \le x \le v: f^{m+2i}(x) = x \big\}$ is a period-{\small $(m+2i)$} point of $f$;

\item[{\rm (c)}] 
(odd periods) $\hat q_{m+2i} = \max \big\{ v \le x \le \hat u_0: f^{m+2i}(x) = x \big\}$ is a period-{\small $(m+2i)$} point of $f$.
\end{itemize}  
{\large This kind of {\it skew} 'V' shape phenomenon will be our basis for building recursively the {\it basic} tower of periodic points of $f$ associated with $P$ later on.}  

\vspace{.1in}

In summary, on each of the 5 intervals: $[\min P, d]$, $[d, u_1]$, $[\breve u_0', v]$ $(\subset [u_1, v])$, $[v, \hat u_0]$ $(\subset [v, \bar u_1'])$, $[\bar u_1', z_0]$, we have obtained 2 monotonic sequences of points (one sequence of periodic points which comprises the same (first) layer of the {\it basic} tower and the other is a sequence of $d$-points which are needed to continue the {\it basic} tower-building process into the next (second) layer of the tower): 
%{\large
\begin{multline*} 
$$
\min P < \tilde q_{m+2} < \tilde \mu_{m,1}' < \tilde q_{m+4} < \tilde \mu_{m,2}' < \cdots < d < \cdots < c_4 < u_2 < c_2 < u_1 \le \breve u_0' < \cdots \\ 
< \breve p_{m+4} < \breve \mu_{m,2} < \breve p_{m+2} < \breve \mu_{m,1} < v < \hat \mu_{m,1}' < \hat q_{m+2} < \hat \mu_{m,2}' < \hat q_{m+4} < \cdots < \hat u_0 \\ 
\le \bar u_1' < \bar c_4 < \bar u_2' < \bar c_6 < \bar u_3' < \bar c_8 < \cdots < z_0
$$
\end{multline*}
such that, $\breve u_0' = \max \big\{ d \le x \le v: f^2(x) = d \big\}$ and $\hat u_0 = \min \big\{ v \le x \le z_0: f^2(x) = d \big\}$ and, for each $n \ge 1$, $$
u_n = \min \big\{ d \le x \le v : f^{2n}(x) = d \big\} \,\,\, \text{and} \,\,\, \bar u_n' = \max \big\{ v \le x \le z_0 : f^{2n}(x) = d \big\},
$$ 
$$
d < f^{2k}(x) < z_0 \,\,\, \text{for all} \,\,\, 1 \le k \le n \,\,\, \text{and} \,\,\, \text{all} \,\,\, x \in (d, u_n) \cup (\bar u_n', z_0),
$$
\begin{multline*}
$$
\tilde \mu_{m,n-1}' = \max \big\{ \min P \le x \le d : f^{m+2n-2}(x) = d \big\}, \\ 
\breve \mu_{m,n} = \min \big\{ \breve u_0' \le x \le v : f^{m+2n}(x) = d \big\}, \\ 
\hat \mu_{m,n}' = \max \big\{ v \le x \le \hat u_0 : f^{m+2n}(x) = d \big\} \,\,\, \text{and}
$$
\end{multline*} 
the point $\tilde q_{m+2n-2} = \max \big\{ \min P \le x \le d: f^{m+2n-2}(x) = x \big\}$ is a period-$(m+2n)$ point of $f$, and on the interval $[d, \hat u_0]$ where the graph of $y = f^2(x)$ looks like a {\it skew} 'V' shape on the three points $d, v, \hat u_0$, we have, with $x$-coordinates moving from point $d$ via point $v$ to point $\hat u_0$,
\begin{itemize}
\item[{\rm (a)}]
(even periods) $c_{2n} = \min \big\{ d \le x \le v: f^{2n}(x) = x \big\}$ is a period-{\small $(2n)$} point of $f$;

\item[{\rm (b)}]
(odd periods) $\breve p_{m+2n} = \min \big\{ \breve u_0' \le x \le v: f^{m+2n}(x) = x \big\}$ is a period-{\small $(m+2n)$} point of $f$;

\item[{\rm (c)}] 
(odd periods) $\hat q_{m+2n} = \max \big\{ v \le x \le \hat u_0: f^{m+2n}(x) = x \big\}$ is a period-{\small $(m+2n)$} point of $f$,
\end{itemize}
and the point $\bar c_{2n+2} = \min \big\{ \bar u_n' \le x \le \bar u_{n+1}': f^{2n+2}(x) = x \big\}$ is a period-$(2n+2)$ point of $f$.

\vspace{.1in}

We call the collection of all these periodic points 
$$
\min P \le \tilde q_m, \,\,\, \tilde q_{m+2n}, \,\,\, c_{2n}, \,\,\, \breve p_{m+2n}, \,\,\, \hat q_{m+2n}, \,\,\, \bar c_{2n+2}, \,\,\, n \ge 1, \, \text{and} \,\,\, z_0
$$
the {\it first layer} of the {\it basic} tower of periodic points of $f$ associated with $P$ (see Figure 1).

\begin{figure}[htbp] 
\centering
\includegraphics[width=6in,height=2in]{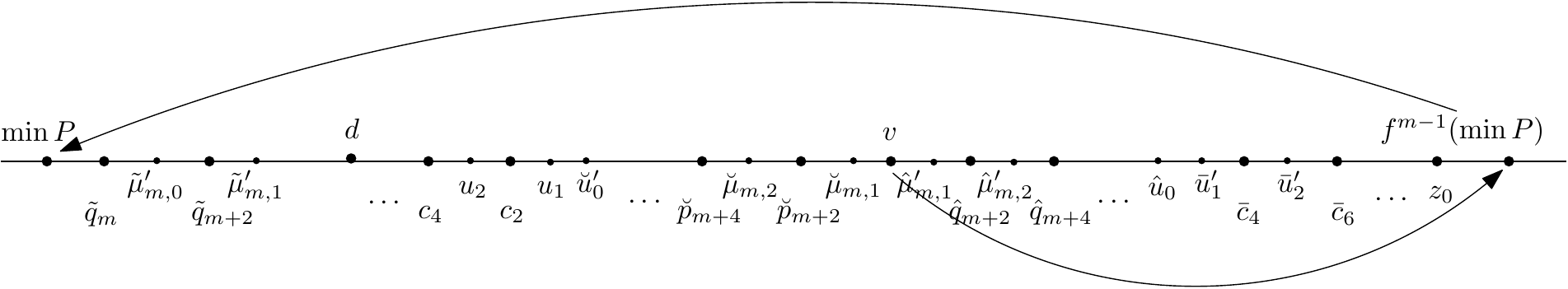} 
\vspace{.3in}

\noindent
{Figure 1: The first layer of the {\it basic} tower of periodic points of $f$ associated with $P$, \\ 
where $v$ is a point in $\big[\min P, f^{m-1}(\min P)\big]$ such that $f(v) = f^{m-1}(\min P)$; \vspace{.05in} \\ 
$z_0 = \min \big\{ v \le x \le f^{m-1}(\min P) : f^2(x) = x \big\}$; $d = \max \big\{ \min P \le x \le v : f^2(x) = z_0 \big\}$; \vspace{.05in} \\ 
and, for each $n \ge 1$, \vspace{.05in} \\ $\tilde q_{m+2n} = \max \big\{ \min P \le x \le d : f^{m+2n}(x) = x \big\}$; \vspace{.05in} \\ 
$\tilde \mu_{m,n}' = \max \big\{ \min P \le x \le d: f^{m+2n}(x) = d \big\}$; \vspace{.05in} \\ 
$u_n = \min \big\{ d \le x \le v : f^{2n}(x) = d \big\}$; \vspace{.05in} \\ 
$c_{2n} = \min \big\{ d \le x \le v : f^{2n}(x) = x \big\}$; \vspace{.05in} \\ 
$\breve u_0' = \max \big\{ d \le x \le v: f^2(x) = d \big\}$; \vspace{.05in} \\ 
$\breve \mu_{m,n} = \min \big\{ \breve u_0' \le x \le v: f^{m+2n}(x) = d \big\}$; \vspace{.05in} \\ 
$\breve p_{m+2n} = \min \big\{ \breve u_0' \le x \le v : f^{m+2n}(x) = x \big\}$;\vspace{.05in}  \\ 
$\hat u_0 = \min \big\{ v \le x \le z_0: f^2(x) = d \big\}$; \vspace{.05in} \\ 
$\hat \mu_{m,n}' = \max \big\{ v \le x \le \hat u_0: f^{m+2n}(x) = d \big\}$; \vspace{.05in} \\ 
$\hat q_{m+2n} = \max \big\{ v \le x \le \hat u_0 : f^{m+2n}(x) = x \big\}$; \vspace{.05in} \\ 
$\bar u_n' = \max \big\{ v \le x \le z_0 : f^{2n}(x) = d \big\}$; \vspace{.05in} \\ 
$\bar c_{2n+2} = \min \big\{ \bar u_n' \le x \le z_0 : f^{2n+2}(x) = x \big\}$.}%\vspace{.5in}}
\end{figure}

\pagebreak

\noindent
{\bf $\mathsection$2. The second layer of the {\it basic} tower of periodic points of $f$ associated with $P$.}

For each $n \ge 0$, let (note the relationship between the subscript of $\tilde \mu_{m,n}'$ and the superscript of $f^{m+2n}$ in the definition of $\tilde \mu_{m,n}'$) 
$$
\tilde \mu_{m,n}' = \max \big\{ \min P \le x \le d : f^{m+2n}(x) = d \big\}.
$$
Let $\tilde \nu_{m,n}$ be a point in $[\tilde \mu_{m,n}', d]$ such that $f^{m+2n}(\tilde \nu_{m,n}) = v$.  Then it turns out that $\tilde \mu_{m,n}' < \tilde \nu_{m,n} < \tilde \mu_{m,n+1}'$.  Furthermore, on the interval $[\tilde \mu_{m,n}', \tilde \mu_{m,n+1}']$, we have 
$$
f^{m+2n+2}(\tilde \mu_{m,n}') = z_0, \, f^{m+2n+2}(\tilde \nu_{m,n}) = \min P, \, f^{m+2n+2}(\tilde \mu_{m,n+1}') = d \,\,\, \text{and}
$$ 
$$
\text{the graph of} \,\,\, y = f^{m+2n+2}(x) \,\,\, \text{looks 'roughly' like a {\it skew} 'V' shape at the 3 points} \,\,\, \tilde \mu_{m,n}',
$$ 
$\tilde \nu_{m,n}, \tilde \mu_{m,n+1}'$ \big(meaning that the graph of $y = f^{m+2n+2}(x)$ passes through the 3 points $\, (\tilde \mu_{m,n}', z_0)$, $(\tilde \nu_{m,n}, \min P), (\tilde \mu_{m,n+1}', d)$ \, on the $x$-$y$ plane with $y$-coordinates $z_0, \min P, d$ respectively\big) just like that of $y = f^2(x)$ (at the points $d, v, \hat u_0$) on the interval $[d, \hat u_0]$.  

Now, for each $n \ge 0$ and all $i \ge 1$, let \big(the number 2 in the superscript $(n,2)$ indicates the second layer\big) 
$$
\quad\,\,\, \tilde p_{m+2n+2i}^{\,\, (n,2)} = \min \big\{ \tilde \mu_{m,n}' \le x \le \tilde \nu_{m,n} : f^{m+2n+2i}(x) = x \big\} \,\,\, \text{and}
$$ 
$$
\tilde \mu_{m,n,i} = \min \big\{ \tilde \mu_{m,n}' \le x \le \tilde \nu_{m,n} : f^{m+2n+2i}(x) = d \big\} \quad\quad\,
$$
\big(note the relationship between the subscript of $\tilde \mu_{m,n,i}$ and the superscript of $f^{m+2n+2i}$ in the definition of $\tilde \mu_{m,n,i}$.  Here the $n$ and $i$ in the subscript of $\tilde \mu_{m,n,i}$ indicate the $i^{th}$ point of the sequence $< \tilde \mu_{m,n,i} >$ in the $n^{th}$ interval $[\tilde \mu_{m,n}', \tilde \mu_{m,n+1}']$.  We need these points $\tilde \mu_{m,n,i}$'s to continue the {\it basic} tower-building process into the third layer\big).  By arguing as those in the proof of (c) in the main text, we obtain that 
$$
\tilde \mu_{m,n}' < \cdots < \tilde \mu_{m,n,3} < \tilde p_{m+2n+6}^{\,\, (n,2)} < \tilde \mu_{m,n,2} < \tilde p_{m+2n+4}^{\,\, (n,2)} < \tilde \mu_{m,n,1} < \tilde p_{m+2n+2}^{\,\, (n,2)} < \tilde \nu_{m,n}.
$$  
In the following, we shall show that, for each $n \ge 0$, with $x$-coordinates moving from point $\tilde \mu_{m,n}'$ via point $\tilde \nu_{m,n}$ to point $\tilde \mu_{m,n+1}'$ (see Figure 2):
\begin{itemize}
\item[{\rm (a)}]
(odd periods) the point $\tilde p_{m+2n+2i}^{\,\, (n,2)} = \min \big\{ \tilde \mu_{m,n}' \le x \le \tilde \nu_{m,n} : f^{m+2n+2i}(x) = x \big\} = \min \big\{ \tilde \mu_{m,n}' \le x \le \tilde \mu_{m,n+1}' : f^{m+2n+2i}(x) = x \big\}$ exists and is a period-$(m+2n+2i)$ point of $f$ for each $i \ge 1$;

\item[{\rm (b)}]
(even periods) the point $\tilde c_{2m+2n+2i}^{\,\, (n,2)} = \min \big\{ \tilde \mu_{m,n,1} \le x \le \tilde \nu_{m,n} : f^{2m+2n+2i}(x) = x \big\} = \min \big\{ \tilde \mu_{m,n,1} \le x \le \tilde \mu_{m,n+1}' : f^{2m+2n+2i}(x) = x \big\}$ exists for each $i \ge 1$ and is a period-$(2m+2n+2i)$ point of $f$ for each $i \ge n+3$;

\item[{\rm (c)}]
(even periods) the point $\tilde c_{2m+2n+2i}'^{\,\, (n,2)} = \max \big\{ \tilde \nu_{m,n} \le x \le \tilde \mu_{m,n+1}' : f^{2m+2n+2i}(x) = x \big\} = \max \big\{ \tilde \mu_{m,n}' \le x \le \tilde \mu_{m,n+1}' : f^{2m+2n+2i}(x) = x \big\}$ exists for each $i \ge 1$ and is a period-$(2m+2n+2i)$ point of $f$ for each $i \ge n+3$.
\end{itemize}

For each fixed $n \ge 0$, the collection of all these periodic points $\tilde p_{m+2n+2i+2}^{\,\,(n,2)}$, $\tilde c_{2m+2n+2i}^{\,\, (n,2)}$, $\tilde c_{2m+2n+2i}'^{\,\, (n,2)}$, $i \ge 1$ \big(note that the periodic point $\tilde p_{m+2n+2}^{\,\,(n,2)}$ is excluded so that the convex hulls of the sets $\big\{ \tilde p_{m+2n+2i}^{\,\,(n,2)}: i \ge 2 \big\}$ and $\big\{ \tilde c_{2m+2n+2i}^{\,\,(n,2)}: i \ge 1 \big\}$ are disjoint\big), is called a {\it compartment} of the second layer of the {\it basic} tower of periodic points of $f$ associated with $P$.  

\begin{figure}[htbp] 
\centering
\includegraphics[width=6in,height=2in]{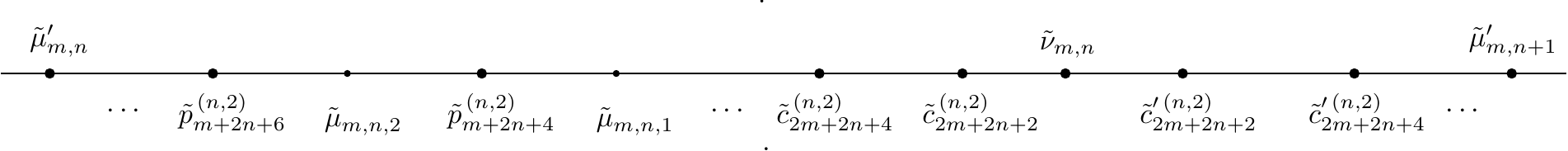}

{Figure 2: A compartment of the second layer of the {\it basic} tower of periodic points of $f$ associated with $P$ in the interval $[\tilde \mu_{m,n}', \tilde \mu_{m,n+1}']$ $(\subset [\min P, d])$, \\  
where $\tilde \nu_{m,n}$ is an {\it auxiliary} point in $[\tilde \mu_{m,n}', \tilde \mu_{m,n+1}']$ such that $f^{m+2n}(\tilde \nu_{m,n}) = v$ \\ 
which is needed only to determine the ordering of the following points: \\ 
and, for each $i \ge 1$, \\ 
$\tilde \mu_{m,n,i} = \min \big\{ \tilde \mu_{m,n}' \le x \le \tilde \mu_{m,n+1}' : f^{m+2n+2i}(x) = d \big\}$; \vspace{.05in} \\ 
$\tilde p_{m+2n+2i}^{\,(n,2)} = \min \big\{ \tilde \mu_{m,n}' \le x \le \tilde \mu_{m,n+1}' : f^{m+2n+2i}(x) = x \big\}$; \vspace{.05in} \\ 
$\tilde c_{2m+2n+2i}^{\, (n,2)} = \min \big\{ \tilde \mu_{m,n,1} \le x \le \tilde \mu_{m,n+1}' : f^{2m+2n+2i}(x) = x \big\}$; \vspace{.05in} \\ 
$\tilde c_{2m+2n+2i}'^{\, (n,2)} = \max \big\{ \tilde \mu_{m,n}' \le x \le \tilde \mu_{m,n+1}' : f^{2m+2n+2i}(x) = x \big\}$. \vspace{.05in}}
\end{figure}
\noindent   

\noindent
{\bf 2.1(a) On the existence of periodic points of $f$ of all odd periods $\ge m+2n+2$ in $[\tilde \mu_{m,n}', \tilde \nu_{m,n}]$ $(\subset [\tilde \mu_{m,n}', \tilde \mu_{m,n+1}']$ $\subset [\min P, d])$ for each $n \ge 0$.}

We can apply Lemma 5(1).  We first prove the following result (see also Lemma 11 below which is related) (Note that, for any periodic point $p$ of $f$, $\ell(p)$ denotes the least period of $p$ with respect to $f$):

\noindent
{\bf Lemma 7.}
{\it All periodic points of $f$ in $[\tilde \mu_{m,n}', d]$ of odd periods have least periods $\ge m+2n+2$.}

\noindent
{\it Proof.}
Suppose $f$ had a periodic point $\tilde p$ in $[\tilde \mu_{m,n}', d]$ of odd period with least period $\le m+2n$. 

\noindent
(1) If $\ell(\tilde p) = 1$, then $f^3(\tilde p) = \tilde p$ and $f^3(d) = f(z_0)$;  

\noindent
(2) If $\ell(\tilde p) \ge 3$, then $f^{\ell(\tilde p)}(\tilde p) = \tilde p$ and, since $\ell(\tilde p)$ is odd and $\ge 3$, $f^{\ell(\tilde p)}(d) = f(z_0)$.

Since $m+2n$ is odd and $\ge \max \{ 3, \ell(\tilde p) \}$, by Lemma 3, there is a periodic point $\tilde w_{m+2n}$ of $f$ in $[\tilde p, d]$ ($\subset (\tilde \mu_{m,n}', d]$) such that $f^{m+2n}(\tilde w_{m+2n}) = \tilde w_{m+2n}$.  Since 
$$
f^{m+2n}(\tilde w_{m+2n}) - d = \tilde w_{m+2n} - d < 0 \,\,\, \text{and} \,\,\, f^{m+2n}(d) - d = f(z_0) - d \ge z_0 - d > 0,
$$
there is a point $\tilde \mu_{m,n}^{\, *}$ in $(\tilde w_{m+2n}, d]$ such that $f^{m+2n}(\tilde \mu_{m,n}^{\, *}) = d$.  Since $\tilde \mu_{m,n}^{\, *} > \tilde w_{m+2n} > \tilde \mu_{m,n}'$. This contradicts the maximality of $\tilde \mu_{m,n}'$ in $[\min P, d]$.  Therefore, we have shown that all periodic points of $f$ in $[\tilde \mu_{m,n}', d]$ of odd periods have least periods $\ge m+2n+2$.
$\hfill\square$

Now, since $f^{m+2n}(\tilde \nu_{m,n}) = v$, we apply Lemma 5(1) with 
{\large 
$$
f^{m+2n+2}(\tilde \mu_{m,n}') = z_0 \,\,\, \text{and} \,\,\, f^{m+2n+2}(\tilde \nu_{m,n}) = \min P \,\,\, \text{and Lemma 7}
$$}
to obtain that, for each $n \ge 0$ and $i \ge 0$, the points 
$$
\tilde p_{m+2n+2+2i}^{\,\, (n,2)} = \min \big\{ \tilde \mu_{m,n}' \le x \le \tilde \mu_{m,n+1}' : f^{m+2n+2+2i}(x) = x \big\} \,\,\, \text{and}
$$ 
$$
\tilde \mu_{m,n,i+1} = \min \big\{ \tilde \mu_{m,n}' \le x \le \tilde \mu_{m,n+1}' : f^{m+2n+2+2i}(x) = d \big\}\qquad\quad\,\,\,
$$
exist, $\tilde \mu_{m,n}' < \cdots < \tilde \mu_{m,n,3} < \tilde p_{m+2n+6}^{\,\, (n,2)} < \tilde \mu_{m,n,2} < \tilde p_{m+2n+4}^{\,\, (n,2)} < \tilde \mu_{m,n,1} < \tilde p_{m+2n+2}^{\,\, (n,2)} < \tilde \mu_{m,n+1}'$ and, for each $n \ge 0$ and $i \ge 0$, the point $\tilde p_{m+2n+2+2i}^{\,\, (n,2)}$ is a period-$(m+2n+2+2i)$ point of $f$, or, equivalently, in $[\tilde \mu_{m,n}', \tilde \nu_{m,n}]$ $(\subset [\tilde \mu_{m,n}', \tilde \mu_{m,n+1}']$ $\subset [\min P, d])$, 
$$
\text{for each} \,\,\, n \ge 0 \,\,\, \text{and} \,\,\, i \ge 1, \, \text{the point} \,\,\, \tilde p_{m+2n+2i}^{\,\, (n,2)} \,\,\, \text{is a period-}(m+2n+2i) \,\,\, \text{point of} \,\,\, f.
$$

\noindent
{\bf 2.1(b) On the existence of periodic points of $f$ of all even periods $\ge 2m+4n+6$ in $[\tilde \mu_{m,n,1}, \tilde \nu_{m,n}]$ $(\subset [\tilde \mu_{m,n}', \tilde \mu_{m,n+1}']$ $\subset [\min P, d])$ for each $n \ge 0$.}

In this case, we shall apply Lemma 6(1).  Since, for each $i \ge 1$, $\tilde \mu_{m,n,i} = \min \big\{ \tilde \mu_{m,n}' \le x \le \tilde \nu_{m,n} : f^{m+2n+2i}(x) = d \big\}$, we have $f^{m+2n+2}(\tilde \mu_{m,n,1}) = d$ and so, 
$$
f^{m+2n+3}(\tilde \mu_{m,n,1}) = f(d).
$$
Furthermore, since $f^{m+2n}(\tilde \nu_{m,n}) = v$, we have $f^{m+2n+2}(\tilde \nu_{m,n}) = f^2(v) = \min P$.  Consequently, we have 
$$
f^{2m+2n+2}(\tilde \nu_{m,n}) = f^m(\min P) = \min P.
$$
Therefore, we can apply Lemma 6(1) with 
{\large 
$$
f^{m+2n+3}(a) = f(d) \in \big\{ z_0, f(z_0), f(d) \big\} \,\,\, \text{and} \,\,\, f^{2(m+n+1)}(b) = \min P
$$}
to obtain that, for each $i \ge 0$, the point
$$
\tilde c_{2(m+n+1)+2i}^{\,\, (n,2)} = \min \big\{ \tilde \mu_{m,n,1} \le x \le \tilde \nu_{m,n}: f^{2(m+n+1)+2i}(x) = x \big\}
$$
exists and 
$$
\tilde \mu_{m,n,1} < \cdots < \tilde c_{2m+2n+6}^{\, (n,2)} < \tilde c_{2m+2n+4}^{\, (n,2)} < \tilde c_{2m+2n+2}^{\, (n,2)} < \tilde \nu_{m,n}.
$$
Furthermore, for each $i \ge 0$ such that $(m+n+1)+i \ge \max \big\{ m+2n+3, m+2 \big\} = m+2n+3$, i.e., for each $i \ge n+2$, the point $\tilde c_{2(m+n+1)+2i}^{\,\, (n,2)}$ is a period-$\big(2(m+n+1)+2i\big)$ point of $f$, or, equivalently, in $[\tilde \mu_{m,n,1}, \tilde \nu_{m,n}]$ $(\subset [\tilde \mu_{m,n}', \tilde \mu_{m,n+1}']$ $\subset [\min P, d])$,
$$
\text{for each} \,\,\, n \ge 0 \,\,\, \text{and} \,\,\, i \ge n+3, \, \text{the point} \,\,\, \tilde c_{2m+2n+2i}^{\,\, (n,2)} \,\,\, \text{is a period-}(2m+2n+2i) \,\,\, \text{point of} \,\,\, f.
$$  

\noindent
{\bf 2.1(c) On the existence of periodic points of $f$ of all even periods $\ge 2m+4n+6$ in $[\tilde \nu_{m,n}, \tilde \mu_{m,n+1}']$ $(\subset [\tilde \mu_{m,n}', \tilde \mu_{m,n+1}']$ $\subset [\min P, d])$ for each $n \ge 0$.}

In this case, we shall apply Lemma 6(2).  By definition of $\tilde \mu_{m,n+1}'$, we have 
$$
f^{m+2n+3}(\tilde \mu_{m,n+1}') = f\big(f^{m+2n+2}(\tilde \mu_{m,n+1}')\big) = f(d).
$$
On the other hand, since $\tilde \nu_{m,n}$ is a point in $[\tilde \mu_{m,n}', \tilde \mu_{m,n+1}']$ such that $f^{m+2n}(\tilde \nu_{m,n}) = v$, we have 
$$
f^{m+2n+2}(\tilde \nu_{m,n}) = f^2(v) = \min P.
$$
Therefore, we can apply Lemma 6(2) with 
{\large 
$$
f^{m+2n+3}(\tilde \mu_{m,n+1}') = f(d) \in \big\{ z_0, f(z_0), f(d) \big\} \,\,\, \text{and} \,\,\, f^{2(m+n+1)}(\tilde \nu_{m,n}) = \min P
$$}
to obtain that, for each $i \ge 0$, the point 
$$
\tilde c_{2(m+n+1)+2i}'^{\,\, (n,2)} = \max \big\{ \tilde \nu_{m,n} \le x \le \tilde \mu_{m,n+1}': f^{2(m+n+1)+2i}(x) = x \big\}
$$ 
exists and
$$
\tilde \nu_{m,n} < \tilde c_{2m+2n+2}'^{\,\, (n,2)} < \tilde c_{2m+2n+4}'^{\,\, (n,2)} < \tilde c_{2m+2n+6}'^{\,\, (n,2)} < \cdots < \tilde \mu_{m,n+1}'.
$$
Furthermore, for each $i \ge 0$ such that $(m+n+1)+i \ge \max \big\{ m+2n+3, m+2 \big\} = m+2n+3$, i.e., for each $i \ge n+2$, the point $\tilde c_{2(m+n+)+2i}'^{\,\, (n,2)}$ is a period-$\big(2(m+n+1)+2i\big)$ point of $f$, or, equivalently, in $[\tilde \nu_{m,n}, \tilde \mu_{m,n+1}']$ $(\subset [\tilde \mu_{m,n}', \tilde \mu_{m,n+1}']$ $\subset [\min P, d])$,
$$
\text{for each} \,\,\, n \ge 0 \,\,\, \text{and} \,\,\, i \ge n+3, \, \text{the point} \,\,\, \tilde c_{2m+2n+2i}'^{\,\, (n,2)} \,\,\, \text{is a period-}(2m+2n+2i) \,\,\, \text{point of} \,\,\, f.
$$
$$\aleph \qquad\qquad\qquad \aleph \qquad\qquad\qquad \aleph \qquad\qquad\qquad \aleph \qquad\qquad\qquad \aleph$$
\indent For each $n \ge 1$, let $u_n = \min \big\{ d \le x \le v : f^{2n}(x) = d \big\}$ be defined as before.  Then we have $d < \cdots < u_3 < u_2 < u_1 < v$, $f^2(d) = z_0$ and $f^{2n}(u_n) = d$.  In particular, we have $f^{2n}\big([d, u_n]\big) = [d, z_0] \supset \{ v \}$.  Let $\nu_n$ be a point in $[d, u_n]$ such that 
$$
f^{2n}(\nu_n) = v.
$$
Then, since $f^2(v) = \min P$, we obtain that $f^{2n+2}\big([d, \nu_n]\big) \supset [\min P, z_0] \supset \{ d \}$.  So, the point $\nu_n$ happens to satisfy that $u_{n+1} < \nu_n < u_n$.  On the interval $[u_{n+1}, u_n]$, we have 
$$
f^{2n+2}(u_{n+1}) = d, \,\, f^{2n+2}(\nu_n) = \min P, \,\, f^{2n+2}(u_n) = z_0 \,\,\, \text{and}
$$ 
$$
\text{the graph of} \,\,\, y = f^{2n+2}(x) \,\,\, \text{looks 'roughly' like a {\it skew} 'V' shape at the points} \,\,\, u_{n+1}, \nu_n, u_n
$$ 
\big(meaning that the graph of $y = f^{2n+2}(x)$ passes through the 3 points $\, (u_{n+1}, d)$, $(\nu_n, \min P)$, $(u_n, z_0)$ \, on the $x$-$y$ plane with $y$-coordinates $d, \min P, z_0$ respectively\big) just like a mirror-symmetric copy of that of $y = f^2(x)$ at the 3 points $d, v, \hat u_0$ on the interval $[d, \hat u_0]$.  

In the following, we show that, for each $n \ge 1$, there exist 4 sequences of points in $[u_{n+1}, u_n]$: \big(the 2 in the superscript $(n, 2)$ below indicates second layer\big)
{\large 
\begin{multline*}
$$
u_{n+1} < \cdots < p_{m+2n+4}^{\,\,(n,2)} < p_{m+2n+2}^{\,\,(n,2)} < \nu_n < q_{m+2n+2}^{\,\,(n,2)} < q_{m+2n+4}^{\,\,(n,2)} < \cdots \\ < u_{n,1}' < c_{2n+2}'^{\,\,(n,2)} < u_{n,2}' < c_{2n+4}'^{\,\,(n,2)} < \cdots < u_n
$$
\end{multline*}}
such that, for each $i \ge 1$, $u_{n,i}' = \max \big\{ u_{n+1} \le x \le u_n : f^{2n+2i}(x) = d \big\}$ and, from point $u_n$ via point $\nu_n$ to point $u_{n+1}$ (see Figure 3):
\begin{itemize}
\item[{\rm (a)}]
(even periods) the point $c_{2n+2i}'^{\,\, (n,2)} = \max \big\{ \nu_n \le x \le u_n: f^{2n+2i}(x) = x \big\} = \max \big\{ u_{n+1} \le x \le u_n: f^{2n+2i}(x) = x \big\}$ is a period-$(2n+2i)$ point of $f$ for each $i \ge 1$ and $i \ne n$;

\item[{\rm (b)}]
(odd periods) the point $q_{m+2n+2i}^{\,\, (n,2)} = \max \big\{ \nu_n \le x \le u_{n,1}': f^{m+2n+2i}(x) = x \big\} = \max \big\{ u_{n+1} \le x \le u_{n,1}': f^{m+2n+2i}(x) = x \big\}$ is a period-$(m+2n+2i)$ point of $f$ for each $i \ge 1$;

\item[{\rm (c)}]
(odd periods) the point $p_{m+2n+2i}^{\,\, (n,2)} = \min \big\{ u_{n+1} \le x \le \nu_n: f^{m+2n+2i}(x) = x \big\} = \min \big\{ u_{n+1} \le x \le u_n : f^{m+2n+2i}(x) = x \big\}$ is a period-$(m+2n+2i)$ point of $f$ for each $i \ge 1$.
\end{itemize}

\begin{figure}[htbp] 
\centering
\includegraphics[width=6in,height=2in]{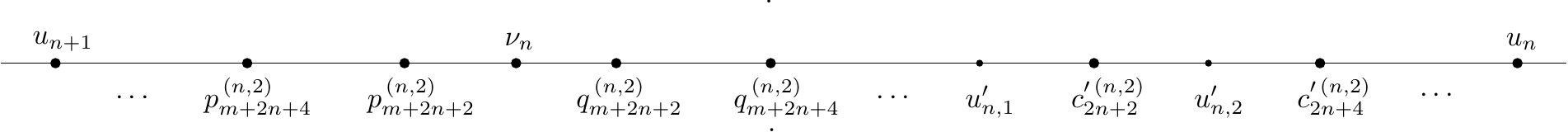} 

{Figure 3: A compartment of the second layer of the {\it basic} tower of periodic points of $f$ associated with $P$ in the interval $[u_{n+1}, u_n]$ $(\subset [d, u_1])$, \\
where $\nu_n$ is an {\it auxiliary} point in $[\min P, u_n]$ such that $f^{2n}(\nu_n) = v$ \\ 
which is needed only to determine the ordering of the following points: \\ 
and, for each $i \ge 1$, \\ 
$u_{n,i}' = \max \big\{ u_{n+1} \le x \le u_n : f^{2n+2i}(x) = d \big\}$; \vspace{.05in} \\ 
$c_{2n+2i}'^{\, (n,2)} = \max \big\{ u_{n+1} \le x \le u_n : f^{2n+2i}(x) = x \big\}$; \vspace{.05in} \\ 
$q_{m+2n+2i}^{\, (n,2)} = \max \big\{ u_{n+1} \le x \le u_{n,1}' : f^{m+2n+2i}(x) = x \big\}$; \vspace{.05in} \\ 
$p_{m+2n+2i}^{\, (n,2)} = \min \big\{ u_{n+1} \le x \le u_n : f^{m+2n+2i}(x) = x \big\}$. \vspace{.05in}} 
\end{figure}
\indent For each fixed $n \ge 1$, the collection of all these periodic points $c_{2m+2n+2i}'^{\,\, (n,2)}$, $q_{m+2n+2i}^{\,\, (n,2)}$, $p_{m+2n+2i}^{\,\, (n,2)}$, $i \ge 1$, is called a {\it compartment} of the second layer of the {\it basic} tower of periodic points of $f$ associated with $P$.    

\noindent
{\bf 2.2(a) On the existence of periodic points of $f$ of all even periods $\ge 2n+2$ (except possibly period $4n$) in $[\nu_n, u_n]$ $(\subset [u_{n+1}, u_n]$ $\subset [d, v])$ for each $n \ge 1$.}  

We can apply Lemma 6(2) to this case.  However, since, in this special case, we can {\it determine their exact periods}, we give a proof in detail.  For higher layers (layer 3 and up) of periodic points, we can determine their exact periods only when the periods are large enough.  

Recall that $f^{2n}(\nu_n) = v$.  On $[\nu_n, u_n]$, since $f^{2n+2}\big([\nu_n, u_n]\big) \supset [\min P, z_0] \supset \{ d \}$, the point 
$$
u_{n,1}' = \max \big\{ \nu_n \le x \le u_n : f^{2n+2}(x) = d \big\} = \max \big\{ u_{n+1} \le x \le u_n : f^{2n+2}(x) = d \big\}
$$
exists and $d < f^{2n+2}(x) < z_0$ for all $u_{n,1}' < x < u_n$.  This, combined with the following fact 
$$
\qquad\quad\qquad d < f^{2i}(x) < z_0 \,\,\, \text{for all} \,\,\, 1 \le i \le n \,\,\, \text{and all} \,\,\, x \in (d, u_n] \,\, (\cup \,\, [\bar u_n', z_0)) \qquad\qquad\quad\,\,\,\,\,\, (\dagger)
$$
implies that 
$$
\qquad\quad\qquad d < f^{2i}(x) < z_0 \,\,\, \text{for all} \,\,\, 1 \le i \le n+1 \,\,\, \text{and all} \,\,\, u_{n,1}' < x < u_n.\qquad\qquad\qquad\quad\,\,\,\,\, (\ddag)
$$
Let $c_{2n+2}'^{\,\,(n,2)} = \max \big\{ u_{n,1}' \le x \le u_n : f^{2n+2}(x) = x \big\} = \max \big\{ u_{n+1} \le x \le u_n : f^{2n+2}(x) = x \big\}$.  Since $f^{2n+2}\big([c_{2n+2}'^{\,\,(n,2)}, u_n]\big) \supset [c_{2n+2}'^{\,\,(n,2)}, z_0] \supset \{ v \}$, let $\nu_{n+1}$ be a point in $[c_{2n+2}'^{\,\,(n,2)}, u_n]$ such that $f^{2n+2}(\nu_{n+1}) = v$.  Then $f^{2n+4}\big([\nu_{n+1}, u_n]\big) \supset [\min P, z_0] \supset \{ d \}$.  So, the point 
$$
u_{n,2}' = \max \big\{ \nu_{n+1} \le x \le u_n : f^{2n+4}(x) = d \big\} = \max \big\{ u_{n+1} \le x \le u_n : f^{2n+4}(x) = d \big\}
$$
exists and $d < f^{2n+4}(x) < z_0$ for all $u_{n,2}' < x < u_n$.  This, combined with ($\ddag$) above, implies that 
$$
d < f^{2i}(x) < z_0 \,\,\, \text{for all} \,\,\, 1 \le i \le n+2 \,\,\, \text{and all} \,\,\, u_{n,2}' < x < u_n.
$$
Let $c_{2n+4}'^{\,\,(n,2)} = \max \big\{ u_{n,2}' \le x \le u_n : f^{2n+4}(x) = x \big\} = \max \big\{ u_{n+1} \le x \le u_n : f^{2n+4}(x) = x \big\}$.  Proceeding in this manner indefinitely, we obtain 2 sequences of points 
$$
u_{n+1} \, < \, \nu_n \, < \, u_{n,1}' \, < \, c_{2n+2}'^{\,\,(n,2)} \, < \, u_{n,2}' \, < \, c_{2n+4}'^{\,\,(n,2)} \, < \, u_{n,3}' \, < \, c_{2n+6}'^{\,\,(n,2)} \, < \, \cdots \, < \, u_n
$$
such that, for each $i \ge 1$, 
$$
u_{n,i}' = \max \big\{ \nu_n \le x \le u_n : f^{2n+2i}(x) = d \big\} \\ = \max \big\{ u_{n+1} \le x \le u_n : f^{2n+2i}(x) = d \big\},
$$ 
{\large 
$$
\qquad\,\, d < f^{2j}(x) < z_0 \,\,\, \text{for all} \,\,\, 1 \le j \le n+i \,\,\, \text{and all} \,\,\, u_{n,i}' < x < u_n, \qquad\quad\,\, (\dagger\dagger)
$$} 
$$
\text{and} \,\,\, c_{2n+2i}'^{\,\,(n,2)} = \max \big\{ u_{n,i}' \le x \le u_n : f^{2n+2i}(x) = x \big\} = \max \big\{ u_{n+1} \le x \le u_n : f^{2n+2i}(x) = x \big\}.
$$

Note that it follows from $(\ddag$) that each $c_{2n+2i}'^{\,\,(n,2)}$ is a periodic point of $f$ with {\it even} period.  To find the least period of $c_{2n+2i}'^{\,\,(n,2)}$ with respect to $f$, we shall use the following result which can be proved by using arguments similar to those at the end of the proof of (c) in the main text (cf. Lemma 4) \big(see also the relevant Lemma 12 below.  Unfortunately, we do not have similar result on the counterpart interval $[\bar u_n', \bar u_{n+1}']$\big):

\noindent
{\bf Lemma 8.}
{\it For each $n \ge 1$ and $i \ge 1$, on the interval $[u_{n,i}', u_n]$ $(\subset [u_{n+1}, u_n])$, $f$ has no periodic points of odd periods $\le 2n+2i+1$, nor has periodic points of even periods $\le 2n+2i$ except period-$(2n+2i)$ points and possibly period-$(2n)$ points.}

\noindent
{\it Proof.}
Suppose $f$ had a period-$(2n+2j)$ point, say $w_{2n+2j}$, in $[u_{n,i}', u_n]$ for some $1 \le j \le i-1$.  Then since $f^{2n+2j+2}\big([w_{2n+2j}, u_n]\big) \supset f^2\big([w_{2n+2j}, z_0]\big) \supset f^2\big([v, z_0]\big) \supset [\min P, z_0] \supset \{ d \}$, there exists a point $u_{n+j+1}^*$ in $(w_{2n+2j}, u_n] \, (\subset [u_{n,i}', u_n])$ such that $f^{2n+2j+2}(u_{n+j+1}^*) = d$.  So, $u_{n,i}' < u_{n+i+1}^*$.  Since $n+j+1 \le n+i$, this contradicts the fact that $(u_{n+j+1}^* \le) \, u_{n,j+1}' \le u_{n,i}'$.  This, combined with Lemma 4(1), implies that $f$ has no periodic points of even periods $\le 2n+2i-2$ except possibly period-$(2n)$ points.  Furthermore, it follows from the above $(\dagger\dagger)$ that $f$ has no periodic points of odd periods $\le 2n+2i+1$.  

On the other hand, since $f^{2n+2i}(u_{n,i}') - u_{n,i}'$ $= d - u_{n,i} < 0$ and $f^{2n+2i}(u_n) - u_n$ $= f^{2i}(d) - u_n$ $= z_0 - u_n > 0$, the point $c_{2n+2i}'^{\,(*)} = \max \big\{ u_{n,i}' \le x \le u_n: f^{2n+2i}(x) = x \big\}$ exists and, since we have just shown in the previous paragraph that $f$ has no periodic points of least (odd or even) periods $\le 2n+2i-1$, the point $c_{2n+2i}'^{\,(*)}$ is a period-$(2n+2i)$ point of $f$ in $[u_{n,i}', u_n]$.
$\hfill\square$

We now use Lemma 8 to find the least period of $c_{2n+2i}'^{\,\,(n,2)}$ with respect to $f$ for each $i \ge 1$:

Let $\ell\big(c_{2n+2i}'^{\,\,(n,2)}\big)$ denote the least period of $c_{2n+2i}'^{\,\,(n,2)}$ with respect to $f$.  Suppose $\ell\big(c_{2n+2j}'^{\,\,(n,2)}\big) < 2n+2j$ for some $j \ge 1$.  Then by Lemma 8, $\ell\big(c_{2n+2j}'^{\,\,(n,2)}\big) = 2n$ and by Lemma 1, $2n+2j$ is a multiple of $2n$.  Suppose $(2n+2j)/(2n) = r > 2$.  Then $(r-1)(2n) > 2n$.  So, $(r-1)(2n) \ge 2n+2$.  Since $f^{2n+2}(u_n) = z_0$, $f^{(r-1)(2n)}\big(c_{2n+2j}'^{\,\,(n,2)}\big) = c_{2n+2j}'^{\,\,(n,2)}$ and $2n+2j > (r-1)(2n)$, it follows from Lemma 3 that there exists a periodic point, say $w_{2n+2j}^{*}$, of $f$ such that $u_{n,j}' <  c_{2n+2j}'^{\,\,(n,2)} < w_{2n+2j}^{*} < u_n$ and $f^{2n+2j}(w_{2n+2j}^{*}) = w_{2n+2j}^{*}$ which contradicts the maximality of $c_{2n+2j}'^{\,\,(n,2)}$ in $[u_{n,j}', u_n]$.  This shows that, for each $i \ge 1$, 
$$
\quad\text{if} \,\,\, (2n+2i)/(2n) \ne 2, \, \text{then} \,\,\, c_{2n+2i}'^{\,\,(n,2)} \,\,\, \text{is a peirod-}(2n+2i) \,\,\, \text{point of} \,\,\, f \,\,\, \text{and,} \qquad\quad \,\,\,\,\,\,
$$ 
$$
\,\,\,\text{if} \,\,\, 2n+2i = 4n, \, \text{then} \,\,\, c_{4n}'^{\,\,(n,2)} \,\,\, \text{is either a period-}(2n) \,\,\, \text{or a period-}(4n) \,\,\, \text{point of} \,\,\, f.
$$  

Note that the point $c_{4n}'^{\,\,(n,2)}$ may be a period-$(2n)$ point of $f$.  For example, when $n = 1$, let $g : [0, 1] \longrightarrow [0, 1]$ be the continuous map defined by putting (i) $g(x) = x + 1/2$ and (ii) $g(x) = 2 - 2x$ and let $P = \{ 0, 1/2, 1 \}$ be the unique period-3 orbit of $g$, then $g$ has exactly 4 period-4 points (which form a period-4 orbit of $g$), i.e., $2/9 \rightarrow 13/18 \rightarrow 5/9 \rightarrow 8/9 \rightarrow 2/9$.  In this case, $\min P = 0$ and $b = g^2(\min P) = 1$ and 
$$
d = \frac 16, \, c_4 = \frac 29, \, u_2 = \frac {11}{48}, \, u_{1,1,2} = \frac 7{24}, \, c_2 = \frac 13, \, u_1 = \frac 5{12}, \, v = \frac 12, \, \hat c_4 = \frac 59, \, z_0 = z = \frac 23.
$$
We see that, in the interval $[u_2, u_1] = [11/48, 5/12] \, (\supset [u_{1,1}, u_1])$, $g$ has only period-2 point, i.e., $c_4' = c_2 = 1/3$, but has no period-4 point.  

\noindent
{\bf 2.2(b) On the existence of periodic points of all odd periods $\ge m+2n+2$ in $[\nu_n, u_{n,1}']$ $(\subset [u_{n+1}, u_n])$ $\subset [d, v])$ for each $n \ge 1$.}

By definition of $u_{n,1}'$, we have $f^{2n+2}(u_{n,1}') = d$ and so, 
$$
f^{2n+3}(u_{n,1}') = f(d).
$$
Since $\nu_n$ is a point in $[u_{n+1}, u_n]$ such that $f^{2n}(\nu_n) = v$, we have 
$$
f^{m+2n+2}(\nu_n) = f^{m+2}(v) = f^m(\min P) = \min P.
$$
On the other hand, it follows from Lemma 4(1) that $$\text{all periodic points of} \,\, f \,\, \text{in} \, [d, u_n] \, \text{with {\it odd} periods have periods} \, \ge 2n+3.\,$$Therefore, by applying Lemma 5(2) with 
{\large 
$$
f^{2n+3}(u_{n,1}') = f(d) \in \big\{ z_0, f(z_0), f(d) \big\} \,\,\, \text{and} \,\,\, f^{m+2(n+1)}(\nu_n) = \min P,
$$}
we obtain that, for each $i \ge 0$, the point 
\begin{multline*}
$$
\qquad\qquad\qquad q_{m+2(n+1)+2i}^{\, (n,2)} = \max \big\{ \nu_n \le x \le u_{n,1}' : f^{m+2(n+1)+2i}(x) = x \big\} \\ 
= \max \big\{ u_{n+1} \le x \le u_{n,1}' : f^{m+2(n+1)+2i}(x) = x \big\}\qquad\qquad\quad\quad\,\,\,\,\,
$$
\end{multline*}
exists and 
$\,\, u_{n+1} < \,\, \nu_n \,\, < \,\, q_{m+2n+2}^{\, (n,2)} \,\, < \,\, q_{m+2n+4}^{\, (n,2)} \,\, < \,\, q_{m+2n+6}^{\, (n,2)} \,\, < \,\, \cdots \,\, < \,\, u_{n,1}' \,\, < \,\, u_n \,\, < \,\, v$.  Furthermore, for each $i \ge 0$, the point $q_{m+2n+2+2i}^{\, (n,2)}$ is a period-$(m+2n+2+2i)$ point of $f$, or,
equivalently,  in $[\nu_n, u_{n,1}']$ $(\subset [u_{n+1}, u_n])$ $\subset [d, v])$,
$$
\text{for each} \,\,\, n \ge 1 \,\,\, \text{and} \,\,\, i \ge 1, \text{the point} \,\,\, q_{m+2n+2i}^{\, (n,2)} \,\,\, \text{is a period-}(m+2n+2i) \,\,\, \text{point of} \,\,\, f.
$$  

\noindent
{\bf 2.2(c) On the existence of periodic points of all odd periods $\ge m+2n+2$ in $[u_{n+1}, \nu_n]$ $(\subset [u_{n+1}, u_n])$ $\subset [d, v])$ for each $n \ge 1$.}

By definition of $u_{n+1}$, we have $f^{2n+2}(u_{n+1}) = d$ and so, $f^{2n+3}(u_{n+1}) = f(d)$.  Since $\nu_n$ is a point in $[u_{n+1}, u_n]$ such that $f^{2n}(\nu_n) = v$, we have $f^{m+2n+2}(\nu_n) = f^{m+2}(v) = f^m(\min P) = \min P$.  On the other hand, it follows from Lemma 4(1) that 
$$
\text{all periodic points of} \,\, f \,\, \text{in} \, [d, u_n] \, \text{with {\it odd} periods have periods} \, \ge 2n+3.\,
$$
Therefore, by Lemma 5(1) with 
{\large 
$$
f^{2n+3}(u_{n+1}) = f(d) \in \big\{ z_0, f(z_0), f(d) \big\} \,\,\, \text{and} \,\,\, f^{m+2(n+1)}(\nu_n) = \min P,
$$}
we obtain that, for each $i \ge 0$, the point 
\begin{multline*}
$$
\qquad\qquad\qquad p_{m+2(n+1)+2i}^{\, (n,2)} = \min \big\{ u_{n+1} \le x \le \nu_n : f^{m+2(n+1)+2i}(x) = x \big\} \\ 
= \min \big\{ u_{n+1} \le x \le u_n : f^{m+2(n+1)+2i}(x) = x \big\},\qquad\qquad\qquad\quad\,\,
$$
\end{multline*}
exists and $\,\, d \,\, < \,\, u_{n+1} \,\, < \,\, \cdots \,\, < \,\, p_{m+2n+6}^{\, (n,2)} \,\, < \,\, p_{m+2n+4}^{\, (n,2)} \,\, < \,\, p_{m+2n+2}^{\, (n,2)} \,\, < \,\, \nu_n \,\, < \,\, u_n \,\, < \,\, v$.  Furthermore, for each $i \ge 0$, the point $p_{m+2n+2+2i}^{\, (n,2)}$ is a period-$(m+2n+2+2i)$ point of $f$, or, equivalently, in $[u_{n+1}, \nu_n]$ $(\subset [u_{n+1}, u_n])$ $\subset [d, v])$,
$$
\text{for each} \,\,\, n\ge 1 \,\,\, \text{and} \,\,\, i \ge 1, \text{the point} \,\,\, p_{m+2n+2i}^{\, (n,2)} \,\,\, \text{is a period-}(m+2n+2i) \,\,\, \text{point of} \,\,\, f.
$$
$$\aleph \qquad\qquad\qquad \aleph \qquad\qquad\qquad \aleph \qquad\qquad\qquad \aleph \qquad\qquad\qquad \aleph$$
\indent Recall that $\breve u_0' = \max \big\{ u_1 \le x \le v: f^2(x) = d \big\}$.  Since $f^{m+2}(\breve u_0') = f^m(d) = f(z_0) \ge z_0$ and $f^{m+2}(v) = f^m(\min P) = \min P$,  we have $f^{m+2}\big([\breve u_0', v]\big) \supset [\min P, z_0] \supset \{ d \}$.  So, the point 
$$
\breve \mu_{m,1} = \min \big\{ \breve u_0' \le x \le v: f^{m+2}(x) = d \big\}
$$
exists.  Since $f^{m+4}\big([\breve u_0', \breve \mu_{m,1}]\big) = f^2\big(f^{m+2}([\breve u_0', \breve \mu_{m,1}])\big) \supset f^2\big([d, z_0]\big) \supset [\min P, z_0] \supset \{ d \}$, the point 
$$
\breve \mu_{m,2} = \min \big\{ \breve u_0' \le x \le \breve \mu_{m,1}: f^{m+4}(x) = d \big\} = \min \big\{ \breve u_0' \le x \le v: f^{m+4}(x) = d \big\}
$$
exists and is $< \breve \mu_{m,1}$.  Inductively, for each $n \ge 1$, let (note the relationship between the subscript of $\breve \mu_{m,n}$ and the superscript of $f^{m+2n}$ in the definition of $\breve \mu_{m,n}$) 
$$
\breve \mu_{m,n} = \min \big\{ \breve u_0' \le x \le v : f^{m+2n}(x) = d \big\}.
$$
Then we have \,\,\, $u_1 \, \le \, \breve u_0' \, < \, \cdots \, < \, \breve \mu_{m,3} \, < \, \breve \mu_{m,2} \, < \, \breve \mu_{m,1} \, < \, v$.  For each $n \ge 1$, let $\breve \nu_{m,n}$ be a point in $[\breve u_0', \breve \mu_{m,n}]$ such that $f^{m+2n}(\breve \nu_{m,n}) = v$.  Then it turns out that $\breve \mu_{m,n+1} < \breve \nu_{m,n} < \breve \mu_{m,n}$.  Furthermore, on the interval $[\breve \mu_{m,n+1}, \breve \mu_{m,n}]$, we have 
$$
f^{m+2n+2}(\breve \mu_{m,n+1}) = d, \, f^{m+2n+2}(\breve \nu_{m,n}) = \min P, \, f^{m+2n+2}(\breve \mu_{m,n}) = z_0 \,\,\, \text{and}
$$ 
$$
\text{the graph of} \,\,\, y = f^{m+2n+2}(x) \,\,\, \text{looks 'roughly' like a {\it skew} 'V' shape at the points} \,\,\, \breve \mu_{m,n+1},
$$ 
$\breve \nu_{m,n}, \breve \mu_{m,n}$ \big(meaning that the graph of $y = f^{m+2n+2}(x)$ passes through the 3 points $\, (\breve \mu_{m,n+1}, d)$, $(\breve \nu_{m,n}, \min P), (\breve \mu_{m,n}, z_0)$ on the $x$-$y$ plane with $y$-coordinates $d, \min P, z_0$ respectively\big) just like a mirror-symmetric copy of that of $y = f^2(x)$ (at the 3 points $d, v, \hat u_0$) on the interval $[d, \hat u_0]$, where $\hat u_0 = \min \{ v \le x \le z_0: f^2(x) = d \}$.  Now, for each $n \ge 1$ and all $i \ge 1$, let \big(the number 2 in the superscript $(n,2)$ indicates the second layer\big) 
$$
\quad\,\,\,\, \breve q_{m+2n+2i}^{\,\, (n,2)} = \max \big\{ \breve \nu_{m,n} \le x \le \breve \mu_{m,n} : f^{m+2n+2i}(x) = x \big\} \,\,\, \text{and}
$$ 
$$
\,\,\,\, \breve \mu_{m,n,i}' = \max \big\{ \breve \nu_{m,n} \le x \le \breve \mu_{m,n} : f^{m+2n+2i}(x) = d \big\} \quad\quad\,\,\,
$$
\big(note the relationship between the subscript of $\breve \mu_{m,n,i}'$ and the superscript of $f^{m+2n+2i}$ in the definition of $\breve \mu_{m,n,i}'$.  Here the $n$ and $i$ in the subscript of $\breve \mu_{m,n,i}'$ indicate the $i^{th}$ point of the sequence $< \breve \mu_{m,n,i}' >$ in the $n^{th}$ interval $[\breve \mu_{m,n+1}, \breve \mu_{m,n}]$.  We need these points $\breve \mu_{m,n,i}'$'s to continue the {\it basic} tower-building process into the third layer\big).  By arguing as before, we obtain that 
$$
\breve \mu_{m,n+1} < \breve \nu_{m,n} < \breve \mu_{m,n,1}' < \breve q_{m+2n+2}^{\,\, (n,2)} < \breve \mu_{m,n,2}' < \breve q_{m+2n+4}^{\,\, (n,2)} < \breve \mu_{m,n,3}' < \breve q_{m+2n+6}^{\,\, (n,2)} < \cdots < \breve \mu_{m,n}.
$$
In the following, we shall show that, for each $n \ge 1$, with $x$-coordinates moving from point $\breve \mu_{m,n}$ via point $\breve \nu_{m,n}$ to point $\breve \mu_{m,n+1}$ (see Figure 4):
\begin{itemize}
\item[{\rm (a)}]
(odd periods) the point $\breve q_{m+2n+2i}^{\,\, (n,2)} = \max \big\{ \breve \mu_{m,n,1}' \le x \le \breve \mu_{m,n} : f^{m+2n+2i}(x) = x \big\} = \max \big\{ \breve \mu_{m,n+1} \le x \le \breve \mu_{m,n}: f^{m+2n+2i}(x) = x \big\}$ exists and is a period-$(m+2n+2i)$ point of $f$ for each $i \ge 1$;

\item[{\rm (b)}]
(even periods) the point $\breve c_{2m+2n+2i}'^{\,\, (n,2)} = \max \big\{ \breve \nu_{m,n} \le x \le \breve \mu_{m,n,1}': f^{2m+2n+2i}(x) = x \big\} = \max \big\{ \breve \mu_{m,n+1} \le x \le \breve \mu_{m,n,1}' : f^{2m+2n+2i}(x) = x \big\}$ exists for each $i \ge 1$ and is a period-$(2m+2n+2i)$ point of $f$ for each $i \ge n+3$; 

\item[{\rm (c)}]
(even periods) the point $\breve c_{2m+2n+2i}^{\,\, (n,2)} = \min \big\{ \breve \mu_{m,n+1} \le x \le \breve \nu_{m,n} : f^{2m+2n+2i}(x) = x \big\} = \min \big\{ \breve \mu_{m,n+1} \le x \le \breve \mu_{m,n} : f^{2m+2n+2i}(x) = x \big\}$ exists for each $i \ge 1$ and is a period-$(2m+2n+2i)$ point of $f$ for each $i \ge n+3$. 
\end{itemize}

\begin{figure}[htbp] 
\centering
\includegraphics[width=6in,height=2in]{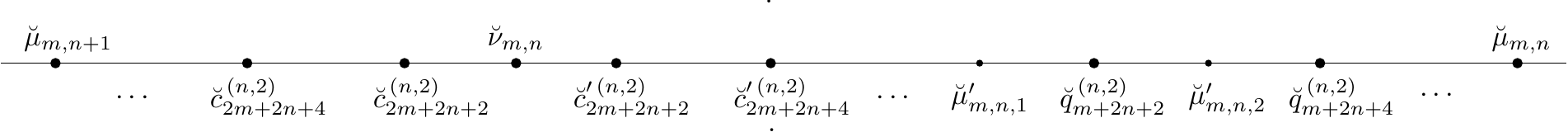} 

{Figure 4: A compartment of the second layer of the {\it basic} tower of periodic points of $f$ associated with $P$ in the interval $[\breve \mu_{m,n+1}, \breve \mu_{m,n}]$ $(\subset [\breve u_0', v] \subset [u_1, v]$), \\
where $\breve \nu_{m,n}$ is an {\it auxiliary} point in $[\breve \mu_{m,n+1}, \breve \mu_{m,n}]$ such that $f^{m+2n}(\nu_n) = v$ \\ 
which is needed only to determine the ordering of the following points: \\ 
and, for each $i \ge 1$, \\ 
$\breve \mu_{m,n,i}' = \max \big\{ \breve \mu_{m,n+1} \le x \le \breve \mu_{m,n} : f^{m+2n+2i}(x) = d \big\}$; \vspace{.05in} \\ 
$\breve q_{m+2n+2i}^{\, (n,2)} = \max \big\{ \breve \mu_{m,n+1} \le x \le \breve \mu_{m,n} : f^{m+2n+2i}(x) = x \big\}$; \vspace{.05in} \\ 
$\breve c_{2m+2n+2i}'^{\, (n,2)} = \max \big\{ \breve \mu_{m,n+1} \le x \le \breve \mu_{m,n,1}' : f^{2m+2n+2i}(x) = x \big\}$; \vspace{.05in} \\ 
$\breve c_{2m+2n+2i}^{\, (n,2)} = \min \big\{ \breve \mu_{m,n+1} \le x \le \breve \mu_{m,n} : f^{2m+2n+2i}(x) = x \big\}$. \vspace{.05in}} 
\end{figure}

For each fixed $n \ge 1$, the collection of all these periodic points $\breve q_{m+2n+2i}^{\,\, (n,2)}$, $\breve c_{2m+2n+2i}'^{\,\, (n,2)}$, $\breve c_{2m+2n+2i}^{\,\, (n,2)}$, $i \ge 1$, is called a {\it compartment} of the second layer of the {\it basic} tower of periodic points of $f$ associated with $P$.  
  
\noindent
{\bf 2.3(a) On the existence of periodic points of $f$ of all odd periods $\ge m+2n+2$ in $[\breve \mu_{m,n,1}', \breve \mu_{m,n}]$ $(\subset [\breve \mu_{m,n+1}, \breve \mu_{m,n}] \subset [\breve u_0', v] \subset [d, v])$ for each $n \ge 1$.}

We can apply Lemma 5(2) to obtain a rough result.  However, in this special case, we can get a result better than that obtained by applying Lemma 5(2).  So, we argue as follows:

We first prove the following result:

\noindent
{\bf Lemma 9.}
{\it For each $n \ge 1$, all periodic points of $f$ in $[\breve u_0', \breve \mu_{m,n}]$ $(\subset [\breve u_0', v])$ of odd periods have least periods $\ge m+2n$.  \big(Note that, in Lemma 7, we have $\ge m+2n+2$\big)}

\noindent
{\it Proof.} 
Suppose $f$ had a periodic point $\breve p$ of {\it odd} period $\ell(\breve p) \le m+2n-2$ in $[\breve u_0', \breve \mu_{m,n}]$ ($\subset [\breve u_0', v]$).  

Since $f(x) > z \ge z_0$ for all $d < x < z_0$ and since $\ell(\breve p)$ is odd, we have $\ell(\breve p) \ge 3$.  Consequently, $f^{\ell(\breve p)}(\breve u_0) \in \big\{ z_0, f(z_0), f(d) \big\}$. 

Since $m+2n-2 \ge \ell(\breve p) \ge 3$, by applying Lemma 3 with $f^3(\breve u_0') = f(d)$ and $f^{\ell(\breve p)}(\breve p) = \breve p$, there is a periodic point $\breve w_{m+2n-2}$ of $f$ in $[\breve u_0', \breve p]$ ($\subset [\breve u_0', \breve \mu_{m,n}]$) such that $f^{m+2n-2}(\breve w_{m+2n-2}) = \breve w_{m+2n-2}$.  Since $f^{m+2n-2}\big([\breve u_0', \breve w_{m+2n-2}]\big) \supset [\breve w_{m+2n-2}, z_0] \supset \{ v \}$, there is a point $\breve \nu_{m,n-1}$ in $[\breve u_0', \breve w_{m+2n-2}]$ such that $f^{m+2n-2}(\breve \nu_{m,n-1}) = v$.  Consequently, since 
$$
f^{m+2n}(\breve u_0') - d \ge z_0 - d > 0 \,\,\, \text{and} \,\,\, f^{m+2n}(\breve \nu_{m,n-1}) - d = f^2(v) - d = \min P - d < 0,
$$
there is a point $\breve \mu_{m,n}^{\, *}$ in $[\breve u_0', \breve \nu_{m,n-1}]$ ($\subset [\breve u_0', \breve \mu_{m,n}]$) such that $f^{m+2n}(\breve \mu_{m,n}^{\, *}) = d$.  Since $\breve u_0' < \breve \mu_{m,n}^{\, *} < \breve \nu_{m,n-1} < \breve \mu_{m,n}$, this contradicts the minimality of $\breve \mu_{m,n}$ in $[\breve u_0', v]$.  So, we have shown that all periodic points of $f$ in $[\breve u_0', \breve \mu_{m,n}]$ of odd periods have least periods $\ge m+2n$.
$\hfill\square$

Since the point $\breve p_{m+2n} = \min \big\{ \breve u_0' \le x \le v: f^{m+2n}(x) = x \big\}$ (defined on the first layer of the {\it basic} tower) exists and satisfies $\breve \mu_{m,n+1} < \breve p_{m+2n} < \breve \mu_{m,n}$, the point 
$$
\breve q_{m+2n} = \max \big\{ \breve \mu_{m,n+1} \le x \le \breve \mu_{m,n}: f^{m+2n}(x) = x \big\}
$$
exists and is $\ge \breve p_{m+2n}$.  Therefore, there is no period-$(m+2n)$ point of $f$ in $(\breve q_{m+2n}, \breve \mu_{m,n}]$.  This, combined with the above Lemma 9, implies that, on the interval $(\breve q_{m+2n}, \breve \mu_{m,n}]$, 
$$
\text{all periodic points of} \,\,\, f \,\,\, \text{of odd periods have least periods} \, \ge m+2n+2.
$$
Since, by the choice of $\breve u_0'$, we have \, $f^2(x) < d$ \, for all \, $\breve u_0' < x < \breve \mu_{m,n}$.  Since $f^{m+2n+2}(\breve q_{m+2n}) - d = f^2(\breve q_{m+2n}) - d < 0$ and $f^{m+2n+2}(\breve \mu_{m,n}) - d = z_0 - d > 0$, the point 
{\small 
$$
\breve \mu_{m,n,1}' = \max \big\{ \breve q_{m+2n} \le x \le \breve \mu_{m,n}: f^{m+2n+2}(x) = d \big\} = \max \big\{ \breve \mu_{m,n+1} \le x \le \breve \mu_{m,n}: f^{m+2n+2}(x) = d \big\}
$$}
satisfies $\breve q_{m+2n} < \breve \mu_{m,n,1}' < \breve \mu_{m,n}$.  Therefore, since $f^{m+2n+2}(\breve \mu_{m,n,1}') - \breve \mu_{m,n,1}' = d - \breve \mu_{m,n,1}' < 0$ and $f^{m+2n+2}(\breve \mu_{m,n}) - \breve \mu_{m,n} = z_0 - \breve \mu_{m,n} > 0$, the point $\breve q_{m+2n+2}^{\,\, (n,2)} = \max \big\{ \breve \mu_{m,n+1} \le x \le \breve \mu_{m,n}: f^{m+2n+2}(x) = x \big\}$ satisfies $\breve q_{m+2n} < \breve \mu_{m,n,1}' < \breve q_{m+2n+2}^{\,\, (n,2)} \,\, (< \breve q_{m+2n+4}^{\,\, (n,2)} < \breve q_{m+2n+6}^{\,\, (n,2)} < \cdots < \breve \mu_{m,n})$. 

Now suppose, for some $j \ge 1$, $\ell\big(\breve q_{m+2n+2j}^{\,\, (n,2)}\big) < m+2n+2j$.  By Lemma 1, $\ell\big(\breve q_{m+2n+2j}^{\,\, (n,2)}\big)$ divides $m+2n+2j$ and so, is odd.  It follows from the previous paragraph that $\ell\big(\breve q_{m+2n+2j}^{\,\, (n,2)}\big) \ge m+2n+2$.  Since $f^{m+2n+2}(\breve \mu_{m,n}) = z_0$ and $f^{\ell(\breve q_{m+2n+2j}^{\,\, (n,2)})}\big(\breve q_{m+2n+2j}^{\,\, (n,2)}\big) = \breve q_{m+2n+2j}^{\,\, (n,2)}$ and since $m+2n+2j > \ell\big(\breve q_{m+2n+2j}^{\,\, (n,2)}\big)$, it follows from Lemma 3 that there is a periodic point $\breve w_{m+2n+2j}$ of $f$ in $\big(\breve q_{m+2n+2j}^{\,\, (n,2)}, \breve \mu_{m,n}\big]$ such that $f^{m+2n+2j}(\breve w_{m+2n+2j}) = \breve w_{m+2n+2j}$.  Since $\breve q_{m+2n+2j}^{\,\, (n,2)} < \breve w_{m+2n+2j} < \breve \mu_{m,n}$, this contradicts the maximality of $\breve q_{m+2n+2j}^{\,\, (n,2)}$ in $[\breve \mu_{m,n+1}, \breve \mu_{m,n}]$.  Therefore, in $[\breve \mu_{m,n,1}', \breve \mu_{m,n}]$, 
$$
\text{for each} \,\,\, n \ge 1 \,\,\, \text{and} \,\,\, i \ge 1, \, \text{the point} \,\,\, \breve q_{m+2n+2i}^{\,\, (n,2)} \,\,\, \text{is a period-}(m+2n+2i) \,\,\, \text{point of} \,\,\, f.
$$

\noindent
{\bf 2.3(b) On the existence of periodic points of $f$ of all even periods $\ge 2m+4n+6$ in $[\breve \nu_{m,n}, \breve u_{m,n,1}']$ $(\subset [\breve \mu_{m,n+1}, \breve \mu_{m,n}]$ $\subset [\breve u_0', v] \subset [u_1, v])$ for each $n \ge 1$.}

We apply Lemma 6(2) with 
{\large 
$$
f^{m+2n+3}(\breve \mu_{m,n,1}') = f(d) \in \big\{ z_0, f(z_0), f(d) \big\} \,\,\, \text{and} \,\,\, f^{2(m+n+1)}(\breve \nu_{m,n}) = \min P
$$}
to obtain that, for each $i \ge 0$, the point 
\begin{multline*}
$$
\qquad\qquad\qquad c_{2(m+n+1)+2i}'^{\,\,(n,2)} = \max \big\{ \breve \nu_{m,n} \le x \le \breve \mu_{m,n,1}' : f^{2(m+n+1)+2i}(x) = x \big\} \\ 
= \max \big\{ \breve \mu_{m,n+1} \le x \le \breve \mu_{m,n,1}' : f^{2(m+n+1)+2i}(x) = x \big\}\qquad\qquad\quad\,
$$
\end{multline*}
exists and $\breve \nu_{m,n} \, < \, \breve c_{2m+2n+2}'^{\,\,(n,2)} \, < \, \breve c_{2m+2n+4}'^{\,\,(n,2)} \, < \, \breve c_{2m+2n+6}'^{\,\,(n,2)} \, < \, \cdots \, < \, \breve \mu_{m,n,1}' \, < \, \breve \mu_{m,n} \, < \, v$.  Furthermore, for each $i \ge 0$ such that $(m+n+1)+i \ge \max \big\{m+2n+3, m+2 \big\} = m+2n+3$, i.e., for each $i \ge n+2$, the point $\breve c_{2(m+n+1)+2i}'^{\,\,(2)}$ is a period-$\big(2(m+n+1)+2i\big)$ point of $f$, or equivalently, in $[\breve \nu_{m,n}, \breve u_{m,n,1}']$ $(\subset [\breve \mu_{m,n+1}, \breve \mu_{m,n}]$ $\subset [\breve u_0', v] \subset [u_1, v])$, 
$$
\text{for each} \,\,\, n \ge 1 \,\,\, \text{and} \,\,\, i \ge n+3, \, \text{the point} \,\,\, c_{2m+2n+2i}'^{\,\,(n,2)} \,\,\, \text{is a period-}(2m+2n+2i) \,\,\, \text{point of} \,\,\,f.
$$  

\noindent
{\bf 2.3(c) On the existence of periodic points of $f$ of all even periods $\ge 2m+4n+6$ in $[\breve \mu_{m,n+1}, \breve \nu_{m,n}]$ $(\subset [\breve \mu_{m,n+1}, \breve \mu_{m,n}]$ $\subset [\breve u_0', v] \subset [u_1, v])$ for each $n \ge 1$.}

We apply Lemma 6(1) with 
{\large 
$$
f^{m+2n+3}(\breve \mu_{m,n+1}) = f(d) \in \big\{ z_0, f(z_0), f(d) \big\} \,\,\, \text{and} \,\,\, f^{2(m+n+1)}(\breve \nu_{m,n}) = \min P
$$}
to obtain that, for each $i \ge 0$, the point 
\begin{multline*}
$$
\qquad\qquad\qquad c_{2(m+n+1)+2i}^{\,\,(n,2)} = \min \big\{ \breve \mu_{m,n+1} \le x \le \breve \nu_{m,n}: f^{2(m+n+1)+2i}(x) = x \big\} \\ 
= \min \big\{ \breve \mu_{m,n+1} \le x \le \breve \mu_{m,n} : f^{2(m+n+1)+2i}(x) = x \big\}\qquad\qquad\quad\,\,\,\,\,
$$
\end{multline*}
exists and $\breve \mu_{m,n+1} \, < \, \cdots \, < \, \breve c_{2m+2n+6}^{\,\,(n,2)} \, < \, \breve c_{2m+2n+4}^{\,\,(n,2)} \, < \, \breve c_{2m+2n+2}^{\,\,(n,2)} \, < \, \breve \nu_{m,n} \, < \, \breve \mu_{m,n} \, < \, v$.  Furthermore, for each $i \ge 0$ such that $(m+n+1)+i \ge \max \big\{m+2n+3, m+2 \big\} = m+2n+3$, i.e., for each $i \ge n+2$, the point $\breve c_{2(m+n+1)+2i}^{\,\,(n,2)}$ is a period-$\big(2(m+n+1)+2i\big)$ point of $f$, or equivalently, in $[\breve \mu_{m,n+1}, \breve \nu_{m,n}]$ $(\subset [\breve \mu_{m,n+1}, \breve \mu_{m,n}]$ $\subset [\breve u_0', v] \subset [u_1, v])$.  
$$
\text{for each} \,\,\, n \ge 1 \,\,\, \text{and} \,\,\, i \ge n+3, \, \text{the point} \,\,\, \breve c_{2m+2n+2i}^{\,\,(n,2)} \,\,\, \text{is a period-}(2m+2n+2i) \,\,\, \text{point of} \,\,\, f.
$$
$$\aleph \qquad\qquad\qquad \aleph \qquad\qquad\qquad \aleph \qquad\qquad\qquad \aleph \qquad\qquad\qquad \aleph$$
\indent Recall that, on $[v, z_0]$, we define $\hat u_1' = \max \big\{ v \le x < z_0: f^2(x) = x \big\}$.  However, for our purpose, we shall use the point 
$$
\hat u_0 = \min \big\{ v \le x \le z_0: f^2(x) = d \big\} \, (\le \bar u_1')
$$
instead of the point $\bar u_1'$ and consider the interval $[v, \hat u_0]$ because on this interval, we have 
{\large 
$$
f^2(x) < d \,\,\, \text{for all} \,\,\, v < x < \hat u_0.
$$}
\indent Since $f^{m+2}\big([v, \hat u_0]\big) \supset [\min P, f(z_0)] \supset [\min P, z_0]) \supset \{ d \}$, the point 
$$
\hat \mu_{m,1}' = \max \big\{ v \le x \le \hat u_0: f^{m+2}(x) = d \big\}
$$
exists.  Since $f^{m+4}\big([\hat \mu_{m,1}', \hat u_0]\big) \supset f^2\big([d, z_0]\big) \supset f^2\big([v, z_0]\big) \supset [\min P, z_0] \supset \{ d \}$, the point 
$$
\hat \mu_{m,2}' = \max \big\{ \hat \mu_{m,1}' \le x \le \hat u_0: f^{m+4}(x) = d \big\} = \max \big\{ v \le x \le \hat u_0: f^{m+4}(x) = d \big\}
$$
exists and is $> \hat \mu_{m,1}'$.  Inductively, for each $n \ge 1$, let \big(note the relationship between the subscript of $\hat \mu_{m,n}'$ and the superscript of $f^{m+2n}$ in the definition of $\hat \mu_{m,n}'$\big) 
$$
\hat \mu_{m,n}' = \max \big\{ v \le x \le \hat u_0: f^{m+2n}(x) = d \big\}.
$$
Then we have $\,\,\, v \, < \, \hat \mu_{m,1}' \, < \, \hat \mu_{m,2}' \, < \, \hat \mu_{m,3}' \, < \, \cdots \, < \, \hat u_0$.  For each $n \ge 1$, let $\hat \nu_{m,n}$ be a point in $[\hat \mu_{m,n}', \hat u_0]$ such that $f^{m+2n}(\hat \nu_{m,n}) = v$.  Then it turns out that $\hat \mu_{m,n}' < \hat \nu_{m,n} < \hat \mu_{m,n+1}'$.  Furthermore, on the interval $[\hat \mu_{m,n}', \hat \mu_{m,n+1}'] \subset [v, \hat u_0]$, we have 
$$
f^{m+2n+2}(\hat \mu_{m,n}') = z_0, \, f^{m+2n+2}(\hat \nu_{m,n}) = \min P, \, f^{m+2n+2}(\hat \mu_{m,n+1}') = d \,\,\, \text{and}
$$
$$
\text{the graph of} \,\,\, y = f^{m+2n+2}(x) \,\,\, \text{looks 'roughly' like a {\it skew} 'V' shape at the points} \,\,\, \hat \mu_{m,n}',
$$ 
$\hat \nu_{m,n}, \hat \mu_{m,n+1}'$ \big(meaning that the graph of $y = f^{m+2n+2}(x)$ passes through the 3 points $\, (\hat \mu_{m,n}', z_0)$, $(\hat \nu_{m,n}, \min P), (\hat \mu_{m,n+1}', d)$ on the $x$-$y$ plane with $y$-coordinates $z_0, \min P, d$ respectively\big) just like that of $y = f^2(x)$ at the points $d, v, \hat u_0$ on the interval $[d, \hat u_0]$.  

Now, for each $n \ge 1$ and all $i \ge 1$, let \big(the number 2 in the superscript $(n,2)$ indicates the second layer\big) 
$$
\quad\,\,\, \hat p_{m+2n+2i}^{\,\, (n,2)} = \min \big\{ \hat \mu_{m,n}' \le x \le \hat \mu_{m,n+1}' : f^{m+2n+2i}(x) = x \big\} \, (< \hat \nu_{m,n})\,\, \text{and}
$$ 
$$
\hat \mu_{m,n,i} = \min \big\{ \hat \mu_{m,n}' \le x \le \hat \mu_{m,n+1}' : f^{m+2n+2i}(x) = d \big\} \, (< \hat \nu_{m,n})\quad\quad
$$
\big(note the relationship between the subscript of $\hat \mu_{m,n,i}$ and the superscript of $f^{m+2n+2i}$ in the definition of $\hat \mu_{m,n,i}$.  Here the $n$ and $i$ in the subscript of $\hat \mu_{m,n,i}$ indicate the $i^{th}$ point of the sequence $< \hat \mu_{m,n,i} >$ in the $n^{th}$ interval $[\hat \mu_{m,n}', \hat \mu_{m,n+1}']$.  We need these points $\hat \mu_{m,n,i}$'s to continue the {\it basic} tower-building process into the third layer\big).  Note that, by the choice of the point $\hat u_0$, we have $f^2(x) < d$ for all $v < x < \hat u_0$.  By arguing as before, we obtain that 
$$
\hat \mu_{m,n}' \, < \, \cdots \, < p_{m+2n+6}^{\,\, (n,2)} \, < \, \hat \mu_{m,n,3} \, < \, \hat p_{m+2n+4}^{\,\, (n,2)} \, < \, \hat \mu_{m,n,2} \, < \hat p_{m+2n+2}^{\,\, (n,2)} \, < \, \hat \mu_{m,n,1} \, < \, \hat \nu_{m,n} \, < \, \hat \mu_{m,n+1}'.
$$
\indent In the following, we shall show that, for each $n \ge 1$, with $x$-coordinates moving from point $\hat \mu_{m,n}'$ via point $\hat \nu_{m,n}$ to point $\hat \mu_{m,n+1}'$ (see Figure 4):

\pagebreak

\begin{itemize}
\item[{\rm (a)}]
(odd periods) the point $\hat p_{m+2n+2i}^{\,\, (n,2)} = \min \big\{ \hat \mu_{m,n}' \le x \le \hat \mu_{m,n,1}: f^{m+2n+2i}(x) = x \big\} = \min \big\{ \hat \mu_{m,n}' \le x \le \hat \mu_{m,n+1}': f^{m+2n+2i}(x) = x \big\}$ exists and is a period-$(m+2n+2i)$ point of $f$ for each $i \ge 1$;

\item[{\rm (b)}]
(even periods) the point $\hat c_{2m+2n+2i}^{\,\, (n,2)} = \min \big\{ \hat \mu_{m,n,1} \le x \le \hat \nu_{m,n}: f^{2m+2n+2i}(x) = x \big\} = \min \big\{ \hat \mu_{m,n,1} \le x \le \hat \mu_{m,n+1}': f^{2m+2n+2i}(x) = x \big\}$ exists for each $i \ge 1$ and is a period-$(2m+2n+2i)$ point of $f$ for each $i \ge n+3$; 

\item[{\rm (c)}]
(even periods) the point $\hat c_{2m+2n+2i}'^{\,\, (n,2)} = \max \big\{ \hat \nu_{m,n} \le x \le \hat \mu_{m,n+1}' : f^{2m+2n+2i}(x) = x \big\} = \min \big\{ \hat \mu_{m,n}' \le x \le \hat \mu_{m,n+1}' : f^{2m+2n+2i}(x) = x \big\}$ exists for each $i \ge 1$ and is a period-$(2m+2n+2i)$ point of $f$ for each $i \ge n+3$. 
\end{itemize}

\begin{figure}[htbp] 
\centering
\includegraphics[width=6in,height=2in]{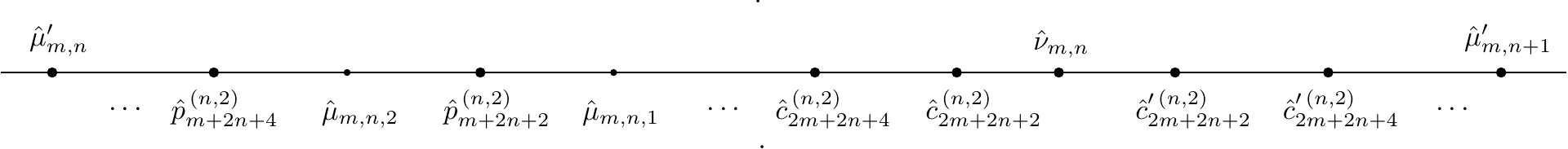} 

{Figure 5: A compartment of the second layer of the {\it basic} tower of periodic points of $f$ associated with $P$ in the interval $[\hat \mu_{m,n}', \hat \mu_{m,n+1}']$ $(\subset [v, \hat u_0])$, \\  
where $\hat \nu_{m,n}$ is an {\it auxiliary} point in $[\hat \mu_{m,n}', \hat \mu_{m,n+1}']$ such that $f^{m+2n}(\hat \nu_{m,n}) = v$ \\ 
which is needed only to determine the ordering of the following points: \\ 
and, for each $i \ge 1$, \\ 
$\hat \mu_{m,n,i} = \min \big\{ \hat \mu_{m,n}' \le x \le \hat \mu_{m,n+1}' : f^{m+2n+2i}(x) = d \big\}$; \vspace{.05in} \\ 
$\hat p_{m+2n+2i}^{\,(n,2)} = \min \big\{ \hat \mu_{m,n}' \le x \le \hat \mu_{m,n+1}' : f^{m+2n+2i}(x) = x \big\}$; \vspace{.05in} \\ 
$\hat c_{2m+2n+2i}^{\, (n,2)} = \min \big\{ \hat \mu_{m,n,1} \le x \le \hat \mu_{m,n+1}' : f^{2m+2n+2i}(x) = x \big\}$; \vspace{.05in} \\ 
$\hat c_{2m+2n+2i}'^{\, (n,2)} = \max \big\{ \hat \mu_{m,n}' \le x \le \hat \mu_{m,n+1}' : f^{2m+2n+2i}(x) = x \big\}$. \vspace{.05in}}  
\end{figure}

For each fixed $n \ge 1$, the collection of all these periodic points $\hat p_{m+2n+2i}^{\,\, (n,2)}$, $\hat c_{2m+2n+2i}^{\,\, (n,2)}$, $\hat c_{2m+2n+2i}'^{\,\, (n,2)}$, $i \ge 1$, is called a {\it compartment} of the second layer of the {\it basic} tower of periodic points of $f$ associated with $P$.  

\noindent
{\bf 2.4(a) On the existence of periodic points of $f$ of all odd periods $\ge m+2n+2$ in $[\hat \mu_{m,n}', \hat \mu_{m,n,1}]$ $(\subset [\hat \mu_{m,n}', \mu_{m,n+1}'] \supset [v, \hat u_0]$ $ \subset [v, z_0])$ for each $n \ge 1$.}

We can apply Lemma 5(2) to obtain a rough result.  However, in this special case, we can get a result better than that obtained by applying Lemma 5(2).  So, we argue as follows:

We first prove the following result:

\noindent
{\bf Lemma 10.}
{\it For each $n \ge 1$, all periodic points of $f$ in $[\hat \mu_{m,n}', \hat u_0]$ $(\subset [v, \hat u_0])$ of odd periods have least periods $\ge m+2n$.  (Note that, in Lemma 7, we have $\ge m+2n+2$)}

\noindent
{\it Proof.}
Suppose $f$ had a periodic point $\hat p$ of {\it odd} period $\ell(\hat p) \le m+2n-2$ in $(\hat \mu_{m,n}', \hat u_0]$ ($\subset [v, \hat u_0]$).  

Since $f(x) > z \ge z_0$ for all $d < x < z_0$ and since $\ell(\hat p)$ is odd, we have $\ell(\hat p) \ge 3$.  Consequently, $f^{\ell(\hat p)}(\hat u_0) \in \big\{ z_0, f(z_0), f(d) \big\}$. 

Since $m+2n-2 \ge \ell(\hat p) \ge 3$, by Lemma 3 with $f^3(\hat u_0) = f(d)$ and $f^{\ell(\hat p)}(\hat p) = \hat p$, there is a point $\hat w_{m+2n-2}$ of $f$ in $[\hat p, \hat u_0]$ ($\subset (\hat \mu_{m,n}', \hat u_0] \subset [v, \hat u_0]$) such that $f^{m+2n-2}(\hat w_{m+2n-2}) = \hat w_{m+2n-2}$.  Since $f^2(x) < d$ for all $v < x < \hat u_0$, we have $f^{m+2n}\big([\hat w_{m+2n-2}, \hat u_0]\big) \supset \big[f^2(\hat w_{m+2n-2}), z_0\big] \supset [d, z_0] \supset \{ d \}$.  So, there is a point $\hat \mu_{m,n}^{\, *}$ in $[\hat w_{m+2n-2}, \hat u_0]$ ($\subset (\hat \mu_{m,n}', \hat u_0]$) such that $f^{m+2n}(\hat \mu_{m,n}^{\, *}) = d$.  Since $\hat \mu_{m,n}' < \hat w_{m+2n-2} < \hat \mu_{m,n}^{\, *} < \hat u_0$, this contradicts the maximality of $\hat \mu_{m,n}'$ in $[v, \hat u_0]$.  Therefore, we have shown that all periodic points of $f$ in $[\hat \mu_{m,n}', \hat u_0]$ of {\it odd} periods have least periods $\ge m+2n$.
$\hfill\square$

Now suppose, for some $j \ge 1$, the least period $\ell\big(\hat p_{m+2n+2j}^{\,\, (n,2)}\big)$ of $\hat p_{m+2n+2j}^{\,\, (n,2)}$ is $< m+2n+2j$.  By Lemma 1, $\ell\big(\hat p_{m+2n+2j}^{\,\, (n,2)}\big)$ divides $m+2n+2j$ and so, is odd.

If $\ell\big(\hat p_{m+2n+2j}^{\,\, (n,2)}\big) \ge m+2n+2$, then, since $f^{m+2n+2}(\hat \mu_{m,n}') = z_0$ and $f^{\ell(\hat p_{m+2n+2j}^{\,\, (n,2)})}\big(\hat p_{m+2n+2j}^{\,\, (n,2)}\big) = \hat p_{m+2n+2j}^{\,\, (n,2)}$ and since $m+2n+2j > \ell\big(\hat p_{m+2n+2j}^{\,\, (n,2)}\big)$, it follows from Lemma 3 that there is a periodic point $\hat w_{m+2n+2j}$ of $f$ in $\big[\hat \mu_{m,n}', \hat p_{m+2n+2j}^{\,\, (n,2)}\big)$ such that $f^{m+2n+2j}(\hat w_{m+2n+2j}) = \hat w_{m+2n+2j}$.  Since $\hat \mu_{m,n}' < \hat w_{m+2n+2j} < \hat p_{m+2n+2j}^{\,\, (n,2)} < \hat \mu_{m,n+1}'$, this contradicts the minimality of $\hat p_{m+2n+2j}^{\,\, (n,2)}$ in $[\hat \mu_{m,n}', \hat \mu_{m,n+1}']$.  

If $\ell\big(\hat p_{m+2n+2j}^{\,\, (n,2)}\big) = m+2n$, then, let $r = (m+2n+2j)/\ell\big(\hat p_{m+2n+2j}^{\,\, (n,2)}\big)$.  Suppose $r > 3$.  Then $f^{3\ell(\hat p_{m+2n+2j}^{\,\, (n,2)})}\big(\hat p_{m+2n+2j}^{\,\, (n,2)}\big) = \hat p_{m+2n+2j}^{\,\, (n,2)}$ and $f^{3\ell(\hat p_{m+2n+2j}^{\,\, (n,2)})}(\hat \mu_{m,n}') = f^{3(m+2n)}(\hat \mu_{m,n}') = z_0$.  Since $m+2n+2j > 3\ell\big(\hat p_{m+2n+2j}^{\,\, (n,2)}\big)$, it follows from Lemma 3 that there is a periodic point $\hat w_{m+2n+2j}^*$ of $f$ in $\big(\hat \mu_{m,n}', \hat p_{m+2n+2j}^{\,\, (n,2)}\big]$ such that $f^{m+2n+2j}(\hat w_{m+2n+2j}^*) = \hat w_{m+2n+2j}^*$ which contradicts the minimality of $\hat p_{m+2n+2j}^{\,\, (n,2)}$ in $[\hat \mu_{m,n}', \hat \mu_{m,n+1}']$.  Therefore, if $m+2n = \ell\big(\hat p_{m+2n+2j}^{\,\, (n,2)}\big) < m+2n+2j$, then $m+2n+2j = 3(m+2n)$.  That is, the point $\hat p_{3(m+2n)}^{\,\, (n,2)}$ is either a period-$(m+2n)$ point or a period-$\big(3(m+2n)\big)$ point of $f$.  

We now show that the point $\hat p_{3(m+2n)}^{\,\, (n,2)}$ is actually a period-$\big(3(m+2n)\big)$ point of $f$.  Since 
$$
f^{m+2n}(\hat \mu_{m,n}') - \hat \mu_{m,n}' = d - \hat \mu_{m,n}' < 0 \,\,\, \text{and} \,\,\, \hat f^{m+2n}(\hat u_0) - \hat u_0 \ge z_0 - \hat u_0 > 0,
$$
the point $\hat p_{m+2n}^* = \min \big\{ \hat \mu_{m,n}' \le x \le \hat u_0: f^{m+2n}(x) = x \big\}$ exists.  By the choice of $\hat u_0$, we have $f^2(x) < d$ for all $v < x < \hat u_0$.  Since 
$$
f^{m+2n+2}(\hat p_{m+2n}^*) - d = f^2(\hat p_{m+2n}^*) - d < 0 \,\,\, \text{and} \,\,\, f^{m+2n+2}(\hat u_0) - d \ge z_0 - d > 0,
$$
we obtain that $\hat \mu_{m,n}' < \hat p_{m+2n}^* < \hat \mu_{m,n+1}'$.  Recall that we have $f^2(x) < x$ for all $v < x < z_0$.  Therefore, we have 
$$
f^{m+2n+2}(\hat p_{m+2n}^*) - \hat p_{m+2n}^* = f^2(\hat p_{m+2n}^*) - \hat p_{m+2n}^* < 0 \,\,\, \text{and} \,\,\, f^{m+2n+2}(\hat \mu_{m,n}') - \hat \mu_{m,n}' = z_0 - \hat \mu_{m,n}' > 0.
$$
Consequently, the point $\hat p_{m+2n+2}^{\,\, (n,2)} = \min \big\{ \hat \mu_{m,n}' \le x \le \hat \mu_{m,n+1}': f^{m+2n+2}(x) = x \big\} = \min \big\{ \hat \mu_{m,n}' \le x \le \hat p_{m+2n}^*: f^{m+2n+2}(x) = x \big\}$ is $< \hat p_{m+2n}^*$.  Since $3(m+2n) > m+2n+2$, it follows from Lemma 3 with 
$$
f^{m+2n+2}(\hat \mu_{m,n}') = z_0 \in \big\{ z_0, f(z_0), f(d) \big\} \,\,\, \text{and} \,\,\, f^{m+2n+2}\big(\hat p_{m+2n+2}^{\,\, (n,2)}\big) =  \hat p_{m+2n+2}^{\,\, (n,2)}
$$
that $\hat \mu_{m,n}' < \cdots < \hat p_{3(m+2n)}^{\,\, (n,2)} < \cdots < \hat p_{m+2n+4}^{\,\, (n,2)} < \hat p_{m+2n+2}^{\,\, (n,2)}$ $(< \hat p_{m+2m}^*)$.  Since the point $\hat p_{m+2m}^*$ is the {\it smallest} point in $[\hat \mu_{m,n}', \hat \mu_{m,n+1}']$ that satisfies $f^{m+2n}(x) = x$ and $\hat p_{3(m+2n)}^{\,\, (n,2)} < \hat p_{m+2m}^*$, we see that $\hat p_{3(m+2n)}^{\,\, (n,2)}$ cannot be a period-$(m+2n)$ point of $f$.  Therefore, $\hat p_{3(m+2n)}^{\,\, (n,2)}$ is a period-$\big(3(m+2n)\big)$ point of $f$.   This, combined with the above, implies that, in $[\hat \mu_{m,n}', \hat \mu_{m,n,1}]$ $(\subset [\hat \mu_{m,n}', \mu_{m,n+1}'] \supset [v, \hat u_0]$ $ \subset [v, z_0])$, 
$$
\text{for each} \,\,\, n \ge 1 \,\,\, \text{and} \,\,\, i \ge 1, \, \text{the point} \,\,\, \hat p_{m+2n+2i}^{\,\, (n,2)} \,\,\, \text{is a period-}(m+2n+2i) \,\,\, \text{point of} \,\,\, f.
$$

\noindent
{\bf 2.4(b) On the existence of periodic points of $f$ of all even periods $\ge 4m+4n+6$ in $[\hat \mu_{m,n,1}, \hat \nu_{m,n}]$ $(\subset [\hat \mu_{m,n}', \hat \mu_{m,n+1}']$ $\subset [v, \hat u_0]$ $\subset [v, z_0])$ for each $n \ge 1$.}  

We apply Lemma 6(1) with {\large 
$$
f^{m+2n+3}(\hat \mu_{m,n,1}) = f(d) \in \big\{ z_0, f(z_0), f(d) \big\} \,\,\, \text{and} \,\,\, f^{2(m+n+1)}(\hat \nu_{m,n}) = \min P
$$}
to obtain that, for each $i \ge 0$, the point 
\begin{multline*}
$$
\qquad\qquad\qquad \hat c_{2(m+n+1)+2i}^{\,\,(n,2)} = \min \big\{ \hat \mu_{m,n,1} \le x \le \hat \nu_{m,n} : f^{2(m+n+1)+2i}(x) = x \big\} \\ 
= \min \big\{ \hat \mu_{m,n,1} \le x \le \hat \mu_{m,n+1}: f^{2(m+n+1)+2i}(x) = x \big\}\qquad\qquad\quad\,\,
$$
\end{multline*}
exists and $\hat \mu_{m,n}' \, < \, \hat \mu_{m,n,1} \, < \, \cdots \, < \, \hat c_{2m+2n+6}^{\,\,(n,2)} \, < \, \hat c_{2m+2n+4}^{\,\,(n,2)} \, < \, \hat c_{2m+2n+2}^{\,\,(n,2)} \, < \hat \nu_{m,n} \, < \, \hat \mu_{m,n+1}' \, < \, \hat u_0$.  Furthermore, for each $i \ge 0$ such that $(m+n+1)+i \ge  m + (m+2n+3) = 2m+2n+3$, i.e., for each $i \ge m+n+2$, the point $\hat c_{2(m+n+1)+2i}^{\,\,(n,2)}$ is a period-$\big(2(m+n+1)+2i\big)$ point of $f$, or, equivalently, in $[\hat u_{m,n,1}, \hat \nu_{m,n}]$ $(\subset [\hat \mu_{m,n}', \hat \mu_{m,n+1}']$ $\subset [v, \hat u_0] \subset [v, z_0])$, 
$$
\text{for each} \,\,\, n \ge 1 \,\,\, \text{and} \,\,\, i \ge m+n+3, \, \text{the point} \,\,\, \hat c_{2m+2n+2i}^{\,\,(n,2)} \,\,\, \text{is a period-}(2m+2n+2i) \,\,\, \text{point of} \,\,\, f.
$$ 

\noindent
{\bf 2.4(c) On the existence of periodic points of $f$ of all even periods $\ge 4m+4n+6$ in $[\hat \nu_{m,n}, \hat u_{m,n+1}']$ $(\subset [\hat \mu_{m,n}', \hat \mu_{m,n+1}']$ $\subset [v, \hat u_0]$ $\subset [v, z_0])$ for each $n \ge 1$.}  

Recall that $\hat \nu_{m,n}$ is a point in $[\hat \mu_{m,n}', \hat \mu_{m,n+1}']$ such that $f^{m+2n}(\hat \nu_{m,n}) = v$.

We apply Lemma 6(2) with 
{\large 
$$
f^{m+2n+3}(\hat \mu_{m,n+1}') = f(d) \in \big\{ z_0, f(z_0), f(d) \big\}, \,\,\, \text{and} \,\,\, f^{2(m+n+1)}(\hat \nu_{m,n}) = \min P
$$}
to obtain that, for each $i \ge 0$, the point 
\begin{multline*}
$$
\qquad\qquad\qquad \hat c_{2(m+n+1)+2i}'^{\,\,(n,2)} = \max \big\{ \hat \nu_{m,n} \le x \le \hat \mu_{m,n+1}': f^{2(m+n+1)+2i}(x) = x \big\} \\ 
= \max \big\{ \hat \mu_{m,n}' \le x \le \hat \mu_{m,n+1}': f^{2(m+n+1)+2i}(x) = x \big\}\qquad\qquad\quad\,\,\,\,
$$
\end{multline*}
exists and $v \, < \, \hat \mu_{m,n}' \, < \, \hat \nu_{m,n} \, < \, \hat c_{2m+2n+2}'^{\,\,(n,2)} \, < \, \hat c_{2m+2n+4}'^{\,\,(n,2)} \, < \, \hat c_{2m+2n+6}'^{\,\,(n,2)} \, < \, \cdots \, < \, \hat \mu_{m,n+1}' \, < \, \hat u_0$.  Furthermore, for each $i \ge 0$ such that $(m+n+1)+i \ge  m + (m+2n+3) = 2m+2n+3$, i.e., for each $i \ge m+n+2$, the point $\hat c_{2(m+n+1)+2i}'^{\,\,(n,2)}$ is a period-$\big(2(m+n+1)+2i\big)$ point of $f$, or, equivalently, in $[\hat \nu_{m,n}, \hat u_{m,n,1}]$ $(\subset [\hat \mu_{m,n}', \hat \mu_{m,n+1}']$ $\subset [v, \hat u_0] \subset [v, z_0])$,
$$
\text{for each} \,\,\, n \ge 1 \,\,\, \text{and} \,\,\, i \ge m+n+3, \, \text{the point} \,\,\, \hat c_{2m+2n+2i}'^{\,\,(n,2)} \,\,\, \text{is a period-}(2m+2n+2i) \,\,\, \text{point of} \,\,\, f.
$$  
$$\aleph \qquad\qquad\qquad \aleph \qquad\qquad\qquad \aleph \qquad\qquad\qquad \aleph \qquad\qquad\qquad \aleph$$
\indent 'Symmetrically', in the interval $[\bar u_1', z_0]$, for each $n \ge 1$, we can find similar sequences on $[\hat u_n', \hat u_{n+1}']$ {\it with some variations}.

Let $\hat u_0 = \min \big\{ v \le x \le z_0: f^2(x) = d \big\}$ and, for each $n \ge 1$, let $\bar u_n' = \max \big\{ v \le x \le z_0 : f^{2n}(x) = d \big\}$ be defined as before.  Then we have $v < \hat u_0 \le \bar u_1' < \bar u_2' < \bar u_3' < \cdots < z_0$, $f^2(z_0) = z_0$ and $f^{2n}(\bar u_n') = d$.  In particular, $f^{2n}\big([\bar u_n', z_0]\big) \supset [d, z_0] \supset \{ v \}$.  Let $\bar \nu_n$ be a point in $[\bar u_n', z_0]$ such that $f^{2n}(\bar \nu_n) = v$.  Then, since $f^2(v) = \min P$, we obtain that $f^{2n+2}\big([\bar \nu_n, z_0]\big) \supset [\min P, z_0] \supset \{ d \}$.  So, the point $\bar \nu_n$ happens to satisfy that $\bar u_{n}' < \bar \nu_n < \bar u_{n+1}'$.  On the interval $[\bar u_{n}', \bar u_{n+1}']$, we have 
$$
f^{2n+2}(\bar u_{n}') = z_0, \,\, f^{2n+2}(\bar \nu_n) = \min P, \,\, f^{2n+2}(\bar u_{n+1}') = d \,\,\, \text{and}
$$ 
$$
\text{the graph of} \,\, y = f^{2n+2}(x) \,\, \text{looks 'roughly' like a {\it skew} 'V' shape (at the 3 points} \,\, \bar u_{n}', \bar \nu_n, \bar u_{n+1}')
$$
\big(meaning that the graph of $y = f^{2n+2}(x)$ passes through the 3 points $(\bar u_n', z_0)$, $(\bar \nu_n, \min P)$, $(\bar u_{n+1}', d)$ on the $x$-$y$ plane with $y$-coordinates $z_0, \min P, d$ respectively\big) just like that of $y = f^2(x)$ at the 3 points $d, v, \hat u_0$ on the interval $[d, \hat u_0]$.  

In the following, we show that, for each $n \ge 1$, there exist 4 sequences of points in $[\bar u_{n}', \bar u_{n+1}']$ \big(the 2 in the supserscript $(n, 2)$ below indicates second layer\big):
\begin{multline*}
$$
\bar u_n' < \cdots < \bar u_{n,3} < \bar u_{n,2} < \bar u_{n,1} < \bar \nu_n < \bar u_{n+1}' \,\,\, \text{and} \\ 
\bar u_{n}' < \cdots < \bar c_{2n+6}^{\,\,(n,2)} < \bar c_{2n+4}^{\,\,(n,2)} < \bar c_{2n+2}^{\,\,(n,2)} < \bar u_{n,1} < \cdots < \bar p_{m+2n+6}^{\,\,(n,2)} < \bar p_{m+2n+4}^{\,\,(n,2)} < \bar p_{m+2n+2}^{\,\,(n,2)} < \bar \nu_n \\ 
< \bar q_{m+2n+2}^{\,\,(n,2)} < \bar q_{m+2n+4}^{\,\,(n,2)} < \bar q_{m+2n+6}^{\,\,(n,2)} < \cdots < \bar u_{n+1}'
$$
\end{multline*}
such that, for each $i \ge 1$, $\bar u_{n,i} = \min \big\{ \bar u_n' \le x \le \bar u_{n+1}' : f^{2n+2i}(x) = d \big\}$ \big(note the relationship between the subscript of $\bar u_{n,i}$ and the superscript of $f^{2n+2i}$, we need these points $\bar u_{n,i}$'s to continue the {\it basic} tower-building process to the third layer\big) and, with $x$-coordinates moving from point $\bar u_n'$ via point $\bar \nu_n$ to point $\bar u_{n+1}'$ (see Figure 6):
\begin{itemize} 
\item[{\rm (a)}]
(even periods) the point $\bar c_{2n+2i}^{\,\, (n,2)} = \min \big\{ \bar u_n' \le x \le \bar \nu_n : f^{2n+2i}(x) = x \big\} = \min \big\{ \bar u_n' \le x \le \bar u_{n+1}': f^{2n+2i}(x) = x \big\}$ exists and is a period-$(2n+2i)$ point of $f$ for each $i \ge 1$;

\item[{\rm (b)}]
(odd periods) the point $\bar p_{m+2n+2i}^{\,\, (n,2)} = \min \big\{ \bar u_{n,1} \le x \le \bar \nu_n : f^{m+2n+2i}(x) = x \big\} = \min \big\{ \bar u_{n,1} \le x \le \bar u_{n+1}': f^{m+2n+2i}(x) = x \big\}$ exists and is a period-$(m+2n+2i)$ point of $f$ for each $i \ge 1$;

\item[{\rm (c)}]
(odd periods) the point $\bar q_{m+2n+2i}^{\,\, (n,2)} = \max \big\{ \bar \nu_{n} \le x \le \bar u_{n+1}' : f^{m+2n+2i}(x) = x \big\} = \max \big\{ \bar u_n' \le x \le \bar u_{n+1}' : f^{m+2n+2i}(x) = x \big\}$ exists and is a period-$(m+2n+2i)$ point of $f$ for each $i \ge 1$.
\end{itemize}

\begin{figure}[htbp] 
\centering
\includegraphics[width=6in,height=2in]{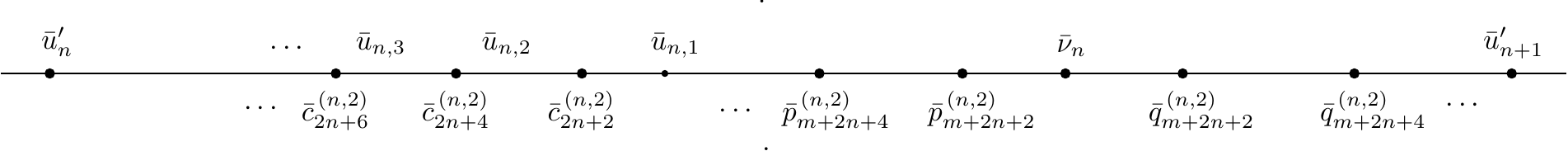} 

{Figure 6: A compartment of the second layer of the {\it basic} tower of periodic points of $f$ associated with $P$ in the interval $[\bar u_n', \bar u_{n+1}']$ $(\subset [\bar u_1', z_0])$, \\  
where $\bar \nu_n$ is an {\it auxiliary} point in $[\bar u_n', \bar u_{n+1}']$ such that $f^{2n}(\bar \nu_n) = v$ \\ 
which is needed only to determine the ordering of the following points: \\ 
and, for each $i \ge 1$, \\ 
$\bar u_{n,i} = \min \big\{ \bar u_n' \le x \le \bar u_{n+1}' : f^{2n+2i}(x) = d \big\}$; \vspace{.05in} \\ 
$\bar c_{2n+2i}^{\,(n,2)} = \min \big\{ \bar u_n' \le x \le \bar u_{n+1}' : f^{2n+2i}(x) = x \big\}$; \vspace{.05in} \\ 
$\bar p_{m+2n+2i}^{\, (n,2)} = \min \big\{ \bar u_{n,1} \le x \le \bar u_{n+1}' : f^{m+2n+2i}(x) = x \big\}$; \vspace{.05in} \\ 
$\bar q_{m+2n+2i}^{\, (n,2)} = \max \big\{ \bar u_n' \le x \le \bar u_{n+1}' : f^{m+2n+2i}(x) = x \big\}$. \vspace{.05in} \\ 
(The relative locations of the periodic points $\bar c_{2n+2i}^{\,(n,2)}$'s with respect to $\bar u_{n,i}$'s are not known.)} 
\end{figure}

For each fixed $n \ge 1$, the collection of all these periodic points $\bar c_{2m+2n+2i}^{\,\, (n,2)}$, $\bar p_{m+2n+2i}^{\,\, (n,2)}$, $\bar q_{m+2n+2i}^{\,\, (n,2)}$, $i \ge 1$, is called a {\it compartment} of the second layer of the {\it basic} tower of periodic points of $f$ associated with $P$.  

\vspace{.2in}

\noindent
{\bf 2.5(a) On the existence of periodic points of $f$ of all even periods $\ge 2n+2$ in $[\bar u_n', \bar u_{n,1}]$ $(\subset [\bar u_n', \bar u_{n+1}']$ $\subset [v, z_0])$ for each $n \ge 1$.}  

In this case, we can apply Lemma 6(1) to get a rough result.  However, since we can obtain a better result, we argue as follows:

For each $n \ge 1$ and $i\ge 1$, we define $\bar \nu_n$, $\bar u_{n,i}$ and $\bar c_{2n+2i}^{\, (n,2)}$ in $[\bar u_n', \bar u_{n+1}']$ 'symmetrically' with respect to $\nu_n$, $u_{n,i}$ and $c_{2n+2i}^{\, (n,2)}$ in $[u_{n+1}, u_n]$ as follows: 

Let $\bar \nu_n$ be a point in $[\bar u_n', \bar u_{n+1}']$ such that $f^{2n}(\bar \nu_n) = v$.  By arguing as those on $[u_{n+1}, u_n]$ and since $f^2(x) < x$ for all $v < x < z_0$, we obtain that, for each $i \ge 1$, the points 
$$
\bar u_{n,i} = \min \big\{ \bar u_n' \le x \le \bar \nu_n : f^{2n+2i}(x) = d \big\} = \min \big\{ \bar u_n' \le x \le \bar u_{n+1}' : f^{2n+2i}(x) = d \big\}\,\,\, \text{and}
$$ 
$$
\bar c_{2n+2i}^{\, (n,2)} = \min \big\{ \bar u_n' \le x \le \bar u_{n,i} : f^{2n+2i}(x) = x \big\} = \min \big\{ \bar u_n' \le x \le \bar u_{n+1}' : f^{2n+2i}(x) = x \big\}\,\,\,\,
$$
exist and it is easy to see that $\bar u_n' < \cdots < \bar u_{n,3} < \bar u_{n,2} < \bar u_{n,1} < \bar \nu_n$ and, by combining with Lemma 4(2), we have 
{\large 
$$
\qquad d < f^{2j}(x) < z_0 \,\,\, \text{for all} \,\,\, \bar u_n' < x < \bar u_{n,i} \,\,\, \text{for all} \,\,\, 1 \le j \le n+i.\qquad\qquad\,\, ({\overline {\dagger\dagger}})
$$}
Consequently, since $f(x) > z \ge z_0$ for all $v < x < z_0$, we obtain that 
$$
\text{all periodic points of} \,\,\, f \,\,\, \text{in} \,\,\, [\bar u_n', \bar u_{n,1}] \,\,\, \text{of odd periods have least periods} \,\,\, \ge 2n+2i+3.
$$
This fact will be used in $\mathsection$3(5).  As for the periodic points $\bar c_{2n+2i}^{\, (n,2)}$'s, we note that $f^{2n+2}\big(\bar c_{2n+2}^{\, (n,2)}\big) = \bar c_{2n+2}^{\, (n,2)} < \bar u_{n,1}$.  Since $f^{2n+2}(\bar u_n') = z_0 \in \big\{ z_0, f(z_0), f(d) \big\}$, it follows from Lemma 3 that \\ $\bar u_n' < \cdots < \bar c_{2n+6}^{\, (n,2)} < \bar c_{2n+4}^{\, (n,2)} < \bar c_{2n+2}^{\, (n,2)} < \bar u_{n,1} < \bar \nu_n < \bar u_{n+1}'$.  

To find the least periods of $\bar c_{2n+2i}^{\, (n,2)}$'s with respect to $f$, our arguments are similar to those of the counterparts $c_{2n+2i}^{\, (n,2)}$'s on $[u_{n+1}, u_n]$.  Let $\ell\big(\bar c_{2n+2i}^{\, (n,2)}\big)$ denote the least period of $\bar c_{2n+2i}^{\, (n,2)}$ with respect to $f$.  Then, by $(\overline{\dagger\dagger})$ above, $\ell\big(\bar c_{2n+2i}^{\, (n,2)}\big)$ is even and, by the following Lemma 4(3), 
$$
f \,\,\, \text{has no periodic points of periods} \, \le 2n+1 \,\,\, \text{in} \,\,\, [\bar u_n', z_0],
$$
we obtain that $\ell\big(\bar c_{2n+2i}^{\, (n,2)}\big)$ is even and $\ge 2n+2$.

Suppose $(2n+2 \le) \,\, \ell\big(\bar c_{2n+2j}^{\, (n,2)}\big) < 2n+2j$ for some $j \ge 1$.  Then since $2n+2j > \ell\big(\bar c_{2n+2j}^{\, (n,2)}\big)$ and $\ell\big(\bar c_{2n+2j}^{\, (n,2)}\big) \ge 2n+2$, it follows from Lemma 3 with $f^{\ell(\bar c_{2n+2j}^{\, (n,2)})}\big(\bar c_{2n+2j}^{\, (n,2)}\big) = \bar c_{2n+2j}^{\, (n,2)}$ and $f^{2n+2}(\bar u_n') = z_0$ that there is a periodic point $\bar c_{2n+2j}^{\, *}$ of $f$ such that $\bar u_n' \le \bar c_{2n+2j}^{\, *} < \bar c_{2n+2j}^{\, (n,2)} < \bar u_{n+1}'$ and $f^{2n+2j}(\bar c_{2n+2j}^{\, *}) = \bar c_{2n+2j}^{\, *}$.  This contradicts the minimality of $\bar c_{2n+2j}^{\, (n,2)}$ in $[\bar u_n', \bar u_{n+1}']$.  Therefore, $\ell\big(\bar c_{2n+2j}^{\, (n,2)}\big) = 2n+2j$ and so, we have shown that, in $[\bar u_n', \bar u_{n,1}] \, (\subset [\bar u_n', \bar u_{n+1}'] \subset [v, z_0])$,
$$
\text{for each} \,\,\, n \ge 1 \,\,\, \text{and} \,\,\, i \ge 1, \, \bar c_{2n+2i}^{\, (n,2)} \,\,\, \text{is a period-}(2n+2i) \,\,\, \text{point of} \,\,\, f.
$$

Note that, for each $n \ge 1$, we have $\bar c_{2n+2}^{\, (n,2)} = \bar c_{2n+2}$ \big(these points $\bar c_{2n+2}$'s are defined in $\mathsection 1$(5)\big).  That is, the second layer periodic points contain some (but not all) periodic points on the first layer.

Unfortunately, here we can not determine the relative locations of the periodic points $\bar c_{2n+2i}^{\, (n,2)}$'s with respect to the points $\bar u_{n,i}$'s which are needed to continue the {\it basic} tower-building process to the third layer.  

\noindent
{\bf 2.5(b) On the existence of periodic points of $f$ of all odd periods $\ge m+2n+2$ in $[\bar u_{n,1}, \bar \nu_n]$ $(\subset [\bar u_n', \bar u_{n+1}']$ $\subset [v, z_0])$ for each $n \ge 1$.}

Since $f^{2n+2}(\bar u_{n,1}) = d$, we have 
$$
f^{2n+3}(\bar u_{n,1}) = f(d).
$$
On the other hand, recall that 
$$
f^{m+2(n+1)}(\bar \nu_n) = f^{m+2}\big(f^{2n}(\bar \nu_n)\big) = f^m\big(f^2(v)\big) = \min P.
$$
Furthermore, by Lemma 4(3), we have 
$$
\text{all periodic points of} \,\, f \,\, \text{in} \, [\bar u_n', \bar u_{n+1}'] \, \text{with {\it odd} periods have least periods} \, \ge 2n+3.\,
$$
Therefore, we can apply Lemma 5(1) with
{\large 
$$
f^{2n+3}(\bar u_{n,1}) = f(d) \in \big\{ z_0, f(z_0), f(d) \big\} \,\,\, \text{and} \,\,\, f^{m+2(n+1)}(\bar \nu_n) = \min P
$$}
to obtain that, for each $i \ge 0$, the point 
\begin{multline*}
$$
\qquad\qquad\qquad \bar p_{m+2(n+1)+2i}^{\, (n,2)} = \min \big\{ \bar u_{n,1} \le x \le \bar \nu_n : f^{m+2(n+1)+2i}(x) = x \big\} \\ 
= \min \big\{ \bar u_{n,1} \le x \le \bar u_{n+1}' : f^{m+2(n+1)+2i}(x) = x \big\},\qquad\qquad\qquad\,\,\,\,
$$
\end{multline*}
exists and $\,\, \bar u_n' \,\, < \,\, \bar u_{n,1} \,\, < \,\, \cdots \,\, < \,\, \bar p_{m+2n+6}^{\, (n,2)} \,\, < \,\, \bar p_{m+2n+4}^{\, (n,2)} \,\, < \,\, \bar p_{m+2n+2}^{\, (n,2)} \,\, < \,\, \bar \nu_n \,\, < \,\, \bar u_{n+1}'$.  Furthermore, for each $i \ge 0$, \, the point $\bar p_{m+2n+2+2i}^{\, (n,2)}$ is a period-$(m+2n+2+2i)$ point of $f$, or, equivalently, in $[\bar u_{n,1}, \bar \nu_n]$ $(\subset [\bar u_n', \bar u_{n+1}']$ $\subset [v, z_0])$, 
$$
\text{for each} \,\,\, n \ge 1 \,\,\, \text{and} \,\,\, i \ge 1, \, \text{the point} \,\,\, \bar p_{m+2n+2i}^{\, (n,2)} \,\,\, \text{is a period-}(m+2n+2i) \,\,\, \text{point of} \,\,\, f.
$$

\noindent
{\bf 2.5(c) On the existence of periodic points of $f$ of all odd periods $\ge m+2n+2$ in $[\bar \nu_n, \bar u_{n+1}']$ $(\subset [\bar u_n', \bar u_{n+1}']$ $\subset [v, z_0])$ for each $n \ge 1$.}

Since $f^{2n+2}(\bar u_{n+1}') = d$, we have 
$$
f^{2n+3}(\bar u_{n+1}') = f(d).
$$
On the other hand, We have 
$$
f^{m+2(n+1)}(\bar \nu_n) = f^{m+2}\big(f^{2n}(\bar \nu_n)\big) = f^m\big(f^2(v)\big) = \min P.
$$
Furthermore, by Lemma 4(3), we have 
$$
\text{all periodic points of} \,\, f \,\, \text{in} \, [\bar u_n', \bar u_{n+1}'] \, \text{with {\it odd} periods have least periods} \, \ge 2n+3.\,
$$
Therefore, we can apply Lemma 5(2) with 
{\large 
$$
f^{m+2n+3}(\bar u_{n+1}') = f(d) \in \big\{ z_0, f(z_0), f(d) \big\} \,\,\, \text{and} \,\,\, f^{m+2(n+1)}(\bar \nu_n) = \min P
$$}
to obtain that, for each $i \ge 0$, the point 
\begin{multline*}
$$
\qquad\qquad\qquad \bar q_{m+2(n+1)+2i}^{\, (n,2)} = \max \big\{ \bar \nu_n \le x \le \bar u_{n+1}' : f^{m+2(n+1)+2i}(x) = x \big\} \\ 
= \max \big\{ \bar u_n' \le x \le \bar u_{n+1}' : f^{m+2(n+1)+2i}(x) = x \big\},\qquad\qquad\qquad\,\,\,\,\,\,\,
$$
\end{multline*}
exists and $\, \bar u_n' \, < \, \bar \nu_n \, < \, \bar q_{m+2n+2}^{\, (n,2)} \, < \, \bar q_{m+2n+4}^{\, (n,2)} \, < \, \bar q_{m+2n+6}^{\, (n,2)} \, < \, \bar q_{m+2n+8}^{\, (n,2)} \, < \, \cdots \, < \, \bar u_{n+1}' \, < \, z_0$.  Furthermore, for each $i \ge 0$, the point $\bar q_{m+2n+2+2i}^{\, (n,2)}$ is a period-$(m+2n+2+2i)$ point of $f$, or, equivalently, in $[\bar \nu_n, \bar u_{n+1}']$ $(\subset [\bar u_n', \bar u_{n+1}']$ $\subset [v, z_0])$, 
$$
\text{for each} \,\,\, n \ge 1 \,\,\, \text{and} \,\,\, i \ge 1, \, \text{the point} \,\,\, \bar q_{m+2n+2i}^{\, (n,2)} \,\,\, \text{is a period-}(m+2n+2i) \,\,\, \text{point of} \,\,\, f.
$$ 
$$\aleph \qquad\qquad\qquad \aleph \qquad\qquad\qquad \aleph \qquad\qquad\qquad \aleph \qquad\qquad\qquad \aleph$$
\indent In summary, on the interval $[d, u_n]$, we have $f^{2n}(d) = z_0, f^{2n}(u_n) = d$ and $f^{2n}(\nu_n) = v$.  So, $f^{2n+2}(d) = z_0 = f^{2n+2}(u_n)$ and $f^{2n+2}(\nu_n) = \min P$.  Thus, 
$$
\text{the graph of} \, y = f^{2n+2}(x) \, \text{looks 'roughly' like a {\it skew} 'V' shape (at the three points} \,\, d, \nu_n, u_n)
$$
and the point $\nu_n$ happens to satisfy that $(d <) \, u_{n+1} < \nu_n < u_n$. 

Similarly, on the interval $[\bar u_n', z_0]$, we have $f^{2n}(\bar u_n') = d, f^{2n}(z_0) = z_0$ and $f^{2n}(\bar \nu_n) = v$.  So, $f^{2n+2}(\bar u_n') = z_0 = f^{2n+2}(z_0)$ and $f^{2n+2}(\bar \nu_n) = \min P$.  Thus, $$
\text{the graph of} \, y = f^{2n+2}(x) \, \text{looks 'roughly' like a {\it skew} 'V' shape (at the three points} \,\, \bar u_n', \bar \nu_n, z_0)
$$
and the point $\bar \nu_n$ happens to satisfy that $\bar u_n' < \bar \nu_n < \bar u_{n+1}' \, (< z_0)$. 

We have shown that there exist 4 monotonic sequences $< u_{n,i}' >$, $< c_{2n+2i}'^{\,\,(n,2)} >$, $< q_{m+2n+2i}^{\,\,(n,2)} >$ and $< p_{m+2n+2i}^{\,\,(n,2)} >$ of points in $[u_{n+1}, u_n]$, where  
\[
\begin{array}{l}
\textrm{$< u_{n,i}' >$ is a monotonic sequence of $d$-points which will be used to continue the}\\
\textrm{{\it basic} tower-building process to the third layer and, $< c_{2n+2i}'^{\,\,(n,2)} >$, $< q_{m+2n+2i}^{\,\,(n,2)} >$ and}\\
\textrm{$< p_{m+2n+2i}^{\,\,(n,2)} >$ are 3 monotonic sequences of periodic points of $f$ whose union}\\
\textrm{constitutes a compartment of the second layer of the {\it basic} tower of periodic points of $f$}
\end{array}  
\]
and 4 monotonic sequences $< \bar u_{n,i} >$, $< \bar c_{2n+2i}^{\,\,(n,2)} >$, $< \bar p_{m+2n+2i}^{\,\,(n,2)} >$ and $< \bar q_{m+2n+2i}^{\,\,(n,2)} >$ of points in $[\bar u_n', \bar u_{n+1}']$, where
\[
\begin{array}{l}
\textrm{$< \bar u_{n,i} >$ is a monotonic sequence of $d$-points which will be used to continue the}\\
\textrm{{\it basic} tower-building process to the third layer and, $< \bar c_{2n+2i}^{\,\,(n,2)} >$, $< \bar p_{m+2n+2i}^{\,\,(n,2)} >$ and}\\
\textrm{$< \bar q_{m+2n+2i}^{\,\,(n,2)} >$ are 3 monotonic sequences of periodic points of $f$ whose union}\\
\textrm{constitutes a compartment of the second layer of the {\it basic} tower of periodic points of $f$}
\end{array}  
\]
such that, for each $n \ge 1$,
{\large
\begin{multline*}
$$
u_{n+1} < \cdots < p_{m+2n+4}^{\,\,(n,2)} < p_{m+2n+2}^{\,\,(n,2)} < \nu_n < q_{m+2n+2}^{\,\,(n,2)} < q_{m+2n+4}^{\,\,(n,2)} < \cdots \\ 
< u_{n,1}' < c_{2n+2}'^{\,\,(n,2)} < u_{n,2}' < c_{2n+4}'^{\,\,(n,2)} < \cdots < u_n,$$\end{multline*}}{\large \begin{multline*}$$\bar u_n' < \cdots < \bar u_{n,3} < \bar u_{n,2} < \bar u_{n,1} \,\, {\it and} \\ 
\bar u_n' < \cdots < \bar c_{2n+6}^{\,\,(n,2)} < \bar c_{2n+4}^{\,\,(n,2)} < \bar c_{2n+2}^{\,\,(n,2)} < \bar u_{n,1} < \qquad\qquad\qquad\qquad\quad \\ 
\cdots < \bar p_{m+2n+4}^{\,\,(n,2)} < \bar p_{m+2n+2}^{\,\,(n,2)} < \bar \nu_n < \bar q_{m+2n+2}^{\,\,(n,2)} < \bar q_{m+2n+4}^{\,\,(n,2)} < \cdots < \bar u_{n+1}',
$$
\end{multline*}}
and, for each $i \ge 1$, 
$$
\, u_{n,i}' = \max \big\{ u_{n+1} \le x \le u_n : f^{2n+2i}(x) = d \big\},
$$ 
$$
\bar u_{n,i} = \min \big\{ \bar u_n' \le x \le \bar u_{n+1}' : f^{2n+2i}(x) = d \big\},
$$ 
$$
d < f^{2j}(x) < z_0 \,\,\, \text{for all} \,\,\, 1 \le j \le n+i \,\,\, \text{and} \,\,\, \text{all} \,\,\, x \in (u_{n,i}', u_n) \cup (\bar u_n', \bar u_{n,i}) \,\,\, \text{and}
$$ 
$$
{\large \text{except possibly}} \,\,\, c_{4n}, \,\,\, \text{all} \,\,\, c_{2n+2i} \,\,\, \text{and} \,\,\, \bar c_{2n+2i}, \,\, i \ne n, \,\, \text{are period-}(2n+2i) \,\,\, \text{points of} \,\,\, f \,\,\, \text{and} 
$$ 
$$
p_{m+2n+2i}^{\, (n,2)}, \,\,\, q_{m+2n+2i}^{\, (n,2)}, \,\,\, \bar p_{m+2n+2i}^{\, (n,2)} \,\,\, \text{and} \,\,\, \bar q_{m+2n+2i}^{\, (n,2)} \,\,\, \text{are all period-}(m+2n+2i) \,\,\, \text{points of} \,\,\, f.
$$

\noindent
Furthermore, on the interval $[\min P, d]$, we have found 4 monotonic sequences \\ $< \tilde \mu_{m,n,i} >$, $< \tilde p_{m+2n+2i}^{\,\, (n,2)} >$, $< \tilde c_{2m+2n+2i}^{\, (n,2)} >$ and $< \tilde c_{2m+2n+2i}'^{\, (n,2)} >$ of points, where
\[
\begin{array}{l}
\textrm{$< \tilde \mu_{m,n,i} >$ is a monotonic sequence of $d$-points which will be used to continue the}\\
\textrm{{\it basic} tower-building process to the third layer and, $< \tilde p_{m+2n+2i}^{\,\, (n,2)} >$, $< \tilde c_{2m+2n+2i}^{\, (n,2)} >$ and}\\
\textrm{$< \tilde c_{2m+2n+2i}'^{\, (n,2)} >$ are 3 monotonic sequences of periodic points of $f$ whose union}\\
\textrm{constitutes a compartment of the second layer of the {\it basic} tower of periodic points of $f$}
\end{array}  
\]
such that, for each $n \ge 0$ (the point $\tilde p_{m+2n+2}^{\,\, (n,2)}$ is excluded),
\begin{multline*}
$$
\tilde \mu_{m,n}' < \cdots < \tilde \mu_{m,n,3} < \tilde p_{m+2n+6}^{\,\, (n,2)} < \tilde \mu_{m,n,2} < \tilde p_{m+2n+4}^{\,\, (n,2)} < \tilde \mu_{m,n,1} \\ 
\cdots < \tilde c_{2m+2n+4}^{\, (n,2)} < \tilde c_{2m+2n+2}^{\, (n,2)} < \tilde \nu_{m,n} < \tilde c_{2m+2n+2}'^{\, (n,2)} < \tilde c_{2m+2n+4}'^{\, (n,2)} < \cdots < \tilde \mu_{m,n+1}'.
$$
\end{multline*}  

\noindent
On the interval $[\breve u_0', v] \, (\subset [u_1, v])$, we have found 4 monotonic sequences \\ 
$< \breve \mu_{m,n,i} >$, $< \breve q_{m+2n+2i}^{\,\, (n,2)} >$, $< \breve c_{2m+2n+2i}'^{\, (n,2)} >$ and $< \breve c_{2m+2n+2i}^{\, (n,2)} >$ of points, where 
\[
\begin{array}{l}
\textrm{$< \breve \mu_{m,n,i} >$ is a monotonic sequence of $d$-points which will be used to continue the}\\
\textrm{{\it basic} tower-building process to the third layer and, $< \breve q_{m+2n+2i}^{\,\, (n,2)} >$, $< \breve c_{2m+2n+2i}'^{\, (n,2)} >$ and}\\
\textrm{$< \breve c_{2m+2n+2i}^{\, (n,2)} >$ are 3 monotonic sequences of periodic points of $f$ whose union}\\
\textrm{constitutes a compartment of the second layer of the {\it basic} tower of periodic points of $f$}
\end{array}  
\]
such that, for each $n \ge 1$, 
\begin{multline*}
$$
\breve \mu_{m,n+1} < \cdots < \breve c_{2m+2n+4}^{\, (n,2)} < \breve c_{2m+2n+2}^{\, (n,2)} < \breve \nu_{m,n} < \breve c_{2m+2n+2}'^{\, (n,2)} < \breve c_{2m+2n+4}'^{\, (n,2)} < \cdots \\  
\breve \mu_{m,n,1}' < \breve q_{m+2n+2}^{\,\, (n,2)} < \breve \mu_{m,n,2}' < \breve q_{m+2n+4}^{\,\, (n,2)} < \breve \mu_{m,n,3}' < \breve \mu_{m,n}.
$$
\end{multline*}

\noindent
On the interval $[v, \hat u_0] \, (\subset [v, \bar u_1'])$, we have found 4 monotonic sequences \\ 
$< \hat \mu_{m,n,i} >$, $< \hat p_{m+2n+2i}^{\,\, (n,2)} >$, $< \hat c_{2m+2n+2i}^{\, (n,2)} >$ and $< \hat c_{2m+2n+2i}'^{\, (n,2)} >$ of points, where 
\[
\begin{array}{l}
\textrm{$< \hat \mu_{m,n,i} >$ is a monotonic sequence of $d$-points which will be used to continue the}\\
\textrm{{\it basic} tower-building process to the third layer and, $< \hat p_{m+2n+2i}^{\,\, (n,2)} >$, $< \hat c_{2m+2n+2i}^{\, (n,2)} >$ and}\\
\textrm{$< \hat c_{2m+2n+2i}'^{\, (n,2)} >$ are 3 monotonic sequences of periodic points of $f$ whose union}\\
\textrm{constitutes a compartment of the second layer of the {\it basic} tower of periodic points of $f$}
\end{array}  
\]
such that, for each $n \ge 1$, 
\begin{multline*}
$$
\hat \mu_{m,n}' < \cdots < \hat \mu_{m,n,3} < \hat p_{m+2n+4}^{\,\, (n,2)} < \hat \mu_{m,n,2} < \hat p_{m+2n+2}^{\,\, (n,2)} < \hat \mu_{m,n,1} \\ 
\cdots < \hat c_{2m+2n+4}^{\, (n,2)} < \hat c_{2m+2n+2}^{\, (n,2)} < \hat \nu_{m,n} < \hat c_{2m+2n+2}'^{\, (n,2)} < \hat c_{2m+2n+4}'^{\, (n,2)} < \cdots < \hat \mu_{m,n+1}'.
$$
\end{multline*}

{\large 
We call the collection of all these periodic points, arranged from $\min P$ (exclusive) to $z_0$ (exclusive), {\it excluding} the periodic points $\tilde p_{m+2n+2}^{\, (n,2)}$, $n \ge 1$ so that the convex hulls of all compartments on the second layer are pairwise disjoint,
$$\quad \tilde p_{m+2n+2i+2}^{\, (n,2)}, \,\, \tilde c_{2m+2n+2i}^{\,(n,2)}, \,\, \tilde c_{2m+2n+2i}'^{\,(n,2)}, \,\, i \ge 1, \, n \ge 1,\,$$
$$p_{m+2n+2i}^{\,\,(n,2)}, \,\, q_{m+2n+2i}^{\,\,(n,2)}, \,\, c_{2m+2n+2i}'^{\,\,(n,2)}, \,\, i \ge 1, \, n \ge 1, \,\,$$
$$\breve c_{2m+2n+2i}^{\, (n,2)}, \,\, \breve c_{2m+2n+2i}'^{\,(n,2)}, \,\, \breve q_{m+2n+2i}^{\,(n,2)}, \,\, i \ge 1, \, n \ge 1,\,$$
$$\hat p_{m+2n+2i}^{\,(n,2)}, \,\, \hat c_{2m+2n+2i}^{\,(n,2)}, \,\, \hat c_{2m+2n+2i}'^{\, (n,2)}, \,\, i \ge 1, \, n \ge 1,$$
$$\bar c_{2m+2n+2i}^{\, (n,2)}, \,\, \bar p_{m+2n+2i}^{\,(n,2)}, \,\, \bar q_{m+2n+2i}^{\,(n,2)}, \,\, i \ge 1, \, n \ge 1,\,\,\,$$
the {\it second layer} of the {\it basic} tower of periodic points of $f$ associated with $P$ (see Figures 2-6).}

\vspace{.2in}

\noindent
{\bf $\mathsection$3. The third layer of the {\it basic} tower of periodic points of $f$ associated with $P$.}

By the way we find the second layer of the {\it basic} tower of periodic points of $f$ associated with $P$, we can write down the formulas for the periodic points of the third layer as follows: 

\noindent
{\bf 3.1 A compartment of the third layer of the {\it basic} tower in $[\tilde \mu_{m,n,k+1}, \tilde \mu_{m,n,k}] \, (\subset [\tilde \mu_{m,n}', \tilde \mu_{m,n+1}'] \subset [\min P, d])$.}

For each $n \ge 0$ and $k \ge 1$, on the interval $[\tilde \mu_{m,n,k+1}, \tilde \mu_{m,n,k}] \, (\subset [\tilde \mu_{m,n}', \tilde \mu_{m,n+1}'] \subset [\min P, d])$, let $\tilde \nu_{m,n,k}$ be a point in $[\tilde \mu_{m,n}', \tilde \mu_{m,n,k}]$ such that $f^{m+2n+2k}(\tilde \nu_{m,n,k}) = v$.  It turns out that $\tilde \mu_{m,n,k+1} < \tilde \nu_{m,n,k} < \tilde \mu_{m,n,k}$.  For each $i \ge 1$, let, with $x$-coordinates moving from point $\tilde \mu_{m,n,k}$ via point $\tilde \nu_{m,n,k}$ to point $\tilde \mu_{m,n,k+1}$, \big(the number 3 in the superscripts indicates the third layer\big), 
{\large 
$$\,\,\,\, \tilde \mu_{m,n,k,i}' = \max \big\{ \tilde \mu_{m,n,k+1} \le x \le \tilde \mu_{m,n,k} : f^{m+2n+2k+2i}(x) = d \big\},\qquad\qquad\qquad$$ 
$$\,\,\,\,\,\,\, \tilde q_{m+2n+2k+2i}^{\, (n,k,3)} = \max \big\{ \tilde \mu_{m,n,k+1} \le x \le \tilde \mu_{m,n,k} : f^{m+2n+2k+2i}(x) = x \big\},\qquad\qquad\quad$$ 
$$\tilde c_{2m+2n+2k+2i}'^{\, (n,k,3)} = \max \big\{ \tilde \mu_{m,n,k+1} \le x \le \tilde \mu_{m,n,k,1}' : f^{2m+2n+2k+2i}(x) = x \big\} \,\,\,\, \text{and}\,$$ 
$$\tilde c_{2m+2n+2k+2i}^{\, (n,k,3)} = \min \big\{ \tilde \mu_{m,n,k+1} \le x \le \tilde \mu_{m,n,k} : f^{2m+2n+2k+2i}(x) = x \big\}.\qquad\,\,\,\,\,\,
$$}

\begin{figure}[htbp] 
\centering
\includegraphics[width=6in,height=2in]{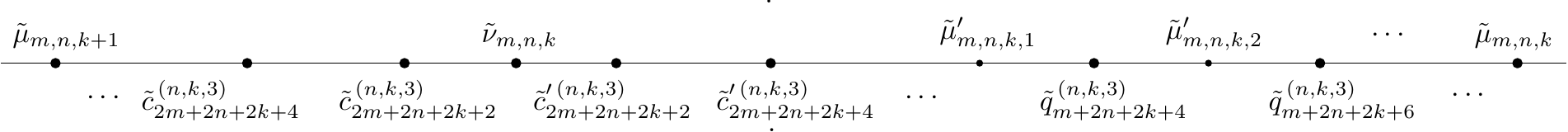} 

{Figure 7: A compartment of the third layer of the {\it basic} tower of periodic points of $f$ associated with $P$ in the interval $[\tilde \mu_{m,n,k+1}, \tilde \mu_{m,n,k}]$ $(\subset [\tilde \mu_{m,n}, \tilde \mu_{m,n+1}]$ $\subset [\min P, d])$, \\
where $\tilde \nu_{m,n,k}$ is an {\it auxiliary} point in $[\tilde \mu_{m,n,k+1}, \tilde \mu_{m,n,k}]$ such that $f^{m+2n+2k}(\tilde \nu_{m,n,k}) = v$ \\ 
which is needed only to determine the ordering of the following points: \\ 
and, for each $i \ge 1$, \\ 
$\tilde \mu_{m,n,k,i}' = \max \big\{ \tilde \mu_{m,n,k+1} \le x \le \tilde \mu_{m,n,k} : f^{m+2n+2k+2i}(x) = d \big\}$; \vspace{.05in} \\ 
$\tilde q_{m+2n+2k+2i}^{\, (n,k,3)} = \max \big\{ \tilde \mu_{m,n,k+1} \le x \le \tilde \mu_{m,n,k} : f^{m+2n+2k+2i}(x) = x \big\}$; \vspace{.05in} \\ 
$\tilde c_{2m+2n+2k+2i}'^{\, (n,k,3)} = \max \big\{ \tilde \mu_{m,n,k+1} \le x \le \tilde \mu_{m,n,k,1}' : f^{2m+2n+2k+2i}(x) = x \big\}$; \vspace{.05in} \\ 
$\tilde c_{2m+2n+2k+2i}^{\, (n,k,3)} = \min \big\{ \tilde \mu_{m,n,k+1} \le x \le \tilde \mu_{m,n,k} : f^{2m+2n+2k+2i}(x) = x \big\}$. \vspace{.05in}} 
\end{figure}

\noindent
For each fixed $n \ge 0$ and $k \ge 1$, the collection of all these periodic points $\tilde q_{m+2n+2k+2i+2}^{\, (n,k,3)}$, $\tilde c_{2m+2n+2k+2i}'^{\, (n,k,3)}$, $\tilde c_{2m+2n+2k+2i}^{\, (n,k,3)}, i \ge 1$ \big(note that the point $\tilde q_{m+2n+2k+2}^{\,\,(n,k,3)}$ is excluded\big), is called a {\it compartment} of the third layer of the {\it basic} tower of periodic points of $f$ associated with $P$.  

It is easy to see that we have 
\begin{multline*}
$$
\tilde \mu_{m,n,k+1} < \cdots < \tilde c_{2m+2n+2k+4}^{\, (n,k,3)} < \tilde c_{2m+2n+2k+2}^{\, (n,k,3)} < \tilde \nu_{m,n,k} < \tilde c_{2m+2n+2k+2}'^{\, (n,k,3)} < \tilde c_{2m+2n+2k+4}'^{\, (n,k,3)} < \cdots \\ < \tilde \mu_{m,n,k,1}' < \tilde q_{m+2n+2k+4}^{\, (n,k,3)} < \tilde \mu_{m,n,k,2}' < \tilde q_{m+2n+2k+6}^{\, (n,k,3)} < \tilde \mu_{m,n,k,3}' < \cdots < \tilde \mu_{m,n,k}
$$
\end{multline*}
with the point $\tilde q_{m+2n+2k+2}^{\, (n,k,3)}$ {\it excluded} so that the convex hulls of the two sets: $\big\{ \tilde c_{2m+2n+2k+2i}'^{\, (n,k,3)}: i \ge 1 \big\}$ and $\big\{ \tilde q_{m+2n+2k+2j}^{\, (n,k,3)}: j \ge 2 \big\}$ are disjoint.

The following result can be viewed as a continuation of Lemma 7 from the interval $[\tilde \mu_{m,n}', d]$ to the subinterval $[\tilde \mu_{m,n}', \tilde \mu_{m,n,k}]$ \big(and later to the subinterval $[\tilde \mu_{m,n,k,i}', \tilde \mu_{m,n,k}]$ where all periodic points of $f$ of odd periods have least periods $\ge m+2n+2k+2i+2$ and then continue on to the higher layers of the {\it basic} tower and so on\big):

\noindent
{\bf Lemma 11.}
{\it For each $n \ge 1$ and $k \ge 1$, all periodic points of $f$ in $[\tilde \mu_{m,n}', \tilde \mu_{m,n,k}]$ of odd periods have least periods $\ge m+2n+2k+2$.}

\noindent
{\it Proof.}
Recall that, by Lemma 7, all periodic points of $f$ in $[\tilde \mu_{m,n}', \tilde \mu_{m,n,k}]$ $(\subset [\tilde \mu_{m,n}', \tilde \mu_{m,n+1}'])$ of odd periods have least periods $\ge m+2n+2$.  Suppose $f$ had a periodic point $\tilde p$ of odd period $\ell(\tilde p) \le m+2n+2k$ in $[\tilde \mu_{m,n}', \tilde \mu_{m,n,k}]$.  Then, since $\ell(\tilde p)$ is odd and $\ell(\tilde p) \ge m+2n+2$, we have $f^{\ell(\tilde p)}(\tilde \mu_{m,n}') = z_0$.  Since $m+2n+2k \ge \ell(\tilde p)$, it follows from Lemma 3 that there is a point $\tilde w_{m+2n+2k}$ in $[\tilde \mu_{m,n}', \tilde \mu_{m,n,k}]$ $(\subset [\min P, d])$ such that $f^{m+2n+2k}(\tilde w_{m+2n+2k}) = \tilde w_{m+2n+2k}$.  Since $f^{m+2n+2k}(\tilde \mu_{m,n}') - d = z_0 - d > 0$ and $f^{m+2n+2k}(\tilde w_{m+2n+2k}) - d = \tilde w_{m+2n+2k} - d < 0$, there is a point $\tilde \mu_{m,n,k}^*$ in $[\tilde \mu_{m,n}', \tilde w_{m+2n+2k}]$ such that $f^{m+2n+2k}(\tilde \mu_{m,n,k}^*) = d$.  Since $\tilde \mu_{m,n}' < \tilde \mu_{m,n,k}^* < \tilde w_{m+2n+2k} < \tilde \mu_{m,n,k} < \tilde \mu_{m,n+1}'$, this contradicts the minimality of $\tilde \mu_{m,n,k}$ in $[\tilde \mu_{m,n}', \tilde \mu_{m,n+1}']$.
$\hfill\square$

Recall that $\tilde \nu_{m,n,k}$ is a point in $[\tilde \mu_{m,n}', \tilde \mu_{m,n,k}]$ such that $f^{m+2n+2k}(\tilde \nu_{m,n,k}) = v$.  We can apply Lemma 5(1) and Lemma 11 with 
$$
f^{m+2n+2k+2}(\tilde \mu_{m,n,k}) = z_0 \in \big\{ z_0, f(z_0), f(d) \big\} \,\,\, \text{and} \,\,\, f^{m+2(n+k+1)}(\tilde \nu_{m,n,k}) = \min P
$$
to obtain that, 
\begin{multline*}
$$
\text{for each} \,\,\, i \ge 1, \, \text{the point} \,\,\, \tilde q_{m+2n+2k+2i}^{\, (n,k,3)} = \max \big\{ \tilde \mu_{m,n,k+1} \le x \le \tilde \mu_{m,n,k} : f^{m+2n+2k+2i}(x) = x \big\}  \\ 
\text{exists and is a period-}(m+2n+2k+2i) \,\,\, \text{point of} \,\,\, f.
$$
\end{multline*}

As for the periods of $\tilde c_{2m+2n+2k+2i}'^{\, (n,k,3)}$'s and $\tilde c_{2m+2n+2k+2i}^{\, (n,k,3)}$'s, we apply Lemma 6 with 
\begin{multline*}
$$
\qquad\quad\quad\,\, f^{m+2n+2k+3}(\tilde \mu_{m,n,k,1}') = f(d) \,\,\, \text{and} \,\,\, f^{2(m+n+k+1)}(\tilde \nu_{m,n,k}) = \min P \,\,\, \text{and} \\ \,\,\, f^{m+2n+2k+3}(\tilde \mu_{m,n,k+1}) = f(d) \,\,\, \text{and} \,\,\, f^{2(m+n+k+1)}(\tilde \nu_{m,n,k}) = \min P \,\,\, \text{respectively}\quad\,\,\,\,
$$
\end{multline*}
to obtain that, for each $(m+n+k+1)+i \ge \max \big\{ m+2n+2k+3, m+2 \big\} = m+2n+2k+3$, i.e., for each $i \ge n+k+2$, both 
\begin{multline*}
$$
\qquad\qquad \tilde c_{2m+2n+2k+2+2i}'^{\, (n,k,3)} = \max \big\{ \tilde \mu_{m,n,k+1} \le x \le \tilde \mu_{m,n,k,1}' : f^{2m+2n+2k+2+2i}(x) = x \big\} \,\,\, \text{and} \,\, \\ 
\,\,\,\, \tilde c_{2m+2n+2k+2+2i}^{\, (n,k,3)} = \min \big\{ \tilde \mu_{m,n,k+1} \le x \le \tilde \mu_{m,n,k} : f^{2m+2n+2k+2+2i}(x) = x \big\}\qquad\quad\,\,\,\,\,\,\,\,\,
$$
\end{multline*}
exist and are period-$(2m+2n+2k+2+2i)$ points of $f$, or, equivalently, 
\begin{multline*}
$$
\text{for each} \,\, i \ge n+k+3, \, \text{both} \,\,\, \tilde c_{2m+2n+2k+2i}'^{\, (n,k,3)} \,\, \text{and} \,\, \tilde c_{2m+2n+2k+2i}^{\, (n,k,3)} \\ \text{are period-}(2m+2n+2k+2i) \,\, \text{points of} \,\, f.
$$
\end{multline*}

\vspace{.1in}

\noindent
{\bf 3.2 A compartment of the third layer of the {\it basic} tower in $[u_{n,k}', u_{n,k+1}'] \, (\subset (u_{n+1}, u_n] \subset [d, u_1])$.}

For each $n \ge 1$ and $k \ge 1$, we consider the interval $[u_{n,k}', u_{n,k+1}'] \, (\subset (u_{n+1}, u_n] \subset [d, u_1])$.  Let $\nu_{n,k}$ be a point in $[u_{n,k}', u_n]$ such that $f^{2n+2k}(\nu_{n,k}) = v$.  It turns out that $u_{n,k}' < \nu_{n,k} < u_{n,k+1}'$.  For each $i \ge 1$, let, with $x$-coordinates moving from point $u_{n,k}'$ via point $\nu_{n,k}$ to point $u_{n,k+1}'$, \big(the number 3 in the superscripts indicates the third layer\big), 
{\large 
$$u_{n,k,i} = \min \{ u_{n,k}' \le x \le u_{n,k+1}' : f^{2n+2k+2i}(x) = d \},\qquad\qquad\quad\,\,\,\,\,$$ 
$$\quad c_{2n+2k+2i}^{\, (n,k,3)} = \min \{ u_{n,k}' \le x \le u_{n,k+1}' : f^{2n+2k+2i}(x) = x \},\qquad\quad\qquad\,$$ 
$$\qquad\quad\, p_{m+2n+2k+2i}^{\, (n,k,3)} = \min \{ u_{n,k,1} \le x \le u_{n,k+1}' : f^{m+2n+2k+2i}(x) = x \} \,\, \text{and}\qquad\qquad$$ 
$$\qquad\quad\, q_{m+2n+2k+2i}^{\, (n,k,3)} = \max \{ u_{n,k}' \le x \le u_{n,k+1}' : f^{m+2n+2k+2i}(x) = x \}.\qquad\qquad\,\,\,\,\,\,\,\,\,
$$}
\indent For each fixed $n \ge 1$ and $k \ge 1$, the collection of all these periodic points $c_{2n+2k+2i}^{\, (n,k,3)}$, $p_{m+2n+2k+2i}^{\, (n,k,3)}$, $q_{m+2n+2k+2i}^{\, (n,k,3)}$, $i \ge 1$ is called a {\it compartment} of the third layer of the {\it basic} tower of periodic points of $f$ associated with $P$.    

It is easy to see that we have 
\begin{multline*}
$$
u_{n,k}' < \cdots < u_{n,k,3} < c_{2n+2k+4}^{\, (n,k,3)} < u_{n,k,2} < c_{2n+2k+2}^{\, (n,k,3)} < u_{n,k,1} < \cdots < \\ 
p_{m+2n+2k+4}^{\, (n,k,3)} < p_{m+2n+2k+2}^{\, (n,k,3)} < \nu_{n,k} < q_{m+2n+2k+2}^{\, (n,k,3)} < q_{m+2n+2k+4}^{\, (n,k,3)} < \cdots < u_{n,k+1}' \,\,\, \text{and}
$$
\end{multline*}
{\large 
$$
\quad d < f^{2j}(x) < z_0 \,\,\, \text{for all} \,\,\, 1 \le j \le n+k+i \,\,\, \text{and all} \,\,\, u_{n,k}' < x < u_{n,k,i}, \quad\,\,\,\,\, (\dagger\dagger)
$$}
and, since $f(x) > z \ge z_0$ for all $d < x < z_0$, we also have (by taking $i = 1)$, 
$$
\text{all periodic points of} \,\,\, f \,\,\, \text{in} \,\,\, [u_{n,k}', u_{n,k,1}] \,\,\, \text{of odd periods have least periods} \,\, \ge 2n+2k+3.
$$
Since $f^{2n+2k+2}\big([u_{n,k}', u_{n,k,1}]\big) \supset [d, z_0] \supset \{ v \}$, we let $\nu_{n,k+1}$ be a point in $[u_{n,k}', u_{n,k,1}]$ such that $f^{2n+2k+2}(\nu_{n,k+1}) = v$.  Then we can apply Lemma 5(1) with 
{\large 
$$
f^{2n+2k+3}(u_{n,k}') = f(z_0) \,\,\, \text{and} \,\,\, f^{m+2n+2k+4}(\nu_{n,k+1}) = \min P
$$}
to obtain that, for each $i \ge 2$ (not $i \ge 0$), the point 
$$
p_{m+2n+2k+2i}^{\,*} = \min \big\{ u_{n,k}' \le x \le \nu_{n,k+1}: f^{m+2n+2k+2i}(x) = x \big\}
$$
exists and is a period-$(m+2n+2k+2i)$ point of $f$.  However, these periodic points $p_{m+2n+2k+2i}^{\,*}$'s are {\it interspersed} with the periodic points $c_{2n+2k+2i}^{\, (n,k,3)}$'s of $f$ of even periods.  To make things simple, we do not count them in the third layer of the {\it basic} tower of periodic points of $f$ associated with $P$.  

As for the periods of $p_{m+2n+2k+2i}^{\, (n,k,3)}$'s and $q_{m+2n+2k+2i}^{\, (n,k,3)}$'s, we apply Lemma 5 with 
$$
f^{2n+2k+3}\big(u_{n,k,1}\big) = f(d) \,\,\, \text{and} \,\,\, f^{m+2(n+k+1)}(\nu_{n,k}) = \min P \,\,\, \text{and}\qquad\qquad\,
$$ 
$$
f^{2n+2k+3}(u_{n,k+1}') = f(d) \,\,\, \text{and} \,\,\, f^{m+2(n+k+1)}(\nu_{n,k}) = \min P \,\,\, \text{respectively}\,\,\,
$$
and Lemma 8 that all periodic points of $f$ in $[u_{n,k}', u_n]$ of odd periods have least periods $\ge 2n+2k+3$ to obtain that, 
$$
\text{for each} \,\,\, i \ge 1, \, \text{both} \,\,\, p_{m+2n+2k+2i}^{\, (n,k,3)} \,\,\, \text{and} \,\,\, q_{m+2n+2k+2i}^{\, (n,k,3)} \,\,\, \text{are period-}(m+2n+2k+2i) \,\,\, \text{points of} \,\,\, f.
$$ 

To find the least periods of $c_{2n+2k+2}^{\, (n,k,3)}$'s with respect to $f$, we shall use the following result which can be viewed as a continuation of Lemma 8 from the interval $[u_{n,k}', u_n]$ to the subinterval $[u_{n,k}', u_{n,k,i}]$ (and later to the subinterval $[u_{n,k,i,j}, u_{n,k,i}]$ and then continue on to the higher layers of the {\it basic} tower.  Unfortunately, we do not have similar result on the counterpart interval $[\bar u_{n,k,i}, \bar u_{n,k}]$ and so on):

\noindent
{\bf Lemma 12.}
{\it For each $n \ge 1$, $k \ge 1$, $i \ge 1$, on the interval $[u_{n,k}', u_{n,k,i}] \, (\subset [u_{n,k}', u_{n,k+1}'])$, $f$ has no periodic points of odd periods $\le 2n+2k+2i+1$, nor has periodic points of even periods $\le 2n+2k+2i$ except period-$(2n+2k+2i)$ points and possibly period-$(2n+2k)$ or period-$(2n)$ points.}

\noindent
{\it Proof.}
Suppose $f$ has a period-$(2n+2k+2j)$ point, say $c_{2n+2k+2j}^*$, in $[u_{n,k}', u_{n,k,i}]$ for some $1 \le j \le i-1$.  Then since $f^{2n+2k+2j+2}\big([u_{n,k}', c_{2n+2k+2j}^*]\big) \supset f^2\big([c_{2n+2k+2j}^*, z_0]\big) \supset f^2\big([v, z_0]\big) \supset [\min P, z_0] \supset \{ d \}$, there exists a point $u_{n,k,j+1}^*$ in $[u_{n,k}', c_{2n+2k+2j}^*) \, (\subset [u_{n,k}', u_{n,k,i}])$ such that $f^{2n+2k+2j+2}(u_{n,k,j+1}^*) = d$.  So, $u_{n,k,j+1}^* < u_{n,k,i}$.  Since $n+k+j+1 \le n+k+i$, this contradicts the fact that $u_{n,k,i} \le u_{n,k,j+1} \,\, (\le u_{n,k,j+1}^*)$.  This, combined with Lemma 4(1) and Lemma 8, implies that $f$ has no periodic points of even periods $\le 2n+2k+2i-2$ except possibly period-$(2n)$ or period-$(2n+2k)$ points.  Furthermore, since $f(x) > z \ge z_0$ for all $d < x < z_0$, it follows from the above $(\dagger\dagger)$ that $f$ has no periodic points of odd periods $\le 2n+2k+2i+1$.  

On the other hand, since $f^{2n+2k+2i}(u_{n,k}') - u_{n,k}' = f^{2i}(d) - u_{n,k}' = z_0 - u_{n,k}' > 0$ and $f^{2n+2k+2i}(u_{n,k,i}) - u_{n,k,i} = d - u_{n,k,i} < 0$, the point $q_{2n+2k+2i} = \max \big\{ u_{n,k}' \le x \le u_{n,k,i}: f^{2n+2k+2i}(x) = x \big\}$ exists and, since we have just shown in the previous paragraph that $f$ has no periodic points of least periods $\le 2n+2k+2i-1$, is a period-$(2n+2k+2i)$ point of $f$ in $[u_{n,k}', u_{n,k,i}]$.
$\hfill\square$

We now use Lemma 12 to find the least period of $c_{2n+2k+2i}^{\,\,(n,k,3)}$ with respect to $f$:

Let $\ell\big(c_{2n+2k+2i}^{\,\,(n,k,3)}\big)$ denote the least period of $c_{2n+2k+2i}^{\,\,(n,k,3)}$ with respect to $f$.  Suppose $\ell\big(c_{2n+2k+2i}^{\,\,(n,k,3)}\big) < 2n+2k+2i$ for some $i \ge 1$.  Then by Lemma 12, $\ell\big(c_{2n+2k+2i}^{\,\,(n,k,3)}\big) = 2n$ or $\ell\big(c_{2n+2k+2i}^{\,\,(n,k,3)}\big) = 2n+2k$.
\begin{multline*}
$$
\text{If} \,\,\, 2n+2k+2i \,\,\, \text{is not a multiple of} \,\,\, 2n \,\,\, \text{nor a multiple of} \,\,\, 2n+2k, \, \text{then} \\ 
c_{2n+2k+2i}^{\,\,(n,k,3)} \,\,\, \text{is a period-}(2n+2k+2i) \,\,\, \text{point of} \,\,\, f.
$$
\end{multline*}

If $\ell\big(c_{2n+2k+2i}^{\,\,(n,k,3)}\big) = 2n+2k$ and $(2n+2k+2i)/(2n+2k) = r > 2$, then $(r-1)(2n+2k) > 2n+2k$.  So, $(r-1)(2n+2k) \ge 2n+2k+2$.  We apply Lemma 3 with 
$$
f^{(r-1)(2n+2k)}(u_{n,k}') \in \big\{ z_0, f(z_0) \big\} \,\,\, \text{and} \,\,\, f^{(r-1)(2n+2k)}\big(c_{2n+2k+2i}^{\,\,(n,,3)}\big) = c_{2n+2k+2i}^{\,\,(n,k,3)}
$$
to obtain a periodic point, say $c_{2n+2k+2i}^{**}$, of $f$ such that $u_{n,k}' < c_{2n+2k+2i}^{**} < c_{2n+2k+2i}^{\,\,(n,,k,3)} < u_{n,k,1}$ and $f^{2n+2k+2i}(c_{2n+2k+2i}^{**}) = c_{2n+2k+2i}^{**}$ which contradicts the minimality of $c_{2n+2k+2i}^{\,\,(n,k,3)}$ in $[u_{n,k}', u_{n,k,1}]$.  

So, if $\ell\big(c_{2n+2k+2i}^{\,\,(n,k,3)}\big) = 2n+2k$ and $\ell\big(c_{2n+2k+2i}^{\,\,(n,k,3)}\big) < 2n+2k+2i$, then we must  have $(2n+2k+2i)/(2n+2k) = 2$.  Thus, $c_{4n+4k}^{\,\,(n,k,3)}$ is either a period-{\small $(2n+2k)$} or a period-{\small $(4n+4k)$} point of $f$.

If $\ell\big(c_{2n+2k+2i}^{\,\,(n,k,3)}\big) = 2n$ and $i \ge n+1$, then $(2n+2k+2i)/(2n) = s > 2$.  So, $2n(s-1) = (2n+2k+2i) - 2n = 2k+2i$ $\ge 2k + 2(n+1) = 2n+2k+2$.  By applying Lemma 3 with
$$
f^{2n(s-1)}(u_{n,k}') = z_0 \,\,\, \text{and} \,\,\, f^{2n(s-1)}\big(c_{2n+2k+2i}^{\,\,(n,k,3)}\big) = c_{2n+2k+2i}^{\,\,(n,k,3)},
$$
there exists a periodic point, say $c_{2n+2k+2i}^{***}$, of $f$ such that $u_{n,k}' < c_{2n+2k+2i}^{***} < c_{2n+2k+2i}^{\,\,(n,k,3)} < u_{n,k,1}$ and $f^{2n+2k+2i}(c_{2n+2k+2i}^{***}) = c_{2n+2k+2i}^{***}$ which contradicts the minimality of $c_{2n+2k+2i}^{\,\,(n,k,3)}$ in $[u_{n,k}', u_{n,k,1}]$.  

If $\ell\big(c_{2n+2k+2i}^{\,\,(n,k,3)}\big) = 2n$ and $1 \le i \le n$, then, there can have exactly one $i$ such that $1 \le i \le n$ and $2n+2k+2i$ is a multiple of $2n$.  For such an $i$, $c_{2n+2k+2i}^{\,\,(n,k,3)}$ is either a period-$(2n)$ or a period-$(2n+2k+2i)$ point of $f$.  For all other $i$ such that $1 \le i \le n$ and $2n+2k+2i$ is not a multiple of $2n$, $c_{2n+2k+2i}^{\,\,(n,k,3)}$ is a period-$(2n+2k+2i)$ point of $f$.  

In summary, we have shown that  
\begin{multline*}
$$
\text{for} \,\,\, i = \imath, \, \text{the {\it unique} integer in} \,\,\, [1, n] \,\,\, \text{such that} \,\,\, 2n+2k+2\imath \,\,\, \text{is a multiple of} \,\,\, 2n, \\ 
c_{2n+2k+2\imath}^{\,\,(n,k,3)} \,\,\, \text{is either a period-}(2n) \,\,\, \text{point or a period-}(2n+2k+2\imath) \,\,\, \text{point of} \,\,\, f, \, \\ 
\text{for} \,\,\, i = n+k, \hspace{5.35in} \\ \qquad\, c_{2n+2k+2i}^{\,\,(n,k,3)} \,\,\, \text{is either a period-}(2n+2k) \,\,\, \text{point or a period-}(2n+2k+2i) \,\,\, \text{point of} \,\,\, f, \\
\text{for each} \,\,\, i \ge 1 \,\,\, \text{and} \,\,\, i \notin \{ \imath, n+k \}, \hspace{3.935in} \\ c_{2n+2k+2i}^{\,\,(n,k,3)} \,\,\, \text{is a period-}(2n+2k+2i) \,\,\, \text{point of} \,\,\, f.\hspace{2.5in}\,\,
$$
\end{multline*}

\vspace{.1in}

\noindent
{\bf 3.3 A compartment of the third layer of the {\it basic} tower in $[\breve \mu_{m,n,k}', \breve \mu_{m,n,k+1}'] \, (\subset [\breve \mu_{m,n+1}, \breve \mu_{m,n}] \subset [\breve u_0, v] \subset [u_1, v])$.}

For each $n \ge 1$ and $k \ge 1$, we consider the interval $[\breve \mu_{m,n,k}', \breve \mu_{m,n,k+1}'] \, (\subset [\breve \mu_{m,n+1}, \breve \mu_{m,n}] \subset [\breve u_0, v] \subset [u_1, v])$.  Let $\breve \nu_{m,n,k}$ be a point in $[\breve \mu_{m,n,k}', \breve \mu_{m,n}]$ such that $f^{m+2n+2k}(\breve \nu_{m,n,k}) = v$.  It turns out that $\breve \mu_{m,n,k}' < \breve \nu_{m,n,k} < \breve \mu_{m,n,k+1}'$.  For each $i \ge 1$, let, with $x$-coordinates moving from point $\breve \mu_{m,n,k}$ via point $\breve \nu_{m,n,k}$ to point $\breve \mu_{m,n,k+1}'$,
{\large 
$$\quad \breve \mu_{m,n,k,i} = \min \{ \breve \mu_{m,n,k}' \le x \le \breve \mu_{m,n,k+1}' : f^{m+2n+2k+2i}(x) = d \},\qquad\qquad\qquad$$ 
$$\quad\,\,\,\, \breve p_{m+2n+2k+2i}^{\, (n,k,3)} = \min \{ \breve \mu_{m,n,k}' \le x \le \breve \mu_{m,n,k+1}' : f^{m+2n+2k+2i}(x) = x \},\qquad\qquad\quad$$ 
$$\, \breve c_{2m+2n+2k+2i}^{\, (n,k,3)} = \min \{ \breve \mu_{m,n,k,1} \le x \le \breve \mu_{m,n,k+1}' : f^{2m+2n+2k+2i}(x) = x \} \,\,\, \text{and}\,\,$$ 
$$\,\,\, \breve c_{2m+2n+2k+2i}'^{\, (n,k,3)} = \max \{ \breve \mu_{m,n,k}' \le x \le \breve \mu_{m,n,k+1}' : f^{2m+2n+2k+2i}(x) = x \}.\qquad\,\,\,\,\,
$$}
For each fixed $n \ge 1$ and $k \ge 1$, the collection of all these periodic points $\breve p_{m+2n+2k+2i}^{\, (n,k,3)}$, $\breve c_{2m+2n+2k+2i}^{\, (n,k,3)}$, $\breve c_{2m+2n+2k+2i}'^{\, (n,k,3)}$, $i \ge 1$ is called a {\it compartment} of the third layer of the {\it basic} tower of periodic points of $f$ associated with $P$.  

\noindent It is easy to see that we have 
\begin{multline*}
$$
\breve \mu_{m,n,k}' < \cdots < \breve p_{m+2n+2k+4}^{\, (n,k,3)} < \breve \mu_{m,n,k,2} < \breve p_{m+2n+2k+2}^{\, (n,k,3)} < \breve \mu_{m,n,k,1} < \\ 
\cdots < \breve c_{2m+2n+2k+4}^{\, (n,k,3)} < \breve c_{2m+2n+2k+2}^{\, (n,k,3)} < \breve \nu_{m,n,k} < \breve c_{2m+2n+2k+2}'^{\, (n,k,3)} < \breve c_{2m+2n+2k+4}'^{\, (n,k,3)} < \cdots < \breve \mu_{m,n,k+1}',
$$
\end{multline*}

We now prove the following result which can be viewed as a continuation of Lemma 9 from the interval $[\breve u_0', \breve \mu_{m,n}]$ to the subinterval $[\breve \mu_{m,n,k}', \breve \mu_{m,n}]$ (and later to the subinterval $[\breve \mu_{m,n,k}', \breve \mu_{m,n,k,i}]$ and then continue on to the higher layers of the {\it basic} tower and so on):

\noindent
{\bf Lemma 13.}
{\it For each $n \ge 1$ and $k \ge 1$, all periodic points of $f$ in $[\breve \mu_{m,n,k}', \breve \mu_{m,n}]$ of odd periods have least periods $\ge m+2n+2k$.  (Note that, in Lemma 11, we have $\ge m+2n+2k+2$).}

\noindent
{\it Proof.} 
Recall that, by Lemma 9, all periodic points of $f$ in $[\breve u_0', \breve \mu_{m,n}]$ $(\subset [\breve \mu_{m,n+1}, \breve \mu_{m,n}])$ of odd periods have least periods $\ge m+2n$.  Since 
$$
f^{m+2n}(\breve u_0') - \breve u_0' = f(z_0) - \breve u_0' \ge z_0 - \breve u_0' > 0 \,\,\, \text{and} \,\,\, f^{m+2n}(\breve \mu_{m,n}) - \breve \mu_{m,n} = d - \breve \mu_{m,n} < 0,
$$
the point
$$
\breve q_{m+2n} = \max \big\{\breve u_0' \le x \le \breve \mu_{m,n}: f^{m+2n}(x) = x \big\}
$$
exists.  By the choice of $\breve u_0$, we have $f^2(x) < d$ for all $\breve u_0 < x < v$.  Consequently, since 
$$
f^{m+2n+2}(\breve q_{m+2n}) - d = f^2(\breve q_{m+2n}) - d < 0 \,\,\, \text{and} \,\,\, f^{m+2n+2}(\breve \mu_{m,n}) - d = z_0 - d > 0,
$$
we obtain that the point $\breve \mu_{m,n,1}' = \max \big\{ \breve u_0' \le x \le \breve \mu_{m,n}: f^{m+2n+2}(x) = d \big\}$ satisfies $\breve q_{m+2n} < \breve \mu_{m,n,1}'$ $(< \breve \mu_{m,n,2}' < \breve \mu_{m,n,3}' < \cdots < \breve \mu_{m,n})$.  This, combined with Lemma 9 and the maximality of $\breve q_{m+2n}$ in $[\breve u_0', \breve \mu_{m,n}]$, implies that 
$$
\text{all periodic points of} \,\,\, f \,\,\, \text{in} \,\,\, [\breve \mu_{m,n,1}', \breve \mu_{m,n}] \,\,\, \text{of odd periods have least periods} \, \ge m+2n+2.
$$
In particular, suppose $f$ had a periodic point $\breve p$ of odd period $\ell(\breve p) \le m+2n+2k-2$ in $[\breve \mu_{m,n,k}', \breve \mu_{m,n}]$.  Then $\ell(\breve p) \ge m+2n+2$ and since $\ell(\breve p)$ is odd, we have $f^{\ell(\breve p)}(\breve \mu_{m,n}) = z_0$.  Since $f^{\ell(\breve p)}(\breve p) = \breve p$ and $m+2n+2k-2 \ge \ell(\breve p)$, it follows from Lemma 3 that $f$ has a period-$(m+2n+2k-2)$ point $\breve w_{m+2n+2k-2}$ in $[\breve p, \breve \mu_{m,n}]$ $(\subset [\breve \mu_{m,n,k}', \breve \mu_{m,n}] \subset [\breve \mu_{m,n+1}, \breve \mu_{m,n}])$.  By the choice of $\breve u_0'$, we have $f^2(x) < d$ for all $\breve u_0' < x < v$.  Since 
$$
f^{m+2n+2k}(\breve w_{m+2n+2k-2}) - d = f^2(\breve w_{m+2n+2k-2}) - d < 0 \,\,\, \text{and} \,\,\, f^{m+2n+2k}(\breve \mu_{m,n}) - d = z_0 - d > 0
$$
there is a point $\breve \mu_{m,n,k}^{\,*}$ in $[\breve w_{m+2n+2k-2}, \breve \mu_{m,n}]$ such that $f^{m+2n+2k}(\breve \mu_{m,n,k}^{\,*}) = d$.  Since $\breve \mu_{m,n+1} < \breve \mu_{m,n,k}' < \breve p \le \breve w_{m+2n+2k-2} < \breve \mu_{m,n,k}^{\,*} < \breve \mu_{m,n}$, this contradicts the maximality of $\breve \mu_{m,n,k}'$ in $[\breve \mu_{m,n+1}, \breve \mu_{m,n}]$.  This completes the proof.
$\hfill\square$

Now since 
$$
f^{m+2n+2k}(\breve \mu_{m,n,k}') - \breve \mu_{m,n,k}' = d - \breve \mu_{m,n,k}' < 0 \,\,\, \text{and} \,\,\, f^{m+2n+2k}(\breve \mu_{m,n}) - \breve \mu_{m,n} = z_0 - \breve \mu_{m,n} > 0
$$
the point 
$$
\breve p_{m+2n+2k}^* = \min \big\{ \breve \mu_{m,n,k}' \le x \le \breve \mu_{m,n}: f^{m+2n+2k}(x) = x \big\}
$$
exists.  Since, by the choice of $\breve u_0$, we have $f^2(x) < d$ for all $\breve u_0 < x < v$, and so, 
$$
f^{m+2n+2k+2}(\breve p_{m+2n+2k}^*) - d = f^2(\breve p_{m+2n+2k}^*) - d < 0 \,\,\, \text{and} \,\,\, f^{m+2n+2k+2}(\breve \mu_{m,n}) - d = z_0 - d > 0.
$$
Therefore, we have
$$
(\breve \mu_{m,n,k}' <) \,\, \breve p_{m+2n+2k}^* < \breve \mu_{m,n,k+1}' = \max \big\{ \breve \mu_{m,n,k}' \le x \le \breve \mu_{m,n}: f^{m+2n+2k+2}(x) = d \big\} < \breve \mu_{m,n}.
$$
This, combined with Lemma 13, implies that all periodic point of $f$ in $[\breve \mu_{m,n,k}', \breve p_{m+2n+2k}^*)$ of odd periods have least periods $\ge m+2n+2k+2$.  

On the other hand, by the choice of $\breve u_0'$, we have $f^2(x) < d$ for all $\breve u_0 < x < v$.  Since 
\begin{multline*}
$$
f^{m+2n+2k+2}(\breve \mu_{m,n,k}') - \breve \mu_{m,n,k}' = z_0 - \breve \mu_{m,n,k}' > 0 \,\,\, \text{and} \\ 
f^{m+2n+2k+2}(\breve p_{m+2n+2k}^*) - \breve p_{m+2n+2k}^* = f^2(\breve p_{m+2n+2k}^*) - \breve p_{m+2n+2k}^* < d - \breve p_{m+2n+2k}^* < 0
$$
\end{multline*}
the point $\breve p_{m+2n+2k+2} = \min \big\{ \breve \mu_{m,n,k}' \le x < \breve p_{m+2n+2k}^*: f^{m+2n+2k+2}(x) = x \big\}$ $(< \breve p_{m+2n+2k}^*)$ exists and must be a period-$(m+2n+2k+2)$ point of $f$.  

Since $f^{m+2n+2k+2}(\breve \mu_{m,n,k}') = z_0$ and $f^{m+2n+2k+2}(\breve p_{m+2n+2k+2}) = \breve p_{m+2n+2k+2}$, and since we have just shown in the above that {\it all periodic point of $f$ in the interval $[\breve \mu_{m,n,k}', \breve p_{m+2n+2k}^*)$ $(\supset [\breve \mu_{m,n,k}', \breve p_{m+2n+2k+2}])$ of odd periods have least periods $\ge m+2n+2k+2$}, it follows from Lemma 3 that, 
$$
\text{for each} \,\,\, i \ge 1, \, \text{the point} \,\,\, \breve p_{m+2n+2k+2i} \,\,\, \text{is a period-}(m+2n+2k+2i) \,\,\, \text{point of} \,\,\, f.
$$

As for the periods of $\breve c_{2m+2n+2k+2i}^{\, (n,k,3)}$'s and $\breve c_{2m+2n+2k+2i}'^{\, (n,k,3)}$'s, we apply Lemma 6 with 
\begin{multline*}
$$
\,\,\,\,\,\qquad\qquad f^{m+2n+2k+3}(\breve \mu_{m,n,k,1}) = f(d) \,\,\, \text{and} \,\,\, f^{2(m+n+k+1)}(\breve \nu_{m,n,k}) = \min P \,\,\, \text{and}\,\, \\ 
f^{m+2n+2k+3}(\breve \mu_{m,n,k+1}') = f(d) \,\,\, \text{and} \,\,\, f^{2(m+n+k+1)}(\breve \nu_{m,n,k}) = \min P \,\,\, \text{respectively}\quad\,
$$
\end{multline*}
to obtain that, for each $i \ge 0$ such that $(m+n+k+1)+i \ge \max \big\{ m+2n+2k+3, m+2 \big\} = m+2n+2k+3$, both $\breve c_{2m+2n+2k+2+2i}^{\, (n,k,3)}$ and $\breve c_{2m+2n+2k+2+2i}'^{\, (n,k,3)}$ are period-$(2m+2n+2k+2+2i)$ points of $f$, or, equivalently, 
{\small 
$$
\text{for each} \,\, i \ge n+k+3, \, \text{both} \,\, \breve c_{2m+2n+2k+2i}^{\, (n,k,3)} \,\, \text{and} \,\, \breve c_{2m+2n+2k+2i}'^{\, (n,k,3)} \,\, \text{are period-}(2m+2n+2k+2i) \,\, \text{points of} \,\, f.\vspace{.1in}
$$}
{\bf 3.4 A compartment of the third layer of the {\it basic} tower in $[\hat \mu_{m,n,k+1}, \hat \mu_{m,n,k}] \, (\subset [\hat \mu_{m,n}', \hat \mu_{m,n+1}'] \subset [v, \hat u_0] \subset [v, z_0])$.}

For each $n \ge 1$ and $k \ge 1$, we consider the interval $[\hat \mu_{m,n,k+1}, \hat \mu_{m,n,k}] \, (\subset [\hat \mu_{m,n}', \hat \mu_{m,n+1}'] \subset [v, \hat u_0] \subset [v, z_0])$.  Let $\hat \nu_{m,n,k}$ be a point in $[\hat \mu_{m,n}', \hat \mu_{m,n,k}]$ such that $f^{m+2n+2k}(\hat \nu_{m,n,k}) = v$.  It turns out that $\hat \mu_{m,n,k+1} < \hat \nu_{m,n,k} < \hat \mu_{m,n,k}$.  For each $i \ge 1$, let, with $x$-coordinates moving from point $\hat \mu_{m,n,k}$ via point $\hat \nu_{m,n,k}$ to point $\hat \mu_{m,n,k+1}$, 
{\large 
$$
\quad \hat \mu_{m,n,k,i}' = \max \{ \hat \mu_{m,n,k+1} \le x \le \hat \mu_{m,n,k} : f^{m+2n+2k+2i}(x) = d \},\qquad\qquad\qquad\,$$ 
$$\quad\,\,\, \hat q_{m+2n+2k+2i}^{\, (n,k,3)} = \max \{ \hat \mu_{m,n,k+1} \le x \le \hat \mu_{m,n,k} : f^{m+2n+2k+2i}(x) = x \},\qquad\qquad\quad$$ 
$$\hat c_{2m+2n+2k+2i}'^{\, (n,k,3)} = \max \{ \hat \mu_{m,n,k+1} \le x \le \hat \mu_{m,n,k,1}' : f^{2m+2n+2k+2i}(x) = x \} \,\,\, \text{and}\,$$ 
$$\hat c_{2m+2n+2k+2i}^{\, (n,k,3)} = \min \{ \hat \mu_{m,n,k+1} \le x \le \hat \mu_{m,n,k} : f^{2m+2n+2k+2i}(x) = x \}.\qquad\,\,\,\,\,
$$}
For each fixed $n \ge 1$ and $k \ge 1$, the collection of all these periodic points $\hat q_{m+2n+2k+2i}^{\, (n,k,3)}$, $\hat c_{2m+2n+2k+2i}'^{\, (n,k,3)}$, $\hat c_{2m+2n+2k+2i}^{\, (n,k,3)}$, $i \ge 1$ is called a {\it compartment} of the third layer of the {\it basic} tower of periodic points of $f$ associated with $P$.

\noindent It is easy to see that we have 
\begin{multline*}
$$
\hat \mu_{m,n,k+1} < \cdots < \hat c_{2m+2n+2k+4}^{\, (n,k,3)} < \hat c_{2m+2n+2k+2}^{\, (n,k,3)} < \hat \nu_{m,n,k} < \hat c_{2m+2n+2k+2}'^{\, (n,k,3)} < \hat c_{2m+2n+2k+4}'^{\, (n,k,3)} < \cdots \\ < \hat \mu_{m,n,k,1}' < \hat q_{m+2n+2k+2}^{\, (n,k,3)} < \hat \mu_{m,n,k,2}' < \hat q_{m+2n+2k+4}^{\, (n,k,3)} < \hat \mu_{m,n,k,3}' < \cdots < \hat \mu_{m,n,k}.
$$
\end{multline*}

We now prove the following result which can be viewed as a continuation of Lemma 10 from the interval $[\hat \mu_{m,n}', \hat u_0]$ to the subinterval $[\hat \mu_{m,n}', \hat \mu_{m,n,k}]$ (and later to the subinterval $[\hat \mu_{m,n,k,i}, \hat \mu_{m,n,k}]$ and then continue on to the higher layers of the {\it basic} tower and so on):

\noindent
{\bf Lemma 14.}
{\it For each $n \ge 1$ and $k \ge 1$, all periodic points of $f$ in $[\hat \mu_{m,n}', \hat \mu_{m,n,k}]$ of odd periods have least periods $\ge m+2n+2k$.  (Note that, in Lemma 11, we have $\ge m+2n+2k+2$).}

\noindent
{\it Proof.}
The proof is similar to that of Lemma 13 but different from that of Lemma 11.  Recall that, by Lemma 10, all periodic points of $f$ in $[\hat \mu_{m,n}', \hat u_0]$ $(\subset [\hat \mu_{m,n}', \hat \mu_{m,n+1}'])$ of odd periods have least periods $\ge m+2n$.  Since 
$$
f^{m+2n}(\hat \mu_{m,n}') - \hat \mu_{m,n}' = d - \hat \mu_{m,n}' < 0 \,\,\, \text{and} \,\,\, f^{m+2n}(\hat u_0) - \hat u_0 = f(z_0) - \hat u_0 > 0,
$$
the point $\hat p_{m+2n} = \min \big\{\hat \mu_{m,n}' \le x \le \hat u_0: f^{m+2n}(x) = x \big\}$ exists.  By the choice of $\hat u_0$, we have $f^2(x) < d$ for all $v < x < \hat u_0$.  Consequently, since 
$$
f^{m+2n+2}(\hat u_{m,n}') - d = z_0 - d > 0 \,\,\, \text{and} \,\,\, f^{m+2n+2}(\hat p_{m+2n}) - d = f^2(\hat p_{m+2n}) - d < 0,
$$
we obtain that the point $\hat \mu_{m,n,1} = \min \big\{ \hat \mu_{m,n}' \le x \le \hat u_0: f^{m+2n+2}(x) = d \big\}$ satisfies $(\hat \mu_{m,n}' < \cdots < \hat \mu_{m,n,3} < \hat \mu_{m,n,2} <)$ $\hat \mu_{m,n,1} < \hat p_{m+2n}$.  This, combined with Lemma 10, implies that 
$$
\text{all periodic points of} \,\,\, f \,\,\, \text{in} \,\,\, [\hat \mu_{m,n}', \hat \mu_{m,n,k}] \,\,\, \text{of odd periods have least periods} \, \ge m+2n+2.
$$
Consequently, suppose $f$ had a periodic point $\hat p$ of odd period $\ell(\hat p) \le m+2n+2k-2$ in $[\hat \mu_{m,n}', \hat \mu_{m,n,k}]$.  Then $\ell(\hat p) \ge m+2n+2$ and since $\ell(\hat p)$ is odd, we have $f^{\ell(\hat p)}(\hat \mu_{m,n}') = z_0$.  Since $f^{\ell(\hat p)}(\hat p) = \hat p$ and $m+2n+2k-2 \ge \ell(\hat p)$, it follows from Lemma 3 that $f$ has a period-$(m+2n+2k-2)$ point $\hat w_{m+2n+2k-2}$ in $[\hat \mu_{m,n}', \hat p]$ $(\subset [\hat \mu_{m,n}', \hat \mu_{m,n,k}] \subset [\hat \mu_{m,n}', \hat \mu_{m,n+1}'])$.  By the choice of $\hat u_0$, we have $f^2(x) < d$ for all $v < x < \hat u_0$.  Since 
$$
f^{m+2n+2k}(\hat w_{m+2n+2k-2}) - d = f^2(\hat w_{m+2n+2k-2}) - d < 0 \,\,\, \text{and} \,\,\, f^{m+2n+2k}(\hat \mu_{m,n}') - d = z_0 - d > 0
$$
there is a point $\hat \mu_{m,n,k}^{\,*}$ in $[\hat \mu_{m,n}', \hat w_{m+2n+2k-2}]$ such that $f^{m+2n+2k}(\hat \mu_{m,n,k}^{\,*}) = d$.  Since $\hat \mu_{m,n}' < \hat \mu_{m,n,k}^{\,*} < \hat w_{m+2n+2k-2} \le \hat p < \hat \mu_{m,n,k} < \hat \mu_{m,n+1}'$, this contradicts the minimality of $\hat \mu_{m,n,k}$ in $[\hat \mu_{m,n}', \hat \mu_{m,n+1}']$.  This completes the proof.
$\hfill\square$

Now since
$$
f^{m+2n+2k}(\hat \mu_{m,n}') - \hat \mu_{m,n}' = z_0 - \hat \mu_{m,n}' > 0 \,\,\, \text{and} \,\,\, f^{m+2n+2k}(\hat \mu_{m,n,k}) - \hat \mu_{m,n,k} = d - \hat \mu_{m,n,k} < 0
$$
the point 
$$
\hat q_{m+2n+2k}^* = \max \big\{ \hat \mu_{m,n}' \le x \le \hat \mu_{m,n,k}: f^{m+2n+2k}(x) = x \big\}
$$
exists.  Since, by the choice of $\hat u_0$, we have $f^2(x) < d$ for all $v < x < \hat u_0$, and so, 
$$
f^{m+2n+2k+2}(\hat \mu_{m,n}') - d = z_0 - d > 0 \,\,\, \text{and} \,\,\, f^{m+2n+2k+2}(\hat q_{m+2n+2k}^*) - d = f^2(\hat q_{m+2n+2k}^*) - d < 0.
$$
Therefore, we have
$$
(\hat \mu_{m,n}' <) \,\, \hat \mu_{m,n,k+1} = \min \big\{ \hat \mu_{m,n}' \le x \le \hat \mu_{m,n+1}': f^{m+2n+2k+2}(x) = d \big\} < \hat q_{m+2n+2k}^* < \hat \mu_{m,n,k}.
$$
This, combined with Lemma 14, implies that all periodic point of $f$ in $(\hat q_{m+2n+2k}^*, \hat \mu_{m,n,k}]$ of odd periods have least periods $\ge m+2n+2k+2$.  

On the other hand, we have $f^2(x) < x$ for all $v < x < z_0$.  Since 
\begin{multline*}
$$
f^{m+2n+2k+2}(\hat q_{m+2n+2k}^*) - \hat q_{m+2n+2k}^* = f^2(\hat q_{m+2n+2k}^*) - \hat q_{m+2n+2k}^* < 0 \,\,\, \text{and} \\ 
f^{m+2n+2k+2}(\hat \mu_{m,n,k}) - \hat \mu_{m,n,k} = z_0 - \hat \mu_{m,n,k} > 0,
$$
\end{multline*}
the point $\hat p_{m+2n+2k+2} = \min \big\{ \hat q_{m+2n+2k}^* \le x \le \hat \mu_{m,n,k}: f^{m+2n+2k+2}(x) = x \big\}$ $(> \hat q_{m+2n+2k}^*)$ exists and must be a period-$(m+2n+2k+2)$ point of $f$.  

Since $f^{m+2n+2k+2}(\hat \mu_{m,n,k}) = z_0$ and $f^{m+2n+2k+2}(\hat p_{m+2n+2k+2}) = \hat p_{m+2n+2k+2}$, and since we have just shown in the above that {\it all periodic point of $f$ in the interval $(\hat q_{m+2n+2k}^*, \hat \mu_{m,n,k}]$ $(\supset [\hat p_{m+2n+2k+2}, \hat \mu_{m,n,k}])$ of odd periods have least periods $\ge m+2n+2k+2$}, it follows from Lemma 3 that, 
$$
\text{for each} \,\,\, i \ge 1, \, \text{the point} \,\,\, \hat q_{m+2n+2k+2i} \,\,\, \text{is a period-}(m+2n+2k+2i) \,\,\, \text{point of} \,\,\, f.
$$
\indent As for the periods of $\hat c_{2m+2n+2k+2i}^{\, (n,k,3)}$'s and $\hat c_{2m+2n+2k+2i}'^{\, (n,k,3)}$'s, we apply Lemma 6 with 
\begin{multline*}
$$
\,\,\,\,\qquad\qquad f^{m+2n+2k+3}(\hat \mu_{m,n,k,1}') = f(d) \,\,\, \text{and} \,\,\, f^{2(m+n+k+1)}(\hat \nu_{m,n,k}) = \min P \,\,\, \text{and}\,\, \\ 
f^{m+2n+2k+3}(\hat \mu_{m,n,k+1}) = f(d) \,\,\, \text{and} \,\,\, f^{2(m+n+k+1)}(\hat \nu_{m,n,k}) = \min P \,\,\, \text{respectively}\quad\,\,
$$
\end{multline*}
to obtain that, for each $i \ge 0$ such that $(m+n+k+1)+i \ge \max \big\{ m+2n+2k+3, m+2 \big\} = m+2n+2k+3$, both $\hat c_{2m+2n+2k+2+2i}'^{\, (n,k,3)}$ and $\hat c_{2m+2n+2k+2+2i}^{\, (n,k,3)}$ are period-$(2m+2n+2k+2+2i)$ points of $f$, or, equivalently, in $[\hat \mu_{m,n,k+1}, \hat \mu_{m,n,k}] \, (\subset [\hat \mu_{m,n}', \hat \mu_{m,n+1}'] \subset [v, \hat u_0] \subset [v, z_0])$,
{\small 
$$
\text{for each} \,\, i \ge n+k+3, \\ \text{both} \,\,\, \hat c_{2m+2n+2k+2i}'^{\, (n,k,3)} \,\,\, \text{and} \,\,\, \hat c_{2m+2n+2k+2i}^{\, (n,k,3)} \,\,\, \text{are period-}(2m+2n+2k+2i) \,\,\, \text{points of} \,\,\, f.
$$}

\vspace{.1in}

\noindent
{\bf 3.5 A compartment of the third layer of the {\it basic} tower in $[\bar u_{n,k+1}, \bar u_{n,k}] \, (\subset [\bar u_n', \bar u_{n+1}'] \subset [v, z_0])$.}

For each $n \ge 1$ and $k \ge 1$, we consider the interval $[\bar u_{n,k+1}, \bar u_{n,k}] \, (\subset [\bar u_n', \bar u_{n+1}'] \subset [v, z_0])$.  Let $\bar \nu_{n,k}$ be a point in $[\bar u_n', \bar u_{n,k}]$ such that $f^{2n+2k}(\bar \nu_{n,k}) = v$.  Then the point $\bar \nu_{n,k}$ happens to satisfy that $\bar u_{n,k+1} < \bar \nu_{n,k} < \bar u_{n,k}$.  For each $i \ge 1$, let, with $x$-coordinates moving from point $\bar u_{n,k}$ via point $\bar \nu_{n,k}$ to point $\bar u_{n,k+1}$, 
{\large 
$$\bar u_{n,k,i}' = \max \{ \bar u_{n,k+1} \le x \le \bar u_{n,k} : f^{2n+2k+2i}(x) = d \},\qquad\qquad\qquad\,\,\,\,$$ 
$$\bar c_{2n+2k+2i}'^{\, (n,k,3)} = \max \{ \bar u_{n,k+1} \le x \le \bar u_{n,k} : f^{2n+2k+2i}(x) = x \},\qquad\qquad\quad$$ 
$$\bar q_{m+2n+2k+2i}^{\, (n,k,3)} = \max \{ \bar u_{n,k+1} \le x \le \bar u_{n,k,1}' : f^{m+2n+2k+2i}(x) = x \} \,\,\, \text{and}\,$$ 
$$\, \bar p_{m+2n+2k+2i}^{\, (n,k,3)} = \min \{ \bar u_{n,k+1} \le x \le \bar u_{n,k} : f^{m+2n+2k+2i}(x) = x \}.\qquad\,\,\,\,\,
$$}
For each fixed $n \ge 1$ and $k \ge 1$, the collection of all these periodic points $\bar c_{2n+2k+2i}'^{\, (n,k,3)}$, $\bar q_{m+2n+2k+2i}^{\, (n,k,3)}$, $\bar p_{m+2n+2k+2i}^{\, (n,k,3)}$, $i \ge 1$ is called a {\it compartment} of the third layer of the {\it basic} tower of periodic points of $f$ associated with $P$.

\noindent 
It is easy to see that
$$
\bar u_{n,k+1} \,\,\, < \,\,\, \bar \nu_{n,k} \,\,\, < \,\,\, \bar u_{n,k,1}' \,\,\, < \,\,\, \bar u_{n,k,2}' \,\,\, < \,\,\, \bar u_{n,k,3}' \,\,\, < \,\,\, \cdots \,\,\, < \,\,\, \bar u_{n,k} \quad \text{and},
$$
combined with the fact {$(\overline {\dagger\dagger})$} that $d < f^{2j}(x) < z_0$ for all $\bar u_n' \le x \le \bar u_{n,k}$ and all $1 \le j \le n+k$,
{\large 
$$
\quad d < f^{2j}(x) < z_0 \,\,\, \text{for all} \,\,\, 1 \le j \le n+k+i \,\,\, \text{and all} \,\,\, \bar u_{n,k,i}' \le x \le \bar u_{n,k}, \,\,\,\, ({\overline {\dagger\dagger\dagger}})
$$}
\begin{multline*}
$$
\hspace{.2in} \bar u_{n,k+1} \,\,\, < \,\,\, \bar u_{n,k,1}' \,\,\, < \,\,\, \bar c_{2n+2k+2}'^{\, (n,k,3)} \,\,\, < \,\,\, \bar c_{2n+2k+4}'^{\, (n,k,3)} \,\,\, < \,\,\, \bar c_{2n+2k+6}'^{\, (n,k,3)} \,\,\, < \,\,\, \cdots \,\,\, < \,\,\, \bar u_{n,k}, \hspace{.2in} \\ 
\bar u_{n,k+1} \,\, < \,\, \cdots \,\, < \,\, \bar p_{m+2n+2k+6}^{\, (n,k,3)} \,\, < \,\, \bar p_{m+2n+2k+4}^{\, (n,k,3)} \,\, < \,\, \bar p_{m+2n+2k+2}^{\, (n,k,3)} \,\, < \,\, \bar \nu_{n,k} \\ 
< \,\, \bar q_{m+2n+2k+2}^{\, (n,k,3)} \,\, < \,\, \bar q_{m+2n+2k+4}^{\, (n,k,3)} \,\, < \,\, \bar q_{m+2n+2k+6}^{\, (n,k,3)} \,\, < \,\, \cdots \,\, < \,\, \bar u_{n,k,1}' \,\, < \,\, \bar u_{n,k}.
$$
\end{multline*}

As for the periods of $\bar c_{2n+2k+2i}'^{\, (n,k,3)}$'s, since we do not have a result similar to Lemma 12, we can apply Lemma 6 with 
$$
f^{2(n+k+1)}(\bar u_{n,k}) = z_0 \,\,\, \text{and} \,\,\, f^{2(n+k+1)}(\bar \nu_{n,k}) = \min P.
$$
However, since, in this case, we have the following fact that 
$$
\qquad\quad\,\,\, d < f^{2j}(x) < z_0 \,\,\, \text{for all} \,\,\, 1 \le j \le n+k+i \,\,\, \text{and all} \,\,\, \bar u_{n,k,i}' \le x \le \bar u_{n,k}, \,\,\,\,\qquad\qquad\, ({\overline {\dagger\dagger\dagger}})
$$
we can get a result better than that obtained by simply applying Lemma 6.  First, by Lemma 6(4), we obtain that, for each $i \ge 1$,

\noindent
if $(n+k+1)+i$ is even and $\ge 2(n+k+1)$, then $\bar c_{2(n+k+1)+2i}'^{\, (n,k,3)}$ is a period-$(2n+2k+2+2i)$ point of $f$,

\noindent
if $(n+k+1)+i$ is odd and $\ge 2(n+k+1)$, then $\bar c_{2(n+k+1)+2i}'^{\, (n,k,3)}$ is either a period-$(2n+2k+2+2i)$ or an {\it odd} period-$(n+k+1+i)$ point of $f$.

Now it follows from the above $({\overline {\dagger\dagger\dagger}})$ that $\bar c_{2(n+k+1)+2i}'^{\, (n,k,3)}$ cannot be an {\it odd} period-$(n+k+1+i)$ point of $f$.  Consequently, we obtain that, for each $i \ge 1$ such that $(n+k+1)+i \ge 2(n+k+1)$, $\bar c_{2(n+k+1)+2i}'^{\, (n,k,3)}$ is a period-$(2n+2k+2+2i)$ point of $f$, or, equivalently, 
$$
\text{for each} \,\,\, i \ge n+k+2, \, \bar c_{2n+2k+2i}'^{\, (n,k,3)} \,\,\, \text{is a period-}(2n+2k+2i) \,\,\, \text{point of} \,\,\, f.
$$
\indent As for the periods of $\bar q_{m+2n+2k+2i}^{\, (n,k,3)}$'s and $\bar p_{m+2n+2k+2i}^{\, (n,k,3)}$'s, we apply Lemma 5 with 
$$
f^{2n+2k+3}(\bar u_{n,k,1}') = f(d) \,\,\, \text{and} \,\,\, f^{m+2(n+k+1)}(\bar \nu_{n,k}) = \min P \,\,\, \text{and}\qquad\quad\,\,\,\,\,
$$ 
$$
\,\,\, f^{2n+2k+3}(\bar u_{n,k+1}) = f(d) \,\,\, \text{and} \,\,\, f^{m+2(n+k+1)}(\bar \nu_{n,k}) = \min P \,\,\, \text{respectively}\,\,\,
$$
and the above fact $({\overline {\dagger\dagger\dagger}})$ to obtain that, 
$$
\text{for each} \,\, i \ge 1, \, \text{both} \,\,\, \bar q_{m+2n+2k+2i}^{\, (n,k,3)} \,\,\, \text{and} \,\,\, \bar p_{m+2n+2k+2i}^{\, (n,k,3)} \,\,\, \text{are period-}(m+2n+2k+2i) \,\,\, \text{points of} \,\,\, f.
$$
Note that the {\it relative} locations of the periodic points $\bar c_{2n+2k+2i}'^{\, (n,k,3)}$'s with respect to the points $\bar u_{n,k,i}'$'s are not known except that $\bar u_{n,k,i}' < \bar c_{2n+2k+2i}'^{\, (n,k,3)}$ for all $i \ge 1$.

{\large 
We call the collection of all these periodic points, from the (smallest) point $\min P$ (exclusive) to the (largest) point $z_0$ (exclusive), excluding $\tilde q_{m+2n+2k+2}^{\, (n,k,3)}$, 
$$\,\,\,\,\, \tilde c_{2m+2n+2k+2i}^{\, (n,k,3)}, \, \tilde c_{2m+2n+2k+2i}'^{\, (n,k,3)}, \, \tilde q_{m+2n+2k+2+2i}^{\, (n,k,3)}, \, i \ge 1, \, k \ge 1, \, n \ge 1;$$
$$c_{2m+2n+2k+2i}^{\, (n,k,3)}, \, p_{m+2n+2k+2i}^{\, (n,k,3)}, \, q_{m+2n+2k+2i}^{\, (n,k,3)}, \, i \ge 1, \, k \ge 1, \, n \ge 1;\,\,$$
$$\, \breve p_{m+2n+2k+2i}^{\, (n,k,3)}, \, \breve c_{2m+2n+2k+2i}^{\, (n,k,3)}, \, \breve c_{2m+2n+2k+2i}'^{\, (n,k,3)}, \, i \ge 1, \, k \ge 1, \, n \ge 1;$$
$$\, \hat c_{2m+2n+2k+2i}^{\, (n,k,3)}, \, \hat c_{2m+2n+2k+2i}'^{\, (n,k,3)}, \, \hat q_{m+2n+2k+2i}^{\, (n,k,3)}, \, i \ge 1, \, k \ge 1, \, n \ge 1;\,$$
$$\, \bar p_{m+2n+2k+2i}^{\, (n,k,3)}, \, \bar q_{m+2n+2k+2i}^{\, (n,k,3)}, \, \bar c_{m+2n+2k+2i}'^{\, (n,k,3)}, \, i \ge 1, \, k \ge 1, \, n \ge 1 \quad\vspace{.1in}$$
the {\it third layer} of the {\it basic} tower of periodic points of $f$ associated with $P$.}

\pagebreak

\noindent
{\bf $\mathsection$4. The higher layers of the {\it basic} tower of periodic points of $f$ associated with $P$.}

We wrap up the results (from the first layer to the third layer of the {\it basic} tower) obtained so far and arrange them in the order from the interval $[\min P, d]$, through $[d, u_1]$, $[\breve u_0, v]$ $(\subset [u_1, v])$, $[v, \hat u_0]$ $(\subset [v, \bar u_1'])$ to $[\bar u_1', z_0]$.

\vspace{.2in}

\noindent
{\bf 4.1 On the interval $[\min P, d]$}: 

\noindent
The first layer, the starting step, for each $n \ge 0$,
\begin{multline*}
$$
\tilde \mu_{m,n}' = \max \big\{ \min P \le x \le d: f^{m+2n}(x) = d \big\},  \\ 
\tilde q_{m+2n} = \max \big\{ \min P \le x \le d: f^{m+2n}(x) = x \big\}. \qquad\qquad\qquad\qquad\qquad\qquad\qquad\qquad
$$
\end{multline*}
\indent For each $n \ge 0$, $\tilde q_{m+2n}$ is a period-$(m+2n)$ point of $f$.

\vspace{.1in}

\noindent
The second layer, for any fixed $n \ge 0$ and for each $k \ge 1$,
\begin{multline*}
$$
\tilde \mu_{m,n,k} = \min \big\{ \tilde \mu_{m,n}' \le x \le \tilde \mu_{m,n+1}': f^{m+2n+2k}(x) = d \big\}, \\ 
\tilde p_{m+2n+2k}^{\, (n,2)} = \min \big\{ \tilde \mu_{m,n}' \le x \le \tilde \mu_{m,n+1}': f^{m+2n+2k}(x) = x \big\}, \qquad\qquad\qquad \\ 
\qquad\qquad\quad \tilde c_{2m+2n+2k}^{\, (n,2)} = \min \big\{ \tilde \mu_{m,n,1} \le x \le \tilde \mu_{m,n+1}': f^{2m+2n+2k}(x) = x \big\}, \\ 
\tilde c_{2m+2n+2k}'^{\, (n,2)} = \max \big\{ \tilde \mu_{m,n}' \le x \le \tilde \mu_{m,n+1}': f^{2m+2n+2k}(x) = x \big\}.
$$
\end{multline*}
\indent For any fixed $n \ge 0$, we have 
\begin{multline*}
$$
\text{for each} \,\,\, k \ge 1, \, \tilde p_{m+2n+2k}^{\, (n,2)} \,\,\, \text{is a period-}(m+2n+2k) \,\,\, \text{point of} \,\,\, f \,\,\, \text{and}, \\ 
\text{for each} \,\,\, k \ge n+3, \, \text{both} \,\, \tilde c_{2m+2n+2k}^{\, (n,2)} \,\, \text{and} \,\, \tilde c_{2m+2n+2k}'^{\, (n,2)} \,\, \text{are period-}(2m+2n+2k) \,\, \text{points of} \,\, f.
$$
\end{multline*}

\noindent
The third layer, for any fixed $n \ge 0, k \ge 1$ and for each $i \ge 1$,
\begin{multline*}
$$
\tilde \mu_{m,n,k,i}' = \max \big\{ \tilde \mu_{m,n,k+1} \le x \le \tilde \mu_{m,n,k} : f^{m+2n+2k+2i}(x) = d \big\}, \\ 
\tilde q_{m+2n+2k+2i}^{\, (n,k,3)} = \max \big\{ \tilde \mu_{m,n,k+1} \le x \le \tilde \mu_{m,n,k} : f^{m+2n+2k+2i}(x) = x \big\}, \qquad\qquad\,\,\, \\ 
\,\,\,\,\quad\qquad \tilde c_{2m+2n+2k+2i}'^{\, (n,k,3)} = \max \big\{ \tilde \mu_{m,n,k+1} \le x \le \tilde \mu_{m,n,k,1}' : f^{2m+2n+2k+2i}(x) = x \big\}, \\ 
\tilde c_{2m+2n+2k+2i}^{\, (n,k,3)} = \min \big\{ \tilde \mu_{m,n,k+1} \le x \le \tilde \mu_{m,n,k} : f^{2m+2n+2k+2i}(x) = x \big\}.
$$
\end{multline*}
\indent For any fixed $n \ge 0$ and $k \ge 1$, we have 
\begin{multline*}
$$
\text{for each} \,\,\, i \ge 1, \, \tilde q_{m+2n+2k+2i}^{\, (n,k,3)} \,\,\, \text{is a period-}(m+2n+2k+2i) \,\,\, \text{point of} \,\,\, f \,\,\, \text{and}, \\ 
\text{for each} \,\,\, i \ge n+k+3, \,\qquad\qquad\qquad\qquad\qquad\qquad\qquad\qquad\qquad\qquad\qquad\qquad\qquad \\ 
\text{both} \,\,\, \tilde c_{2m+2n+2k+2i}'^{\, (n,k,3)} \,\,\, \text{and} \,\,\, \tilde c_{2m+2n+2k+2i}^{\, (n,k,3)} \,\,\, \text{are period-}(2m+2n+2k+2i) \,\,\, \text{points of} \,\,\, f.
$$
\end{multline*}

\pagebreak

\noindent
{\bf 4.2 On the interval $[d, u_1]$:}

The first layer, the starting step, for each $n \ge 1$,
\begin{multline*}
$$
u_n = \min \big\{ d \le x \le v: f^{2n}(x) = d \big\} \\ 
c_{2n} = \min \big\{ d \le x \le u_n: f^{2n}(x) = x \big\}.\qquad\qquad\qquad\qquad\qquad\qquad\qquad\qquad\qquad\qquad
$$
\end{multline*}
\indent For each $n \ge 1$, $c_{2n}$ is a period-$(2n)$ point of $f$.

\vspace{.1in}

The second layer, for any fixed $n \ge 1$ and for each $k \ge 1$,
\begin{multline*}
$$
u_{n,k}' = \max \big\{ u_{n+1} \le x \le u_n: f^{2n+2k}(x) = d \big\}, \\ 
c_{2n+2k}'^{\, (n,2)} = \max \big\{ u_{n+1} \le x \le u_n: f^{2n+2k}(x) = x \big\},\qquad\qquad\qquad\qquad\,\,\,\,\,\,\,\,\,\,\,\, \\ \,\,\,\qquad\qquad q_{m+2n+2k}^{\, (n,2)} = \max \big\{ u_{n+1} \le x \le u_{n,1}': f^{m+2n+2k}(x) = x \big\}, \\ 
p_{m+2n+2k}^{\, (n,2)} = \min \big\{ u_{n+1} \le x \le u_n: f^{m+2n+2k}(x) = x \big\}.
$$
\end{multline*}
\indent For any fixed $n \ge 1$, we have 
\begin{multline*}
$$
\text{for each} \,\,\, k \ge 1, \, c_{2n+2k}'^{\, (n,2)} \,\,\, \text{is a period-}(2n+2k) \,\,\, \text{(except possibly period-}(4n)) \,\,\, \text{point of} \,\,\, f \,\,\, \text{and}, \\ 
\text{for each} \,\,\, k \ge 1, \, \text{both} \,\,\, q_{m+2n+2k}^{\, (n,2)} \,\,\, \text{and} \,\,\, p_{m+2n+2k}^{\, (n,2)} \,\,\, \text{are period-}(m+2n+2k) \,\,\, \text{points of} \,\,\, f.
$$
\end{multline*}

The third layer, for any fixed $n \ge 1$, $k \ge 1$ and for each $i \ge 1$,
\begin{multline*}
$$
u_{n,k,i} = \min \big\{ u_{n,k} \le x \le u_{n,k+1} : f^{2n+2k+2i}(x) = d \big\}, \\ 
c_{2n+2k+2i}^{\, (n,k,3)} = \min \big\{ u_{n,k} \le x \le u_{n,k+1} : f^{2n+2k+2i}(x) = x \big\},\qquad\qquad\qquad\quad\,\,\,\,\,\,\, \\ \quad\qquad\qquad\qquad\quad\, p_{m+2n+2k+2i}^{\, (n,k,3)} = \min \big\{ u_{n,k,1} \le x \le u_{n,k+1} : f^{m+2n+2k+2i}(x) = x \big\}, \qquad\qquad \\ 
q_{m+2n+2k+2i}^{\, (n,k,3)} = \max \big\{ u_{n,k} \le x \le u_{n,k+1} : f^{m+2n+2k+2i}(x) = x \big\}.
$$
\end{multline*}
\indent For any fixed $n \ge 1$ and $k \ge 1$, we have 
\begin{multline*}
$$
\text{for each} \,\,\, i \ge 1 \,\,\, \text{such that} \,\,\, i \notin \{ \imath, n+k \},  \, \text{where} \,\,\, \imath \,\,\, \text{is the {\it unique} integer such that} \\ 
1 \le \imath \le n \,\,\, \text{and} \,\,\, 2n+2k+2\imath \,\,\, \text{is a multiple of} \,\,\, 2n, \hspace{1.6in} \\ 
c_{2n+2k+2i}^{\, (n,k,3)} \,\,\, \text{is a period-}(2n+2k+2i) \,\,\, \text{point of} \,\,\, f \,\,\, \text{and}, \\ 
\text{for each} \,\, i \ge 1, \, \text{both} \,\, p_{m+2n+2k+2i}^{\, (n,k,3)} \,\, \text{and} \,\, q_{m+2n+2k+2i}^{\, (n,k,3)} \,\, \text{are period-}(m+2n+2k+2i) \,\, \text{points of} \,\, f.
$$
\end{multline*}

\pagebreak

\noindent
{\bf 4.3 On the interval $[\breve u_0, v]$ $(\subset [u_1, v])$}:

\noindent
The first layer, the starting step, for each $n \ge 1$,
\begin{multline*}
$$
\breve \mu_{m,n} = \min \big\{ \breve u_0 \le x \le v: f^{m+2n}(x) = d \big\}, \\ 
\breve p_{m+2n} = \min \big\{ \breve u_0 \le x \le v: f^{m+2n}(x) = x \big\}. \qquad\qquad\qquad\qquad\qquad\qquad
$$
\end{multline*}
\indent For each $n \ge 1$, $\breve p_{m+2n}$ is a period-$(m+2n)$ point of $f$.

\vspace{.1in}

\noindent
The second layer, for any fixed $n \ge 1$ and for each $k \ge 1$,
\begin{multline*}
$$
\breve \mu_{m,n,k}' = \max \big\{ \breve \mu_{m,n+1} \le x \le \breve \mu_{m,n}: f^{m+2n+2k}(x) = d \big\}, \\ 
\breve q_{m+2n+2k}^{\, (n,2)} = \max \big\{ \breve \mu_{m,n+1} \le x \le \breve \mu_{m,n}: f^{m+2n+2k}(x) = x \big\}, \qquad\qquad\qquad \\ \qquad\qquad\quad \breve c_{2m+2n+2k}'^{\, (n,2)} = \max \big\{ \breve \mu_{m,n+1} \le x \le \breve \mu_{m,n,1}': f^{2m+2n+2k}(x) = x \big\}, \\ 
\breve c_{2m+2n+2k}^{\, (n,2)} = \min \big\{ \breve \mu_{m,n+1} \le x \le \breve \mu_{m,n}: f^{2m+2n+2k}(x) = x \big\}.
$$
\end{multline*}
\indent For any fixed $n \ge 1$, we have 
\begin{multline*}
$$
\text{for each} \,\,\, k \ge 1, \, \breve q_{m+2n+2k}^{\, (n,2)} \,\,\, \text{is a period-}(m+2n+2k) \,\,\, \text{point of} \,\,\, f \,\,\, \text{and}, \\ 
\text{for each} \,\, k \ge n+3, \, \text{both} \,\, \breve c_{2m+2n+2k}'^{\, (n,2)} \,\, \text{and} \,\, \breve c_{2m+2n+2k}^{\, (n,2)} \,\, \text{are period-}(2m+2n+2k) \,\, \text{points of} \,\, f.
$$
\end{multline*}

\noindent
The third layer, for any fixed $n \ge 1$, $k \ge 1$ and for each $i \ge 1$,
\begin{multline*}
$$
\quad \breve \mu_{m,n,k,i} = \min \big\{ \breve \mu_{m,n,k}' \le x \le \breve \mu_{m,n,k+1}' : f^{m+2n+2k+2i}(x) = d \big\}, \\ 
\,\,\,\,\, \breve p_{m+2n+2k+2i}^{\, (n,k,3)} = \min \big\{ \breve \mu_{m,n,k}' \le x \le \breve \mu_{m,n,k+1}' : f^{m+2n+2k+2i}(x) = x \big\},\qquad\qquad\quad \\ 
\qquad\quad \breve c_{2m+2n+2k+2i}^{\, (n,k,3)} = \min \big\{ \breve \mu_{m,n,k,1} \le x \le \breve \mu_{m,n,k+1}' : f^{2m+2n+2k+2i}(x) = x \big\}, \\
\breve c_{2m+2n+2k+2i}'^{\, (n,k,3)} = \max \big\{ \breve \mu_{m,n,k}' \le x \le \breve \mu_{m,n,k+1}' : f^{2m+2n+2k+2i}(x) = x \big\}.
$$
\end{multline*}
\indent For any fixed $n \ge 1$ and $k \ge 1$, we have 
\begin{multline*}
$$
\text{for each} \,\,\, i \ge 1, \, \breve p_{m+2n+2k+2i}^{\, (n,k,3)} \,\,\, \text{is a period-}(m+2n+2k+2i) \,\,\, \text{point of} \,\,\, f \,\,\, \text{and}, \\ 
\text{for each} \,\,\, i \ge n+k+3, \,\qquad\qquad\qquad\qquad\qquad\qquad\qquad\qquad\qquad\qquad\qquad\qquad\qquad \\ 
\text{both} \,\,\, \breve c_{2m+2n+2k+2i}^{\, (n,k,3)} \,\,\, \text{and} \,\,\, \breve c_{2m+2n+2k+2i}'^{\, (n,k,3)} \,\,\, \text{are period-}(2m+2n+2k+2i) \,\,\, \text{points of} \,\,\, f.
$$
\end{multline*}

\pagebreak

\noindent
{\bf 4.4 On the interval $[v, \hat u_0]$ ($\subset [v, \bar u_1']$)}: 

\noindent
The first layer, the starting step, for each $n \ge 1$,
\begin{multline*}
$$
\hat \mu_{m,n}' = \max \big\{ v \le x \le \hat u_0: f^{m+2n}(x) = d \big\}, \\ 
\hat q_{m+2n} = \max \big\{ v \le x \le \hat u_0: f^{m+2n}(x) = x \big\}. \qquad\qquad\qquad\qquad\qquad\qquad\qquad\qquad
$$
\end{multline*}
\indent For each $n \ge 1$, $\hat q_{m+2n}$ is a period-$(m+2n)$ point of $f$.

\vspace{.1in}

\noindent
The second layer, for any fixed $n \ge 1$ and for each $k \ge 1$,
\begin{multline*}
$$
\hat \mu_{m,n,k} = \min \big\{ \hat \mu_{m,n}' \le x \le \hat \mu_{m,n+1}': f^{m+2n+2k}(x) = d \big\}, \\ 
\hat p_{m+2n+2k}^{\, (n,2)} = \min \big\{ \hat \mu_{m,n}' \le x \le \hat \mu_{m,n+1}': f^{m+2n+2k}(x) = x \big\}, \qquad\qquad\qquad \\ 
\qquad\qquad\quad \hat c_{2m+2n+2k}^{\, (n,2)} = \min \big\{ \hat \mu_{m,n,1} \le x \le \hat \mu_{m,n+1}': f^{2m+2n+2k}(x) = x \big\}, \\ 
\hat c_{2m+2n+2k}'^{\, (n,2)} = \max \big\{ \tilde \mu_{m,n}' \le x \le \hat \mu_{m,n+1}': f^{2m+2n+2k}(x) = x \big\}.
$$
\end{multline*}
\indent For any fixed $n \ge 1$, we have 
\begin{multline*}
$$
\text{for each} \,\,\, k \ge 1, \, \hat p_{m+2n+2k}^{\, (n,2)} \,\,\, \text{is a period-}(m+2n+2k) \,\,\, \text{point of} \,\,\, f \,\,\, \text{and}, \\ 
\text{for each} \,\, k \ge n+3, \, \text{both} \,\, \hat c_{2m+2n+2k}^{\, (n,2)} \,\, \text{and} \,\, \hat c_{2m+2n+2k}'^{\, (n,2)} \,\, \text{are period-}(2m+2n+2k) \,\, \text{points of} \,\, f.
$$
\end{multline*}

\noindent
The third layer, for any fixed $n \ge 1$, $k \ge 1$ and for each $i \ge 1$,
\begin{multline*}
$$
\hat \mu_{m,n,k,i}' = \max \big\{ \hat \mu_{m,n,k+1} \le x \le \hat \mu_{m,n,k} : f^{m+2n+2k+2i}(x) = d \big\}, \\ 
\hat q_{m+2n+2k+2i}^{\, (n,k,3)} = \max \big\{ \hat \mu_{m,n,k+1} \le x \le \hat \mu_{m,n,k} : f^{m+2n+2k+2i}(x) = x \big\}, \qquad\qquad\,\,\, \\ \,\,\,\,\quad\qquad \hat c_{2m+2n+2k+2i}'^{\, (n,k,3)} = \max \big\{ \hat \mu_{m,n,k+1} \le x \le \hat \mu_{m,n,k,1}' : f^{2m+2n+2k+2i}(x) = x \big\}, \\ 
\hat c_{2m+2n+2k+2i}^{\, (n,k,3)} = \min \big\{ \hat \mu_{m,n,k+1} \le x \le \hat \mu_{m,n,k} : f^{2m+2n+2k+2i}(x) = x \big\}.
$$
\end{multline*}
\indent For any fixed $n \ge 1$ and $k \ge 1$, we have 
\begin{multline*}
$$
\text{for each} \,\,\, i \ge 1, \, \hat q_{m+2n+2k+2i}^{\, (n,k,3)} \,\,\, \text{is a period-}(m+2n+2k+2i) \,\,\, \text{point of} \,\,\, f \,\,\, \text{and}, \\ 
\text{for each} \,\,\, i \ge n+k+3, \,\qquad\qquad\qquad\qquad\qquad\qquad\qquad\qquad\qquad\qquad\qquad\qquad\qquad \\ 
\text{both} \,\,\, \hat c_{2m+2n+2k+2i}'^{\, (n,k,3)} \,\,\, \text{and} \,\,\, \hat c_{2m+2n+2k+2i}^{\, (n,k,3)} \,\,\, \text{are period-}(2m+2n+2k+2i) \,\,\, \text{points of} \,\,\, f.
$$
\end{multline*}

\pagebreak

\noindent
{\bf 4.5 On the interval $[\bar u_1', z_0]$}:

\noindent
The first layer, the starting step, for each $n \ge 1$ (note that there is a {\it twist} in the definition of $\bar c_{2n+2}$),
\begin{multline*}
$$
\bar u_n' = \max \big\{ v \le x \le z_0: f^{2n}(x) = d \big\}, \\ 
\bar c_{2n+2} = \min \big\{ \bar u_n' \le x \le z_0: f^{2n+2}(x) = x \big\}.\quad\qquad\qquad\qquad\qquad\qquad\qquad\qquad\qquad
$$
\end{multline*}
\indent For each $n \ge 1$, $\bar c_{2n+2}$ is a period-$(2n+2)$ point of $f$.

\vspace{.1in}

\noindent
The second layer, for any fixed $n \ge 1$ and for each $k \ge 1$,
\begin{multline*}
$$
\bar u_{n,k} = \min \big\{ \bar u_n' \le x \le \bar u_{n+1}': f^{2n+2k}(x) = d \big\}, \\ 
\bar c_{2n+2k}^{\, (n,2)} = \min \big\{ \bar u_n' \le x \le \bar u_{n+1}': f^{2n+2k}(x) = x \big\},\qquad\qquad\qquad\qquad \,\,\,\,\,\,\,\,\,\,\,\, \\ 
\,\,\,\qquad\qquad \bar p_{m+2n+2k}^{\, (n,2)} = \min \big\{ \bar u_{n,1} \le x \le \bar u_{n+1}': f^{m+2n+2k}(x) = x \big\}, \\ 
\bar q_{m+2n+2k}^{\, (n,2)} = \max \big\{ \bar u_n' \le x \le \bar u_{n+1}': f^{m+2n+2k}(x) = x \big\}.
$$
\end{multline*}
\indent For any fixed $n \ge 1$, we have 
\begin{multline*}
$$
\text{for each} \,\,\, k \ge 1, \, \bar c_{2n+2k}^{\, (n,2)} \,\,\, \text{is a period-}(2n+2k) \,\,\, \text{point of} \,\,\, f \,\,\, \text{and}, \\ 
\text{for each} \,\,\, k \ge 1, \, \text{both} \,\,\, p_{m+2n+2k}^{\, (n,2)} \,\,\, \text{and} \,\,\, q_{m+2n+2k}^{\, (n,2)} \,\,\, \text{are period-}(m+2n+2k) \,\,\, \text{points of} \,\,\, f.
$$
\end{multline*}

\noindent
The third layer, for any fixed $n \ge 1$, $k \ge 1$ and for each $i \ge 1$,
\begin{multline*}
$$
\bar u_{n,k,i}' = \max \big\{ \bar u_{n,k+1} \le x \le \bar u_{n,k} : f^{2n+2k+2i}(x) = d \big\}, \\ 
\bar c_{2n+2k+2i}'^{\, (n,k,3)} = \max \big\{ \bar u_{n,k+1} \le x \le \bar u_{n,k} : f^{2n+2k+2i}(x) = x \big\},\qquad\qquad\qquad\quad \,\,\,\,\,\,\, \\ 
\quad\qquad\qquad\qquad\quad\, \bar q_{m+2n+2k+2i}^{\, (n,k,3)} = \max \big\{ \bar u_{n,k+1} \le x \le \bar u_{n,k,1}' : f^{m+2n+2k+2i}(x) = x \big\}, \qquad\qquad \\ 
\bar p_{m+2n+2k+2i}^{\, (n,k,3)} = \min \big\{ \bar u_{n,k+1} \le x \le \bar u_{n,k} : f^{m+2n+2k+2i}(x) = x \big\}.
$$
\end{multline*}
\indent For any fixed $n \ge 1$ and $k \ge 1$, we have 
\begin{multline*}
$$
\text{for each} \,\,\, i \ge n+k+2, \, \bar c_{2n+2k+2i}^{\, (n,k,3)} \,\,\, \text{is a period-}(2n+2k+2i) \,\,\, \text{point of} \,\,\, f \,\,\, \text{and}, \\ 
\text{for each} \,\, i \ge 1, \, \text{both} \,\, \bar p_{m+2n+2k+2i}^{\, (n,k,3)} \,\, \text{and} \,\, \bar q_{m+2n+2k+2i}^{\, (n,k,3)} \,\, \text{are period-}(m+2n+2k+2i) \,\, \text{points of} \,\, f.
$$
\end{multline*}
$$\aleph \qquad\qquad\qquad \aleph \qquad\qquad\qquad \aleph \qquad\qquad\qquad \aleph \qquad\qquad\qquad \aleph$$
\indent Now a pattern emerges which allows us to write down the formulas for the periodic points of the fourth and higher layers of the {\it basic} tower of periodic points of $f$ associated with $P$ and then we can apply Lemmas 5 $\&$ 6 to determine the least periods of all except possibly finitely many of these periodic points.  We leave the details to the interested readers.

\pagebreak

{\bf $\mathsection$5. An example.}

Let $P$ denote a period-3 orbit of $f$ and let $m = 3$.  We now describe the first and the second layers of the {\it basic} tower of periodic points of $f$ associated with $P$ and leave the third and the higher layers of the {\it basic} tower to the interested readers.  We shall follow the notations introduced before.

On the first layer,
\begin{multline*}
$$
\min P \le \tilde q_3 < \tilde q_5 < \tilde q_7 < \cdots < d < \cdots < c_6 < c_4 < c_2 < u_1 \le \breve u_0 < \cdots < \breve p_9 < \breve p_7 < \breve p_5 < v < \\ 
\hat q_5 < \hat q_7 < \hat q_9 < \cdots < \hat u_0 \le \bar u_1' < \bar c_4 < \bar c_6 < \bar c_8 < \cdots < z_0;
$$
\end{multline*}
\indent On the second layer, (we deliberately leave out the points \, $\tilde p_5^{(0,2)}, \, \tilde p_7^{(1,2)}, \, \tilde p_9^{(2,2)}, \, \tilde p_{11}^{(3,2)}, \cdots)$ 
$$
\,\, (\min P) <\qquad\qquad\qquad\qquad\qquad\qquad
$$
{\footnotesize $\tilde \mu_{m,0}' < \cdots < \tilde p_{11}^{(0,2)} < \tilde p_9^{(0,2)} < \tilde p_7^{(0,2)} < \cdots < \tilde c_{12}^{(0,2)} < \tilde c_{10}^{(0,2)} < \tilde c_8^{(0,2)} < \tilde c_8'^{(0,2)} < \tilde c_{10}'^{(0,2)} < \tilde c_{12}'^{(0,2)} < \cdots < $ \\
\indent $\tilde \mu_{m,1}' < \cdots < \tilde p_{13}^{(1,2)} < \tilde p_{11}^{(1,2)} < \tilde p_9^{(1,2)} < \cdots < \tilde c_{14}^{(1,2)} < \tilde c_{12}^{(1,2)} < \tilde c_{10}^{(1,2)} < \tilde c_{10}'^{(1,2)} < \tilde c_{12}'^{(1,2)} < \tilde  c_{14}'^{(1,2)} < \cdots < $ \\
\indent\indent $\tilde \mu_{m,2}' < \cdots < \tilde p_{15}^{(2,2)} < \tilde p_{13}^{(2,2)} < \tilde p_{11}^{(2,2)} < \cdots < \tilde c_{16}^{(2,2)} < \tilde c_{14}^{(2,2)} < \tilde c_{12}^{(2,2)} < \tilde c_{12}'^{(2,2)} < \tilde c_{14}'^{(2,2)} < \tilde  c_{16}'^{(2,2)} < \cdots < $}  
$$
\,\,\,\, d < \qquad\qquad\qquad\qquad\qquad\qquad
$$ 
{\footnotesize $u_4 < \cdots < p_{15}^{(3,2)} < p_{13}^{(3,2)} < p_{11}^{(3,2)} < q_{11}^{(3,2)} < q_{13}^{(3,2)} < q_{15}^{(3,2)} < \cdots < c_{8}'^{(3,2)} < c_{10}'^{(3,2)} < c_{12}'^{(3,2)} < \cdots <$ \\
\indent $u_3 < \cdots < p_{13}^{(2,2)} < p_{11}^{(2,2)} < p_{9}^{(2,2)} < q_{9}^{(2,2)} < q_{11}^{(2,2)} < q_{13}^{(2,2)} < \cdots < c_{6}'^{(2,2)} < c_{8}'^{(2,2)} < c_{10}'^{(2,2)} < \cdots <$ \\
\indent\indent $u_2 < \cdots < p_{11}^{(1,2)} < p_{9}^{(1,2)} < p_{7}^{(1,2)} < q_{7}^{(1,2)} < q_{9}^{(1,2)} < q_{11}^{(1,2)} < \cdots < c_{4}'^{(1,2)} < c_{6}'^{(1,2)} < c_{8}'^{(1,2)} < \cdots <$} 
$$
\,\, u_1 \le \breve u_0 < \qquad\qquad\qquad\qquad\qquad\quad
$$ 
{\footnotesize $\breve \mu_{m,4} < \cdots < \breve c_{18}^{(3,2)} < \breve c_{16}^{(3,2)} < \breve c_{14}^{(3,2)} < \breve c_{14}'^{(3,2)} < \breve c_{16}'^{(3,2)} < \breve c_{18}'^{(3,2)} < \cdots < \breve q_{11}^{(3,2)} < \breve q_{13}^{(3,2)} < \breve q_{15}^{(3,2)} < \cdots <$ \\
\indent $\breve \mu_{m,3} < \cdots < \breve c_{16}^{(2,2)} < \breve c_{14}^{(2,2)} < \breve c_{12}^{(2,2)} < \breve c_{12}'^{(2,2)} < \breve c_{14}'^{(2,2)} < \breve c_{16}'^{(2,2)} < \cdots < \breve q_{9}^{(2,2)} < \breve q_{11}^{(2,2)} < \breve q_{13}^{(2,2)} < \cdots <$ \\
\indent\indent $\breve \mu_{m,2} < \cdots < \breve c_{14}^{(1,2)} < \breve c_{12}^{(1,2)} < \breve c_{10}^{(1,2)} < \breve c_{10}'^{(1,2)} < \breve c_{12}'^{(1,2)} < \breve c_{14}'^{(1,2)} < \cdots < \breve q_{7}^{(1,2)} < \breve q_{9}^{(1,2)} < \breve q_{11}^{(1,2)} < \cdots <$} 
$$
\breve \mu_{m,1} < v < \qquad\qquad\qquad\qquad\qquad
$$ 
{\footnotesize $\hat \mu_{m,1}' < \cdots < \hat p_{11}^{(1,2)} < \hat p_9^{(1,2)} < \hat p_7^{(1,2)} < \cdots < \hat c_{14}^{(1,2)} < \hat c_{12}^{(1,2)} < \hat c_{10}^{(1,2)} < \hat c_{10}'^{(1,2)} < \hat c_{12}'^{(1,2)} < \hat c_{14}'^{(1,2)} < \cdots < $ \\
\indent $\hat \mu_{m,2}' < \cdots < \hat p_{13}^{(2,2)} < \hat p_{11}^{(2,2)} < \hat p_9^{(2,2)} < \cdots < \hat c_{16}^{(2,2)} < \hat c_{14}^{(2,2)} < \hat c_{12}^{(2,2)} < \hat c_{12}'^{(2,2)} < \hat c_{14}'^{(2,2)} < \hat  c_{16}'^{(2,2)} < \cdots < $ \\
\indent\indent $\hat \mu_{m,3}' < \cdots < \hat p_{15}^{(3,2)} < \hat p_{13}^{(3,2)} < \hat p_{11}^{(3,2)} < \cdots < \hat c_{18}^{(3,2)} < \hat c_{16}^{(3,2)} < \hat c_{14}^{(3,2)} < \hat c_{14}'^{(3,2)} < \hat c_{16}'^{(3,2)} < \hat  c_{18}'^{(3,2)} < \cdots < $}
$$
\hat u_0 \le \qquad\qquad\qquad\qquad\quad\quad\,\,\,
$$ 
{\footnotesize $\bar u_1' < \cdots < \bar c_{8}^{(1,2)} < \bar c_6^{(1,2)} < \bar c_4^{(1,2)} < \cdots < \bar p_{11}^{(1,2)} < \bar p_{9}^{(1,2)} < \bar p_{7}^{(1,2)} < \bar q_{7}^{(1,2)} < \bar q_{9}^{(1,2)} < \bar q_{11}^{(1,2)} < \cdots < $ \\
\indent $\bar u_2' < \cdots < \bar c_{10}^{(2,2)} < \bar c_{8}^{(2,2)} < \bar c_6^{(2,2)} < \cdots < \bar p_{13}^{(2,2)} < \bar p_{11}^{(2,2)} < \bar p_{9}^{(2,2)} < \bar q_{9}^{(2,2)} < \bar q_{11}^{(2,2)} < \bar  q_{13}^{(2,2)} < \cdots < $ \\
\indent\indent $\bar u_3' < \cdots < \bar c_{12}^{(3,2)} < \bar c_{10}^{(3,2)} < \bar c_{8}^{(3,2)} < \cdots < \bar p_{15}^{(3,2)} < \bar p_{13}^{(3,2)} < \bar p_{11}^{(3,2)} < \bar q_{11}^{(3,2)} < \bar q_{13}^{(3,2)} < \bar  q_{15}^{(3,2)} < \cdots < $}
$$
(z_0). \qquad\qquad\qquad\qquad\qquad\,\,\,
$$
\end{document}